\documentclass{amsart}
\usepackage[left=2.5cm,right=2cm,
    top=2.5cm,bottom=2cm,bindingoffset=0cm]{geometry}
\usepackage{amsmath, amsthm, amssymb, amstext}
\usepackage{hyperref}
\usepackage{array,amsfonts,mathrsfs,bm}
\usepackage{cancel}
\newtheorem{theorem}{Theorem}[section]
\newtheorem{lemma}[theorem]{Lemma}
\newtheorem{prop}[theorem]{Proposition}
\theoremstyle{definition}
\newtheorem{definition}[theorem]{Definition}

\theoremstyle{remark}
\newtheorem{remark}[theorem]{Remark}
\newtheorem*{nota}{Notation}
\usepackage{graphicx,color}
\usepackage{oubraces}
\usepackage[shortlabels]{enumitem}

\numberwithin{equation}{section}

\begin{document}

\title[Determinants of spinor Laplacians on translation surfaces]{Laplacians in spinor bundles over translation surfaces: self-adjoint extentions and regularized determinants}

\author{Alexey Kokotov}
\address{Department of Mathematics \& Statistics, Concordia University, 1455 De Maisonneuve Blvd. W. Montreal, QC  H3G 1M8, Canada}
\email{alexey.kokotov@concordia.ca}
\thanks{The research of the first author was supported by NSERC}

\author{Dmitrii Korikov}
\address{Department of Mathematics \& Statistics, Concordia University, 1455 De Maisonneuve Blvd. W. Montreal, QC  H3G 1M8, Canada}
\email{dmitrii.v.korikov@gmail.com}
\thanks{The research of the second author was supported by Fonds de recherche du Qu\'ebec.}

\date{\today}

%58J52 Global analysis, analysis on manifolds->Determinants and determinant bundles, analytic torsion
%35P99 	Partial differential equations->None of the above, but in this section
%32G15 Several complex variables and analytic spaces->Moduli of Riemann surfaces, Teichmüller theory
%32G08 Several complex variables and analytic spaces->Deformations of fiber bundles
%30F10 Functions of a complex variable->Compact Riemann surfaces and uniformization
%30F45 Functions of a complex variable->Conformal metrics (hyperbolic, Poincaré, distance functions)

\begin{abstract}
 We study the regularized determinants ${\rm det}\, \Delta$  of various self-adjoint extensions of symmetric Laplacians acting in spinor bundles over compact Riemann surfaces with flat singular metrics $|\omega|^2$, where $\omega$ is a holomorphic one form on the Riemann surface.
 We find an explicit expression for ${\rm det}\, \Delta$ for the so-called self-adjoint Szeg\"o extension through the Bergman tau-function on the moduli space of Abelian differentials and the theta-constants (corresponding to the spinor bundle). This expression can be considered as a version of the well-known spin-$1/2$ bosonization formula of  Bost-Nelson for the case of flat conformal metrics with conical singularities and a higher genus generalization of the Ray-Singer formula for flat elliptic curves. We establish comparison formulas for the determinants of two different extensions (e. g., the Szeg\"o extension and the Friedrichs one). The paper answers a question raised by D'Hoker and Phong \cite{DH-P} more than thirty years ago. We also reconsider the results from \cite{DH-P} on the regularization of diverging determinant ratio for Mandelstam metrics (for any spin) proposing (and computing) a new regularization of this ratio.      
\end{abstract}

\maketitle

\section{Introduction}
Let $\omega$ be a holomorphic one-form with $2g-2$ simple zeros $P_1,\dots,P_{2g-2}$ on a compact Riemann surface $X$ of genus $g>1$. Using the one-form $\omega$ one can construct 
\begin{itemize} 
\item a hermitian metric $h$ in the $n$-th power ($n=\pm 1, \pm 2, \dots$), $K^n$, of the canonical bundle $K$ over
$X$:
$$h=\omega^{-n}\bar\omega^{-n}\,$$
and, in particular,
 a conformal Riemannian metric $$|\omega|^2:=\rho^{-2}(z, \bar z)|dz|^2$$ on $X$ (or, equivalently,  
a hermitian metric $h=\omega \bar\omega$ in the holomorphic line bundle $K^{-1}$).

\item a hermitian metric  $h$ in the $(n+1/2)$-spinor bundle 
$C\otimes K^n$, where $C$ is one of the $4^g$ holomorphic line bundles with property $C\otimes C=K$:
$$h=|\omega|^{-1}\omega^{-n}\bar\omega^{-n}\,.$$ 

\end{itemize} 
All these metrics have singularities at the zeroes of $\omega$. 

 For $L=K^n$ with $n\geq 2$ (or $L=K^n\otimes C$ with $n\geq 0$) introduce the Dolbeault Laplacian
\begin{equation}
\label{Dol}
\Delta_L=-\rho^2h^{-1}\partial(h\bar\partial)
\end{equation}
%old
%=-|\omega|^{-2}h^{-1}\partial(h\bar \partial)\,.}
%\begin{equation}\label{Dol}\Delta_L=4\rho^2h^{-1}\partial(h\bar\partial)=4|\omega|^{-2}h^{-1}\partial(h\bar \partial)\,.\end{equation}

\begin{remark}
Notice that for $m=n$ {\rm(}or $m=n+1/2${\rm)} the operator $\Delta^{-}_m$ from {\rm(2.3)} of {\rm\cite{DH-P}} coincides
 with $\Delta_L$ up to a numerical factor. 
\end{remark}

In \cite{DH-P} the authors, being motivated by the need to compute various partition functions appearing in the string theory, stated the problem of finding a reasonable  regularization
of the formal expression
\begin{equation}\label {QUANT} \frac{{\rm det}' \Delta_L}{{\rm det}\, <e_i, e_j>_L\,{\rm det}\,<e^*_k, e_l^*>_{L^{-1}\otimes K}}\,, \end{equation} 
where ${\rm det}'$ stands for (modified, i. e. with excluded zero mode) determinant and
$\{e_j\},\ \{e^*_k\}$ are the bases of the spaces of holomorphic sections of $L$ and $L^{-1}\otimes K$ respectfully.

(Let us remind the reader that for smooth metrics the operators $\Delta^{-}_m$ and $\Delta^{-}_{1-m}$ have the same (non-zero) eigenvalues (\cite{DH-P}) and, therefore, in this case \ref{QUANT} is invariant under the change $m\mapsto 1-m$; it is reasonable to postulate the same invariance for the regularization in question, hence the conditions $n\geq 2$ or $n\geq 0$ above. The scalar case (i. e. $\Delta^{-}_0$ and, simultaneously, $\Delta^{-}_{1-0}=\Delta^{-}_1$) is excluded as well-known
(see \cite{KokKorot}).)

 Notice that the quantities entering "the Quillen metric"  (\ref{QUANT}) are meaningless due to the following  reasons: 
 \begin{itemize}
\item it is unclear how to naturally associate to singular differential expression (\ref{Dol}) an operator in Hilbert space, so  neither of definitions of the determinant through the spectrum is immediately applicable. 

\item while the second Gram determinant in the denominator is absent in all the cases except  $L=C$ (bundles of negative degree do not have holomorphic sections), the first Gram determinant is infinite (again, in all the cases except $L=C$).   
 \end{itemize} 

In \cite{DH-P} a regularization of (\ref{QUANT}) was proposed (see also (\cite{Sonoda}) for a similar result and the discussion in the Appendix to \cite{HKKK}), and question of the spectral theoretic interpretation of  (\ref{QUANT}) in case $L=C$ was raised. 

The main idea of D'Hoker-Phong's regularization (as well as Sonoda's)
is to apply Polyakov's formula (that compares the determinants of Laplacians in two conformally equivalent smooth metrics) to the pair (singular metric, Arakelov metric) and then to make use of the known results on the Arakelov determinant. Of course, in case of singular metric the Polyakov's formula can be applied only formally (logically, this application constitutes a part of definition), moreover,
it now contains diverging integral which should be cleverly regularized. 

The present paper consists of two parts. In the short first part (Section 2) we 
  propose a new regularization of (\ref{QUANT}) and compute it explicitly (as in \cite{DH-P} no spectral properties of (\ref{Dol}) are needed for that). It turns out that this part is relatively easy: the key idea of the regularization (going back to \cite{EKZ})  was already used in a similar context in \cite{Kok-MRL}, all the results lying in the background  were elaborated in \cite{HK}, \cite{HK1}, \cite{KL} and \cite{KokKorot}. 
 In the much longer second part (Sections 3, 4, 5, 6) we address the D'Hoker-Phong question concerning the case $K=C$. The main results of the second part are described in Section 3. 

\

\paragraph{\it Acknowledgements.} The first author acknowledges numerous important discussions with Luc Hillairet during our stay at Max Planck Institute in Bonn in 2018, the first preliminary plan of the present work appeared as a result of these discussions. Several important ideas (in particular, the introduction of the holomorphic extension and the observation on the almost isospectrality  of the holomorphic and the Friedrichs extensions) used below go back to L. Hillairet's papers on the scalar conical Laplacians. 
 
\section{Regularization of Quillen metric}  
\label{sec reqularization}
 
In this section we closely follow $\S$ 3.2 from \cite{Kok-MRL}, where the case of $L=K^2$ and the metric given by $|Q|$, where $Q$ is a holomorphic quadratic differential with simple zeros was considered. The method of regularization for (\ref{QUANT}) used in \cite{Kok-MRL} was taken from \cite{EKZ} (where it was applied to scalar conical Laplacians) and remains the same in our case.  Explicit computation of the resulting quantity is based on the same ingredients that were used in \cite{Kok-MRL}: the reference to Theorem 5.8 from \cite{Fay} and the statement on the essential coincidence of the two regularizations of the scalar conical Laplacians (cf. Proposition 6 from \cite{Kok-MRL}).

Let $L=K^n$ (with $n\geq 2$) or $L=K^n\otimes C$ (with $n\geq 1$). Following \cite{KZ}, introduce the moduli space $H_g(1, \dots, 1)$ of pairs $(X,\omega)$, where $X$ is a compact Riemann surface of genus $g$ and $\omega$ is a holomorphic differential on $X$ with $2g-2$ simple zeroes. 
In \cite{DH-P} D'Hoker and Phong put into correspondence to any pair $(X, \omega)\in H_g(1, \dots, 1)$ 
(with a canonical basis of cycles on $X$ chosen)
a basis of holomorphic sections $\{e_1, \dots e_N\}, N={\rm dim}\,H^0(L)$ of the line bundle $L$ over $X$; the latter basis can be constructed explicitly (through $\omega$, theta-functions, prime-form, etc).
Say, for the simplest case $K^2$ the basis of $N=3g-3$ holomorphic quadratic differentials on $X$ corresponding
to a pair $(X, \omega)\in H_g(1, \dots, 1)$ is given by
$$\omega v_k; \ \ k= 1, 2, \dots, g; \ \ \ \ \omega \Omega_{P_{2g-2}, P_k}, \ \ k=1, \dots, 2g-3\,,$$
where $\{v_k\}_{k=1}^g$ is a basis of normalized holomorphic differentials, $P_1, \dots, P_{2g-2}$ are the zeros of $\omega$ and $\Omega_{a-b}$ is the Abelian  differential of the third kind
with zero $a$-periods and simple poles at $a, b$ with residues $1$ and $-1$.
We refer the reader to p. 641 of \cite{DH-P} for the construction of the basis in the general case $K^n$ or $K^n\otimes C$.  We will call this basis {\it the D'Hoker-Phong basis}.

In a vicinity of a simple zero $P_k$ of a holomorphic differential $\omega$ the metric $|\omega|^2$ coincides with metric of the standard round cone of the angle $4\pi$
$$|\omega|^2=|x_k|^2||dx_k|^2\,$$
in the so-called distinguished holomorphic coordinate
\begin{equation}
\label{dist coord}
x_k(Q)=\left(2 \int_{P_k}^Q\omega\right)^{1/2}\,.
\end{equation}
Smoothing all these standard cones in $\epsilon$-vicinities of their tips (the zeros of differential $\omega$) one gets a smooth conformal metric $\rho_{(\epsilon)}^{-2}|dz|^2=|\omega|^2_{ (\epsilon)}$ on $X$.  One can perform this smoothing without changing the area, thus, we assume that $\int_X|\omega|^2=\int_X|\omega|^2_{ (\epsilon)}$.  Simultaneously, one gets a smooth hermitian metric $h_{(\epsilon)}$ in $L$. 
   
Denote by $\mathcal{G}(\{e_k\}; X, |\omega|^2_{(\epsilon)}, h_{(\epsilon)})$ the Gram determinant of the hermitian products of the elements of the D'Hoker-Phong basis w. r. t. the metrics $(|\omega|^2_{(\epsilon)}, h_{(\epsilon)})$. 

Notice, that, since the metrics $|\omega|^2_{(\epsilon)}$ and  $h_{(\epsilon)}$ are smooth, one can introduce the $\zeta$-regularized modified (i. e. with zero modes excluded) determinant of the Dolbeault Laplacian ${\rm det'} \Delta_{|\omega|^2_{(\epsilon)}}(X, \omega)$ in $L$ corresponding to these metrics.  Moreover, the Gram determinant for basic holomorphic sections $e_k$ of the bundle $L$ equipped  with metric $h_{(\epsilon)}$ is finite.

Choose a compact Riemann surface $X_0$ of genus $g$ and a holomorphic differential $\omega_0$ with $2g-2$ simple zeros on $X_0$ (in other words we choose an element $(X_0, \omega_0)$ of the the moduli space $H_g(1, \dots, 1)$). %of holomorphic differentials with $2g-2$ simple zeros on compact Riemann surfaces of genus $g$. 

Now for any $(X, \omega)\in H_g(1, \dots, 1)$ introduce the ratio 
\begin{equation*}
%\label{Ratio} 
R(X, \omega; \epsilon)=\frac{ {\rm det'} \Delta_{|\omega|^2_{(\epsilon)}}(X, \omega)\left(   \mathcal{G}(\{e_k\}; X, |\omega|^2_{(\epsilon)}, h_{(\epsilon)}) \right)^{-1}             }{{\rm det'} \Delta_{|\omega_0|^2_{(\epsilon)}}(X_0, \omega_0)\left(   \mathcal{G}(\{e_k\}; X_0, |\omega_0|^2_{(\epsilon)}, h^0_{(\epsilon)}) \right)^{-1} }\,.
\end{equation*}
The following key observation (in case of scalar conical Laplacians) was first made in \cite{EKZ}.
The proof repeats the proof of Lemma 2 from \cite{Kok-MRL} verbatim. Essentially, it reduces to a straightforward application of Proposition 3.8 from \cite{Fay} (a version of Polyakov's formula for Dolbeault Laplacians which is in fact a specification  of the general   Bismut-Gillet-Soul\'e  anomaly formula \cite{BGS}). 
  
\begin{lemma} For sufficiently small $\epsilon$ the quantity $R(X, \omega; \epsilon)$ is $\epsilon$-independent. 
\end{lemma}
\begin{definition}
 Let $(X, \omega)$ be an element of $H_g(1, \dots, 1)$. Let $L=K^n$ with $n\geq 2$ or $K=K^n\otimes C$ with $n\geq 1$ and let the metrics on $X$ and in $L$ are constructed via the holomorphic differntial $\omega$. 
Assign to  the formal expression from {\rm(\ref{QUANT})} the following value:
\begin{equation}\label{VALUE}\frac{{\rm det}' \Delta_L}{{\rm det}\, <e_i, e_j>_L}:=
\lim_{\epsilon\to 0}R(X, \omega; \epsilon)\end{equation}
\end{definition}
Thus, quantity (\ref{QUANT}) is correctly defined up to moduli (i. e. coordinates in the space $H_g(1, \dots, 1)$) independent overall constant (ruled by the choice of basic element $(X_0, \omega_0)$); it is also worth noticing that \ref{QUANT} provides a well-defined
metric on the determinant line bundle.  

It turns out that one can explicitly compute this quantity. To this end let us specify the choice of the spinor bundle $C$.  
Let ${\bf \Delta}$ be the line bundle of degree $g-1$ (equivalently, the linear equivalence class of divisors on X)
with $\mathcal{A}_{P_0}({\bf \Delta})=-K_{P_0}$. Here $P_0$ is a basic point, $\mathcal{A}_{P_0}$ is the Abel map, $K_{P_0}$ is the vector of Riemann constants.
Let $p, q\in \{0, \frac{1}{2}\}^g$ and let $\chi^{p, q}$ be the line bundle of degree zero (simultaneously the linear equivalence class of divisors) with $\mathcal{A}_{P_0}(\chi^{p, q})={\mathbb B}p+q$, where ${\mathbb B}$ is the matrix of $b$-periods of $X$. Let 
$$C={\bf \Delta}\otimes \chi^{p, q}\,.$$
Also, let ${\rm det'}\Delta_0(|\omega|^2_{(\epsilon)}, X)$ be the $\zeta$-regularized determinant of the scalar Laplacian in the smooth metric $|\omega^2|^2_{(\epsilon)}$
and let
$$c_0(X, |\omega|^2_{(\epsilon)})^{-2}=\frac{{\rm det'}\Delta_0(|\omega|^2_{(\epsilon)}, X)}
{{\rm det\,}\Im \mathbb B\,{\rm Area\,}(X, |\omega|^2_{(\epsilon)})}$$
(we remind the reader that ${\rm Area\,}(X, |\omega|^2_{(\epsilon)})={\rm Area\,}(X, |\omega|^2)$).

Let $m=2n+1$ if $L=K^n\otimes C$ and $m=2n$ if $L=K^n$. Then 
\begin{equation}
\label{schjot}
\begin{split}
R(X, \omega, \epsilon&)=c_0^{-(3m^2-6m+2)}(X, |\omega|^2_{(\epsilon)})\,\times \\
\times &\frac{\left[ {\rm det'} \Delta_{|\omega|^2_{(\epsilon)}}(X, \omega)\left(   \mathcal{G}(\{e_k\}; X, |\omega|^2_{(\epsilon)}, h_{(\epsilon)}) \right)^{-1}c_0^{(3m^2-6m+2)}(X, |\omega|^2_{(\epsilon)})\right]}
{{\rm det'} \Delta_{|\omega_0|^2_{(\epsilon)}}(X_0, \omega_0)\left(   \mathcal{G}(\{e_k\}; X_0, |\omega_0|^2_{(\epsilon)}, h^0_{(\epsilon)}) \right)^{-1} } 
.
\end{split}
\end{equation}
According to Theorem 5.8 from \cite{Fay}, the expression in the square brackets in (\ref{schjot}), coincides (up to a moduli independent constant) for $L=C\otimes K^n$, $m=2n+1$ with
\begin{equation*}
\left|\frac{\theta\left[^p _q\right]\left((m-1)K_{P_0}+\sum_{i=1}^{N}\mathcal{A}_{P_0}(x_1+\dots+x_N)\right)\prod_{i<j}^N E(x_i, x_j)}
{ {\rm det}|| e_i(x_j)||\prod_{k=1}^{N}{\bf c}(x_i)^{\frac{N-1}{g-1}}}\right|^2,
\end{equation*}
where $x_1, \dots, x_{N}$ are arbitrary points of $X$, $\theta\left[^p _q\right]$ is theta function defined in (1.8) in \cite{Fay}, $E(x, y)$ is the prime form, and ${\bf c}(\cdot)$ is defined in (1.17) in \cite{Fay}. In case $L=K^n$, $m=2n$ the expression in the square brackets in (\ref{schjot}) coincides with $|{\frak M}_{g, n}|^2$, where
 $${\frak M}_{g, n}=\frac{\theta((2n-1)K_{P_0})} {W[e_1, e_2, \dots, e_N](P_0){\bf c}(P_0)^{(2n-1)^2}}$$
 with $N=(2n-1)(g-1)$ and $W$ being the Wronskian for the D'Hoker-Phong basis 
 (cf. Theorems 5.8 and 1. 6 from \cite{Fay}).  
   
 On the other hand, according to Proposition 6 from \cite{Kok-MRL} (reformulated for the case of metric $|\omega|^2$ with conical angles $4\pi$; the proof requires no changes), the determinant $\allowbreak {\rm det'}\Delta_0(|\omega|^2_{(\epsilon)}, X)$ coincides up to moduli independent constant with
 the  determinant  ${\rm det'}\Delta_0(|\omega|^2, X)$ of the Friedrichs extension of the symmetric scalar Laplacian in the conical metric $|\omega|^2$, and, therefore, due to the results of \cite{KokKorot} one has the relation
 $$c_0(X, |\omega|^2_{(\epsilon)})^{-2}={\rm const\,}|\tau_B(X, \omega)|^2\,,$$
 where $\tau_B$ is the Bergman tau-function on the space $H_g(1, \dots, 1)$ (see (3.1), \cite{KokKorot} for an explicit expression for $\tau_B$).  
 Since the expression in the denominator of (\ref{schjot}) is moduli independent,
 all this leads to an explicit expression for (\ref{VALUE}) through $\omega$ and holomorphic invariants of the surface $X$ up to moduli independent constant. 
 
\subsection*{Case of the spinor bundle $C={\bf \Delta}\otimes\chi^{p, q}$}
The regularization (\ref{VALUE}) with no changes is applicable in the case of the spinor bundle $L=C$. However, in this case Theorem 5.8 from \cite{Fay} is no longer relevant, and to compute the regularized quantity 
(\ref{VALUE}) one instead uses  the bosonization formula of Bost and Nelson (\cite{BN}) (cf. \cite{Fay}, Theorem 4.9). For simplicity, we restrict ourselves to the case of even characteristics (i. e. even $4<p, q>$). In this case  $N=h^0(\Delta\otimes\chi^{p, q})=0$ generically.

According to the Bost-Nelson result one has
$${\det}'\Delta(|\omega|^2_{(\epsilon)}, X)={\rm const}\,\left(\frac{{\rm det'}\Delta_0(|\omega|^2_{(\epsilon)}, X)}{{\rm Area\,}(X, |\omega|^2_{(\epsilon)}) {\rm det}\,\Im {\mathbb B}}         \right)^{-1/2}|\theta[ ^p _q](0)|^2$$ 
and formula (\ref{schjot}) (with no Gram determinants!) together with Proposition 6 from \cite{Kok-MRL} and results from \cite{KokKorot}  give
\begin{equation}\label{result}R(X, \omega, \epsilon)=\lim_{\epsilon\to 0} R(X, \omega, \epsilon)={\rm const}|\theta\left[ ^p _q\right](0)|^2|\tau_B(X, \omega)|^{-1}\end{equation} 

In the remaining part of the paper we clarify to what extent the regularization (\ref{result}) of the spinor conical determinant is compatible with possible spectral regularizations of the determinant of the Laplacian acting in ${\bf \Delta}\otimes\chi^{p, q}$ with
conical metric.

\section{Expressions for determinants}
As emphasized after (\ref{QUANT}), the problem arising with the spectral regularization  of Laplacians on manifolds with conical points is that there are many different ways to associate to differential expression (\ref{Dol}) a self-adjoint operator.

\begin{nota}
\nonumber
We will use the following notation:
\begin{itemize}
\item $L_2(X;L)$ is the space of square integrable sections of $L$ with scalar product
\begin{equation}
\label{L2 scalar}
(a,b)_{L_2(X;L)}:=\int_X ah\overline{b}dS,
\end{equation}
where $h$ is a metric on $L$ and $dS:=\rho^{-2}d\overline{z}dz/2i$ is the area form.

\item $\dot{X}$ is the surface obtained by deleting all conical points from $X$ and $C_c^\infty(\dot{X};L)$ is the set of smooth sections of $L$ with supports contained in $\dot{X}$.

\item $c_{\dots}^{\dots}(\dots)$ are constants depending only on the indices and the parameters in brackets.

\item $z\mapsto v(z)$ as a rule means the section $v$ of a line bundle $L$ {\rm(}which can be identified from the context{\rm)} over $X$ written in the local parameter $z$. Also, we denote $v(P_k):=v(x_k)|_{x_k=0}$.
\end{itemize}
\end{nota}

In what follows, we deal with the spinor bundle $L=C$ and denote by $\Delta$ the Laplacian $\Delta_C$ given by (\ref{Dol}). We start with the minimal closed extension $\Delta_{min}$ which is obtained by closure in $L_2(X,C)$ of the operator $\Delta$ with  domain $C_c^\infty(\dot{X};C)$. Denote by $\Delta_{min}^*$ the adjoint to $\Delta_{min}$ in $L_2(X;C)$. In view of the Green formula
\begin{equation}
\label{Green formula}
(\Delta_L f,f')_{L_2(U;L)}-(f,\Delta_L f')_{L_2(U;L)}=\frac{1}{2i}\int_{\partial U}\big(fh\partial\overline{f'}dz+\overline{f'}h\overline{\partial}fd\overline{z}\big)
\end{equation}
with $U=X$ and $L=C$, we have $\Delta_{min}\subset\Delta_{min}^*$. 

Next, we describe all the self-adjoint extensions of $\Delta_{min}$. To this end, we consider the quotient space ${\rm Dom}\Delta_{min}^*/{\rm Dom}\Delta_{min}$ endowed with the symplectic bilinear form
\begin{equation}
\label{symplectic form}
\mathcal{G}(\pi(u),\pi(v)):=(\Delta_{min}^*u,v)_{L_2(X;C)}-(u,\Delta_{min}^*v)_{L_2(X;C)} \qquad (u,v\in{\rm Dom}\Delta_{min}^*),
\end{equation}
where $\pi$ is the canonical projection from ${\rm Dom}\Delta_{min}^*$ onto ${\rm Dom}\Delta_{min}^*/{\rm Dom}\Delta_{min}$. Let $\Delta\subset\Delta_{min}^*$ (i.e., $\Delta$ is the restriction of $\Delta_{min}^*$ on the domain ${\rm Dom}\Delta\subset{\rm Dom}\Delta_{min}^*$); then $\Delta$ is self-adjoint if and only if $\pi({\rm Dom}\Delta)$ is a lagrangian subspace of ${\rm Dom}\Delta_{min}^*/{\rm Dom}\Delta_{min}$ (see, e.g., Section 14.2 and Proposition 14.7, \cite{Schm}). These lagrangian subspaces are described via the following proposition (proved in paragraph \ref{par domain of adjoint}). 
\begin{prop}
\label{prop Darboux basis}
The Darboux basis on ${\rm Dom}\Delta_{min}^*/{\rm Dom}\Delta_{min}$ can be chosen in the form
\begin{equation}
\label{Darboux basis}
\{\pi(\chi_k f_{k,m,+})\}_{k=1,m=-1}^{2g-2, 2}, \ \ \{\pi(\chi_k f_{k,m,-})\}_{k=1,m=-1}^{2g-2, 2},
\end{equation}
where each $\chi_k$ is a smooth cut-off function equal to one near $P_k$, the support of $\chi_k$ is sufficiently small, and $f_{k,m,\pm}$ are local sections of $C$ given by
\begin{equation}
\label{anzats}
f_{k,m,+}(x_k):=x^m_k, \qquad f_{k,m,-}(x_k):=\frac{\overline{x_k}^{\frac{1}{2}-m} x_k^{\frac{1}{2}}}{\pi(m-\frac{1}{2})}
\end{equation}
in local coordinates {\rm(\ref{dist coord})}.
\end{prop}
We restrict ourselves to considering the following three self-adjoint extensions of $\Delta_{min}$. 

\smallskip

{\bf 1) Friedrichs extension} $\Delta_F$. The following proposition is proved in paragraph \ref{par Friedrichs}. It is a version of the corresponding result of Mooers (\cite{Mooers}) who dealt with Laplacians acting on $k$-forms.
\begin{prop}
\label{prop Friedrichs}
The domain of $\Delta_F$ consists of all sections $u\in{\rm Dom}\Delta_{min}^*$ admitting the asymptotics
\begin{equation}
\label{asymptotics Friedrichs}
u(x_k)=(-2/\pi)c_{k,0,-}|x_k|+c_{k,1,+}x_k+c_{k,2,+}x_k^2+(-2/3\pi)c_{k,-1,-}\overline{x_k}|x_k|+\tilde{u}(x_k)
\end{equation}
near each $P_k$, where $\tilde{u}\in{\rm Dom}\Delta_{min}$. 
\end{prop}

In Subsection \ref{ssec heat}, we derive the short-time asymptotics for the heat kernel $\mathscr{H}(\cdot,\cdot,t)$ corresponding to the equation $(\partial_t+\Delta_F)u=0$. As a corollary, we prove the following statement.
\begin{prop}
\label{prop zeta F wd}
The zeta function $s\mapsto\zeta(s|\Delta_F)$ of $\Delta_F$ admits analytic continuation on the whole complex plane except the pole $s=1$.
\end{prop}

\smallskip

{\bf 2) Szeg\"o extension} $\Delta_S$. Recall that, for smooth metrics $\rho^{-2}$ and $h$, the Szeg\"o kernel $\mathcal{S}$ is defined by 
\begin{equation}
\label{Green to Szego}
\mathcal{S}(z,z'):=-\pi h(z')\partial_{z'}\mathcal{G}(z,z'),
\end{equation}
where $\mathcal{G}$ is the Green function for the (self-adjoint) Laplacian $\Delta=-\rho^{2}h^{-1}\partial h\overline{\partial}$ (see p.25, \cite{Fay}). Also, the inversion formula
\begin{equation}
\label{Szego to Green}
\mathcal{G}(z,z'')=-\frac{1}{\pi^2}\int_{X}\mathcal{S}(z,z')h(z')\overline{\mathcal{S}(z',z'')}dS(z')
\end{equation}
holds. 
At the same time, for the case of spinor bundle $C$ obeying $h^0(C)=0$, the Szeg\"o kernel is independent of the choice of metrics $\rho^{-2}$ and $h$ and is given by the explicit formula (see p.29, \cite{Fay})
\begin{equation}
\label{Szego kernel expression}
\mathcal{S}(z,z')=\frac{\theta\left[^p _q\right]\big(\mathcal{A}_{P_0}(z'-z)\big)}{\theta\left[^p _q\right](0)E(z,z')}.
\end{equation}

In the non-smooth case, we define the Szeg\"o extension $\Delta_S$ as a self-adjoint extension of $\Delta_{min}$ whose Green function $\mathcal{G}=\mathcal{G}^{\mathcal{S}}$ is related to Szeg\"o kernel (\ref{Szego kernel expression}) via formulas (\ref{Green to Szego}) and (\ref{Szego to Green}). It is easy to see that such an extension is unique. Indeed, formula (\ref{Szego to Green}) uniquely determines the Green function $\mathcal{G}=\mathcal{G}^{\mathcal{S}}$ and, therefore, the kernel ${\rm Ker}\Delta_{S}$ and the inverse operator $\Delta_{S}^{-1}$ on the orthogonal complement ${\rm Ker}\Delta_{S}$ in $L_2(X;C)$. Thus, $\Delta_S$ is uniquely determined by $\mathcal{S}$ via (\ref{Szego to Green}).

The existence of the Szeg\"o extension and the description of its domain follow from Proposition \ref{prop Szego ext} proved in paragraph \ref{par Szego}.
\begin{prop}
\label{prop Szego ext}
The Szeg\"o extension $\Delta_S$ exists and its domain consists of all $u\in{\rm Dom}\Delta_{min}^*$ admitting the asymptotics 
\begin{equation}
\label{Szego domain}
u(x_k)=c_{k,0,+}+c_{k,1,+}x_k+c_{k,2,+}x_k^2+(-2/3\pi)c_{k,-1,-}\overline{x_k}|x_k|+\tilde{u}(x_k)
\end{equation}
near each $P_k$, where $\tilde{u}\in{\rm Dom}\Delta_{min}$. 
\end{prop}

\smallskip

{\bf 3) Holomorphic extension} $\Delta_h$. The domain of $\Delta_h$ consists of all $u\in{\rm Dom}\Delta_{min}^*$ admitting the asymptotics 
\begin{equation}
\label{holomorphic asymp}
u(x_k)=c_{k,-1,+}x_k^{-1}+c_{k,0,+}+c_{k,1,+}x_k+c_{k,2,+}x_k^2+\tilde{u}(x_k)
\end{equation}
near each $P_k$, where $\tilde{u}\in{\rm Dom}\Delta_{min}$. Note that the principal part of the above expansion contains only holomorphic local sections (\ref{anzats}). 

The peculiarity of the holomorphic extension $\Delta_h$ is that it coincides, up to unitary equivalence, with the self-adjoint operator
which can serve as a natural spectral theoretic counterpart of the  D'Hoker-Phong's Laplacian,
$\Delta^{\left(\frac{3}{2}\right)}$, 
 acting on spin-$\frac{3}{2}$ fields \cite{DH-P}.

Briefly, this can be explained as follows. The multiplication by $\omega$ is a unitary operator acting from $L_2(X;C)$ onto $L_2(X;C\otimes K)$. Thus, the minimal operators (with domains consisting of smooth sections with support in $\dot{X}$)  corresponding to $\Delta^{\left(\frac{3}{2}\right)}$ and $\Delta$ are unitary equivalent. The self-adjoint extension $\Delta^{\left(\frac{3}{2}\right),h}$ of the minimal operator $\Delta^{\left(\frac{3}{2}\right)}$ with domain containing the D'Hoker-Phong basis of holomorphic sections (we consider this requirement as the most natural)  turns out to be unitary equivalent to the holomorphic extension $\Delta_h$, 
\begin{equation*}
%\label{spin 3-2 holom}
\Delta^{\left(\frac{3}{2}\right),h}=\omega\Delta_h\omega^{-1}.
\end{equation*}
Thus, one can define
\begin{equation}
\label{Holom 3-2 and 1-2}
{\rm det}'\Delta^{\left(\frac{3}{2}\right)}:={\rm det}'\Delta^{\left(\frac{3}{2}\right),h}={\rm det}'\Delta_h.
\end{equation}
In paragraph \ref{par holom}, we prove the following statement.
\begin{prop}
\label{prop isospectrality of F and h}
The operators $\Delta_F$ and $\Delta_h$ are almost isospectral {\rm(}i.e., the non-zero eigenvalues of $\Delta_F$ and $\Delta_h$ coincide as well as their multiplicities{\rm)}. As a corollary, we have 
\begin{equation}
\label{Holom Fried 1-2}
{\rm det}'\Delta_h={\rm det}\Delta_F.
\end{equation}
\end{prop}

\smallskip

Next, we derive the comparison formulas for the $\zeta$-regularized determinants of $\Delta_F$, $\Delta_h$, $\Delta_S$. First, we have ${\rm det}'\Delta_h={\rm det}\Delta_F$ due to Proposition \ref{prop isospectrality of F and h}. In Section \ref{sec comparing det}, we compare the determinants of $\Delta_F$ and $\Delta_S$. To this end, we apply the method developed in \cite{HK} and based on the study of the behavior of the $S$-matrix $S(\lambda)$ associated with $\Delta_F-\lambda$ and $\Delta_S-\lambda$. Since many technicalities are different from those from \cite{HK} (where only the scalar Laplacians are considered), we give all the details making the presentation self-contained. As a result, we obtain the following theorem.
\begin{theorem}
\label{prop comparing determinants F S}
The zeta function $s\mapsto\zeta(s|\Delta_S)$ of $\Delta_S$ admits analytic continuation on the whole complex plane except the pole $s=1$. The determinants of $\Delta_S$ and $\Delta_F$ are related via
\begin{equation}
\label{comparing determinants F S}
{\rm det}\Delta_F=\Gamma(3/4)^{4(g-1)}{\rm det}T(0)\,{\rm det}\Delta_S,
\end{equation}
where $T(0)$ is the $(2g-2)\times(2g-2)$-matrix with entries
\begin{equation}
\label{expression for T(0)}
T_{ij}(0)=-\frac{1}{\pi^2|\theta\left[^p _q\right](0)|^2}\int\limits_{X}\frac{\theta\left[^p _q\right]\left(\mathcal{A}_{P_0}(P_i-z)\right)\overline{\theta\left[^p _q\right]\left(\mathcal{A}_{P_0}(z-P_j)\right)}}{E(P_i,z)\overline{E(z,P_j)}}\frac{dS_z}{|\omega(z)|}.
\end{equation}
\end{theorem}
Let us compare formula (\ref{comparing determinants F S}) with the results of D'Hoker-Phong's regularization \cite{DH-P}. Following \cite{DH-P}, introduce the ratios (\ref{QUANT}) for Laplacians acting on spin-$\frac{1}{2}$ and spin-$\frac{3}{2}$ fields as follows
$$W_{1/2}={\rm det}\Delta_S, \qquad W_{3/2}=\frac{{\rm det}'\Delta^{\left(\frac{3}{2}\right)}}{{\rm det}\{(\phi^{(3/2)}_i,\phi^{(3/2)}_j)_{L_2(X;S^3)}\}_{ij}},$$
where $\phi^{(3/2)}_j=\omega \mathcal{S}(P_j,\cdot)$ are zero modes of $\Delta^{\left(\frac{3}{2}\right)}$ (compare with formula (5.10), \cite{DH-P}) while the kernel of $\Delta_S$ is trivial. Since the multiplication by $\omega$ is a unitary operator acting from $L_2(X;C)$ onto $L_2(X;C\otimes K)$, we have 
$$W_{3/2}=\frac{{\rm det}'\Delta^{\left(\frac{3}{2}\right)}}{{\rm det}\mathfrak{S}},$$
where $\mathfrak{S}$ is the $(2g-2)\times(2g-2)$-matrix with entries 
\begin{equation}
\label{Bergman matrix}
\mathfrak{S}_{kj}=(\mathcal{S}(P_k,\cdot),\mathcal{S}(P_j,\cdot))_{L_2(X;C)}.
\end{equation}
which is related to $T(0)$ via
\begin{equation}
\label{S and T matrices}
\mathfrak{S}=\pi^2 T(0)
\end{equation}
(see formula (\ref{S and T matrices 1}) below). In view of Theorem \ref{prop comparing determinants F S}, the ratio
\begin{equation}\label{DPhF}
\begin{split}
\frac{W_{3/2}}{W_{1/2}}=\frac{{\rm det}'\Delta^{\left(\frac{3}{2}\right)}}{{\rm det}\Delta_S}\frac{1}{{\rm det}\mathfrak{S}}=&\big[(\ref{Holom 3-2 and 1-2}), (\ref{Holom Fried 1-2})\big]=\\
=&\frac{{\rm det}\Delta_F}{{\rm det}\Delta_S}\frac{1}{{\rm det}\mathfrak{S}}=\big[(\ref{comparing determinants F S}), (\ref{S and T matrices})\big]=\left(\frac{\Gamma(3/4)}{\pi}\right)^{4(g-1)}
\end{split}
\end{equation}
is moduli independent, which is, in a sense, compatible with D'Hoker-Phong formula (6.2), \cite{DH-P}. Namely, 
in \cite{DH-P} the left hand side of (\ref{DPhF}) is defined and computed as a tensor-like object depending on the choice of local parameters at the conical points. Naively choosing these parameters to be distinguished (see (\ref{dist coord})), one gets a moduli independent constant (which, unfortunately, differs from the right hand side of (\ref{DPhF})).   

\smallskip

Finally, we find the explicit expressions for ${\rm det}\Delta_S$ through the Bergman tau-function on the moduli space of Abelian differentials and the theta-constants (corresponding to the spinor bundle). To this end, we apply the variational technique developed in \cite{KokKorot}. First, we introduce the family of the surfaces corresponding to the variation of the coordinates $\nu$ in moduli space $H_g(1, \dots, 1)$. Next, we study the dependence of the resolvent kernels $\mathscr{R}$ of the (Szeg\"o) Laplacians on $\nu$ and define and calculate the derivatives $\mathscr{R}$ with respect to $\nu$. Based on this information, we calculate the derivative of the Laplacian eigenvalues and zeta functions with respect to $\nu$. Due to this, we compute the ${\rm det}\Delta_S$ up to a multiplicative moduli independent constant $c$. 
\begin{theorem}
\label{main theorem}
We have
\begin{equation}
\label{determinant expression}
{\rm det}\Delta_S(|\omega|,X)=c\left|\theta\left[^p _q\right](0)\right|^{2}|\tau_B(X,\omega)|^{-1}.
\end{equation}
\end{theorem}
The dependence of ${\rm det}\Delta_F={\rm det}'\Delta_h$ on moduli and spin structure is more sophisticated and is provided by the comparison of formulas (\ref{comparing determinants F S}) and (\ref{determinant expression}).
  
Comparing (\ref{determinant expression}) with (\ref{result}), we conclude that the method of the regularization of ${\rm det}\Delta$ described in Section \ref{sec reqularization} is compatible with the spectral regularization. Also, comparing (\ref{determinant expression}) with formula (1.10) from \cite{KokKorot}, we obtain the formula
$${\rm det}\Delta_S(|\omega|,X)=c\,\left(\frac{{\rm det'}\Delta_0(|\omega|^2,X)}{{\rm Area\,}(X, |\omega|^2) {\rm det}\,\Im {\mathbb B}}         \right)^{-1/2}|\theta[ ^p _q](0)|^2$$
which extends the Bost-Nelson bosonization result to the case of (non-smooth) metrics $|\omega|^2$ on X and $h=|\omega|$ on $C$. 

\begin{remark} It would be tempting to try to directly derive the above analog of Bost-Nelson bosonization formula for conical metrics from the standard one for smooth metrics via methods that use smoothing of the singularities (similar to those from \cite{Freixas}, \cite{Sher}, \cite{Finski} and \cite{Kok-PAMS}). However, we do not see how the proper self-adjoint extension (i. e. Szeg\"o one) can enter the game on this way if the latter is possible. This should become a subject of further study. \end{remark}
	
\begin{remark} It would be also challenging to search for possible higher dimensional generalizations of conical bosonization formula in the spirit of \cite{BL}.	
\end{remark}

\section{Preliminaries}
In this section, we present auxiliary facts which are used in the proofs of Theorems \ref{prop comparing determinants F S} and \ref{main theorem}. In particular, we prove Propositions \ref{prop Darboux basis}--\ref{prop isospectrality of F and h}. 

\subsection{Local properties of solutions to $(\Delta-\lambda)u=f$}
\label{subs Elliptic theory}
Let $u,f\in L_2(X;C)$. We say that $u$ is a generalized solution to the equation $(\Delta-\lambda)u=f$ if $(f,v)_{L_2(X;C)}=(u,(\Delta-\overline{\lambda})v)_{L_2(X;C)}$ for any $v\in C_c^\infty(X;C)$. Note that, in local coordinates
\begin{equation}
\label{z coordinate}
z(Q):=\int_{Q_0}^Q\omega,
\end{equation}
we have $\omega=dz$, $g=h=1$, and
\begin{equation}
\label{elipticity}
\Delta=-\frac{1}{4}(\partial^2_{\Re z}+\partial_{\Im z}^2).
\end{equation}
In this subsection, we describe the well-known properties (such as smoothness and asymptotics near conical points) of the solution $u$ arising from the ellipticity (in local coordinates (\ref{z coordinate})) of the Laplacian $\Delta$. The most used result is the following statement.
\begin{prop}
\label{prop elliptic asymptotics}
Let $u\in L_2(X;C)$ be a generalized solution to the equation $(\Delta-\lambda)u=f$ where $f\in L_2(X;C)$,
 $f\in C^{l}(\dot{X})$,  
 and
\begin{equation}
\label{ell asym rhs condition}
|\partial^{p}_{\overline{x_k}}\partial^q_{x_k}f(x_k)|\le c_f|x_k|^{T-p-q-4} \qquad (p+q\le l, \ T>5/2)
\end{equation}
near $P_k$, where $x_k$ is the distinguished coordinate {\rm (\ref{dist coord})}. Then $u\in C^{l}(\dot{X})$
and it admits the asymptotics 
\begin{equation}
\label{elliptic asymptotics}
u(x_k)=\sum_{(m,\pm)}c^{\lambda}_{k,m,\pm}(u)f_{k,m,\pm}(x_k)\sum_{n\ge 0}d(n,\pm i\mu_m)\left(\frac{\lambda|x_k|^{4}}{4}\right)^n+\tilde{u}(x_k)
\end{equation}
near $P_k$ with the remainder $\tilde{u}$ satisfying
\begin{equation}
\label{remainder estimate ellip}
|\partial^{p}_{\overline{x_k}}\partial^q_{x_k}\tilde{u}(x_k)|\le c_{\lambda,\varepsilon}\big(c_f+\|u\|_{L_2(X;C)}\big)|x_k|^{T-p-q-\varepsilon} \qquad (p+q\le l)
\end{equation}
for any $\varepsilon>0$. Here $f_{k,m,\pm}$ are the local sections of $C$ given by {\rm (\ref{anzats})}, the coefficients $c^{\lambda}_{k,m,\pm}(u)$ are linear functionals of $u$ given by
\begin{equation}
\label{elliptic asymptotics coeff}
c^{\lambda}_{k,m,\pm}(u)=(f,\chi_kf_{k,m,\mp})_{L_2(X;C)}-(u,(\Delta-\overline{\lambda})[\chi_k f_{k,m,\mp}])_{L_2(X;C)},
\end{equation}
where the cut-off function $\chi_k$ is defined before {\rm(\ref{anzats})} and $(\cdot,\cdot)_{L_2(X;C)}$ denotes the extension of the scalar product {\rm(\ref{L2 scalar})}. Also, $\mu_m=i(1-2m)/4$, and the numbers $d(n,q)$ are defined by the chain of equations 
\begin{equation}
\label{chain of equations pencil}
d(0,q)=1, \qquad -n(q+n)d(n,q)=d(n-1,q) \quad (n>0).
\end{equation}
The right hand side of {\rm(\ref{elliptic asymptotics})} includes those terms that are square integrable {\rm(}in metrics $\rho^{-2},h${\rm)} near the vertex $P_k$ and decrease slower than $O(|x_k|^{N})$ as $x_k\to 0$. The constant $c_{\lambda,\varepsilon}$ in {\rm(\ref{remainder estimate ellip})} has at most polynomial growth as $|\lambda|\to \infty$.
\end{prop}
Proposition \ref{prop elliptic asymptotics} is analogous to the statements describing  properties of solutions to elliptic boundary value problems in piecewise smooth domains presented in \cite{Kon,NP}.
Unfortunately, there is no appropriate reference for the case of closed manifolds,  and,  
 for the convenience of the reader, we give the full proof, using the key elements of the reasoning from \cite{Kon,NP} in a simplified and self-contained form. Next, applying the same reasoning to describe the asymptotics of solutions $\Delta_{min}^* u=f$, we prove Proposition \ref{prop Darboux basis}. Finally, we prove that the ellipticity of $\Delta$ implies the discreteness of spectra of self-adjoint expansions of $\Delta_{min}$, a well-known fact in the smooth case, but nevertheless needed to be verified in the non-smooth one.

\begin{nota} In the sequel, the neighborhoods $U_k$, the coordinates $z_k$, and the cut-off functions $\chi_k$, $\tilde{\chi}_k$ {\rm(}$k=1,\dots,K${\rm)} satisfy the following properties
\begin{itemize}
\item $U_k$-s constitute the finite open cover of $X$ and $P_k\in U_j$ if and only if $k=j$.
\item $z_k$ is given by {\rm(\ref{z coordinate})}, where the integration path lies in $U_k$ and $Q_0=P_k$ for $k\le 2g-2$.
\item the map $z_k: \ U_k\to\mathbb{C}$ is injective for $k>2g-2$.
\item the map $x_k: \ U_k\to\mathbb{C}$ is injective for $k\le 2g-2$, where $x_k$ is the distinguished coordinate {\rm(\ref{dist coord})}. Note that $z_k=x_k^2/2$.
\item $\chi_k,\tilde{\chi_k}$ are smooth cut-off functions on $X$ obeying 
$${\rm supp}\tilde{\chi}_k\subset U_k, \qquad \tilde{\chi}_k\chi_k=\chi_k, \qquad \sum_{k=1}^{K}\chi_k=1.$$
\item the $\chi_k(Q),\tilde{\chi_k}(Q)$ with $k\le 2g-2$ depend only on the distance between $Q$ and $P_k$.
\end{itemize}

\smallskip

For any section $v$ any line bundle $L$, we denote by $v^{(k)}$ its representatives in coordinates $z_k$. The representatives $v^{(k)}$ {\rm(}$k\le 2g-2${\rm)} near conical points $P_k$ are considered as functions of the polar coordinates 
$$r_k=|x_k|^2/2, \quad \varphi_k=2{\rm arg}x_k$$
or coordinates $s_k:={\rm log}r_k,\varphi_k$. Note that $z_k=r_k e^{i\varphi_k}=e^{s_k+i\varphi_k}$.

\smallskip

In addition, we use the cut-off functions $\chi_{k,\epsilon}$ on $X$ with supports shrinking to $P_k$ as $\epsilon\to 0$. Each $\chi_{k,\epsilon}$ vanish outside $U_k$ and is given by 
\begin{equation}
\label{shrinkin cutoff}
\chi_{k,\epsilon}(r_k,\varphi_k)=\chi(\varepsilon^{-1}r_k) 
\end{equation}
in $U_k$, where $\chi$ is a smooth cut-off function on $\mathbb{R}$ equal to one near zero.
\end{nota}

\subsubsection{Increasing smoothness outside conical points.} 
Let $k> 2g-2$ i.e. $U_k$ is separated from conical points. Since operator (\ref{elipticity}) is elliptic, the inclusion of $\tilde{\chi}_k f^{(k)}$ into the Sobolev space $H^l\equiv H^l(\mathbb{C})$ ($l=0,1,\dots$) implies the inclusion $\chi_k u^{(k)}\in H^{l+2}$ and the estimate
\begin{equation}
\label{H l smooth incr}
\|\chi_k  u^{(k)}\|_{H^{l+2}}\le c_l(\|\tilde{\chi}_k f^{(k)}\|_{H^{l}}+\|\tilde{\chi}_k u\|_{L_2(X;C)}).
\end{equation}
Also, due to Morrey's inequality we have $u^{(z_k)}\in C^{l}\equiv C^{l}(\mathbb{C})$ and
\begin{equation}
\label{Morrey}
\|\chi_k u^{(k)}\|_{C^{l}}\le c_l\|\tilde{\chi}_k u^{(k)}\|_{H^{l+2}}.
\end{equation}
Combining the two last estimates, we obtain
\begin{equation}
\label{C l smoth incr}
\|\chi_k u^{(k)}\|_{C^{l}}\le c_l(\|\tilde{\chi}_k f^{(k)}\|_{C^{l}}+\|\tilde{\chi}_k u\|_{L_2(X;C)}).
\end{equation}
In addition, if $f$ is smooth outside conical points, then so is $u$ and the equation $(\Delta-\lambda)u=f$ holds pointwise on $\dot{X}$.

\subsubsection{Asymptotics near conical points.} 
\paragraph{\it Reduction to a model problem.} Let $k\le 2g-2$ i.e. $U_k$ is a neighborhood of $P_k$. Throughout the paragraph, the subscript $k$ is usually omitted in the notation of coordinates $x_k,z_k,r_k,\varphi_k$ and cut-off functions $\chi_k$, $\tilde{\chi}_k$. Also, we assume that $\lambda=0$ (the general case will be considered later). Put 
\begin{equation}
\label{ellip local}
y(r,\varphi):=\chi u(z), \quad w(r,\varphi):=-4r^{2}\big(\chi f(z)+[\Delta,\chi]u(z)\big).
\end{equation}
In view of (\ref{elipticity}), the equation $\Delta u=f$ implies
\begin{equation}
\label{equation near c p}
[\partial_{s}^2+\partial_{\varphi}^2]y=[(r\partial_{r})^2+\partial_{\varphi}^2]y=w.
\end{equation}
Also, since the transition function of the spinor bundle $C=\sqrt{K}$ corresponding to the change of variables $x\mapsto z$ is equal to $$\pm\sqrt{\partial x/\partial z}=\pm 2^{-\frac{1}{4}} r^{-\frac{1}{4}}e^{-\frac{i\varphi}{4}},$$ 
the representatives %local trivializations 
of $u,f$ in coordinate $z$ are $4\pi$-antiperiodic in $\varphi$, i.e.,
\begin{equation}
\label{antiperiodicity}
v(r,\varphi+\alpha)=-v(r,\varphi) \qquad (v=y,w, \quad \alpha=4\pi).
\end{equation}
Applying the complex Fourier transform
$$\hat{v}(\mu,\varphi):=\frac{1}{\sqrt{2\pi}}\int_{-\infty}^{+\infty}e^{-i\mu s}v(s,\varphi)ds,$$
one can rewrite (\ref{equation near c p}), (\ref{antiperiodicity}) as
\begin{equation}
\label{operator pencil}
[\partial_{\varphi}^2-\mu^2]\hat{y}(\mu,\varphi)=\hat{w}(\mu,\varphi), \qquad \hat{y}(\mu,\varphi+\alpha)=-\hat{y}(\mu,\varphi).
\end{equation}
The resolvent kernel of problem (\ref{operator pencil}) is given by
\begin{equation}
\label{operator pencil resolvent}
\mathfrak{R}_{\alpha}(\varphi,\varphi',\mu^2):=\frac{{\rm sh}(\mu(|\varphi-\varphi'|-\alpha/2))}{2\mu{\rm ch}(\alpha\mu/2)} \qquad (\varphi,\varphi'\in [0,\alpha)).
\end{equation}
(Here ${\rm sh}(\mu\dots)$ is a solution to $(\partial_{\varphi}^2-\mu^2)\hat{y}=0$ outside $\varphi=\varphi'$, the absolute value of $\varphi-\varphi'$ is taken to make the derivative of the numerator have a jump at $\varphi-\varphi'$ so that $(\partial_{\varphi}^2-\mu^2){\rm sh}(\mu(|\varphi-\varphi'|-\alpha/2))$ is proportional to $\delta(\varphi-\varphi')$ with the coefficient canceled by the denominator; the argument shift by $-\alpha/2$ is introduced to ensure the condition $\hat{y}(\mu,\varphi+\alpha)=-\hat{y}(\mu,\varphi))$ for the numerator at $\varphi=0$.)
It is meromorphic function of $\mu$ with simple poles
\begin{equation}
\label{OP resolvent poles}
\mu_m=i\pi \alpha^{-1}(1-2m)=i(1-2m)/4 \qquad (m\in\mathbb{Z}).
\end{equation}
The estimate
\begin{equation}
\label{OP resolvent at infinity}
\|\mathfrak{R}_{\alpha}(\cdot,\cdot,\mu^2)\|_{L_2([0,\alpha]^2)}=O(|\mu|^{-2}), \qquad \|\partial_\varphi\mathfrak{R}_{\alpha}(\cdot,\cdot,\mu^2)\|_{L_2([0,\alpha]^2)}=O(|\mu|^{-1}),
\end{equation}
is valid for large $|\Re\mu|$.

\

\paragraph{\it Weighted spaces.} Denote by $H_\nu^l$ the space of $\alpha$-antiperiodic (with respect to $\varphi$) functions with finite norms
\begin{align}
\label{weighted norm}
\begin{split}
\|v\|_{H_\nu^l}=\left(\sup_{\phi_0}\sum_{p+q\le l}\int\limits_{0}^{+\infty}rdr\int\limits_{\phi_0}^{\phi_0+\alpha} d\varphi |\partial^q_\varphi (r\partial_r)^p v(r,\varphi)|^2 r^{2(\nu-1)}\right)^{1/2}=\\
=\left(\sup_{\phi_0}\sum_{p+q\le l}\int\limits_{-\infty}^{+\infty}ds\int\limits_{\phi_0}^{\phi_0+\alpha} d\varphi |\partial^q_\varphi\partial_s^p v(s,\varphi)|^2 e^{2\nu s}\right)^{1/2}.
\end{split}
\end{align}
In view of the Parseval identity, we have the following equivalence of the norms
\begin{equation}
\label{Parseval equiv}
\|v\|_{H_\nu^l}\asymp\left(\sup_{\phi_0}\sum_{p+q\le l}\int\limits_{\nu i-\infty}^{\nu i+\infty}|\mu|^{2p}\|\partial^q_\varphi\hat{v}(\mu,\varphi)\|_{L_2([\phi_0,\phi_0+\alpha])}^2 d\mu\right)^{1/2}.
\end{equation}
As it easy to see, the convergence 
$$\|(1-\chi_\epsilon)\chi_{1/\epsilon}v-v\|_{H_\nu^l}\to 0, \qquad \epsilon\to 0$$
is valid for any $v\in H_\nu^l$, where the cut-off function $\chi_{\epsilon}$ is given by (\ref{shrinkin cutoff}) with omitted $k$. Hence, we obtain the following fact.
\begin{remark}
\label{density in cylinder}
The set of smooth $\alpha$-antiperiodic {\rm(}with respect to $\varphi${\rm)} functions vanishing for sufficiently small and large $r$ is dense in $H_\nu^l$.
\end{remark}

Let us provide the analogues of inequality (\ref{H l smooth incr}), (\ref{Morrey}) for weighted spaces $H^l_\nu$. Let $\kappa,\varkappa$ be smooth cut-off functions on $\mathbb{R}$ such that $\kappa\varkappa=\kappa$, $\kappa(s)=1$ for $s\in[-3/2,3/2]$ and $\varkappa(s)=0$ for $s\in[-2,2]$. For $j\in\mathbb{Z}$, denote $\kappa_j(s):=\kappa(s-j)$ and $\varkappa_j(s):=\varkappa(s-j)$. Since the operator $\partial_s^2+\partial_\varphi^2$ is elliptic, the estimate
$$\|\kappa_j v\|_{H^{l+2}(\Pi)}\le c(\|\varkappa_j v\|_{L_2(\Pi)}+\|\varkappa_j (\partial_s^2+\partial_\varphi^2)v\|_{H^{l}(\Pi)})$$
is valid, where $\Pi=\mathbb{R}\times[\phi_0,\phi_0+\alpha]$. Multiplying both parts by $e^{2\nu j}$ and making the summation over integer $j$, we arrive at the estimate
\begin{equation}
\label{smoth inc weight}
\|v\|_{H_\nu^{l+2}}\le c(\|v\|_{H_\nu^{0}}+\|(\partial_s^2+\partial_\varphi^2)v\|_{H_\nu^{l}}).
\end{equation}
Also, due to the Morrey's inequality we have
\begin{equation}
\label{weighted Morrey}
\begin{split}
|\partial^q_\varphi (r\partial_r)^p v(r,\varphi)|=|\partial^q_\varphi \partial_s^p v(r,\varphi)|\le c\|\kappa_j v\|_{H^{l+2}(\Pi)}\le \\ \le c\|v\|_{H^{l+2}_\nu}e^{-\nu k}\le c \|v\|_{H^{l+2}_\nu}r^{-\nu}, \qquad (p+q\le l)
\end{split}
\end{equation}
(here $j$ is chosen in such a way that $r\in [j-3/2,j+3/2]$).

\

\paragraph{\it Problem {\rm(\ref{equation near c p}), (\ref{antiperiodicity})} in the scale of weighted spaces.} 
With problem (\ref{equation near c p}), (\ref{antiperiodicity}), we associate the continuous operator $A_{\nu,l}: y\mapsto (\partial_s^2+\partial_{\varphi}^2)y$ acting from $H^{l+2}_\nu$ to $H^{l}_\nu$.

Suppose that the line $\Im\mu=\nu$ does not intersect the poles $\mu_n$ given by (\ref{OP resolvent poles}). If $y\in H^{l+2}_\nu$, then $w=A_{\nu,l} y\in H^{l}_\nu$ and the complex Fourier transforms $\hat{y}$, $\hat{w}$ are defined and obey (\ref{operator pencil}) for $\Im\mu=\nu$ and almost all $\Re\mu\in\mathbb{R}$. Expressing $\hat{y}$ in terms of $\hat{w}$ and the resolvent kernel $\mathfrak{R}_{\alpha}$ (given by (\ref{operator pencil resolvent})) and then applying the inverse Fourier transform, we obtain
\begin{equation}
\label{solution in weighted space}
y(s,\varphi)=\frac{1}{\sqrt{2\pi}}\int\limits_{\nu i-\infty}^{\nu i+\infty}d\mu \int\limits_{0}^{\alpha}d\varphi' e^{i\mu s}\mathfrak{R}_{\alpha}(\varphi,\varphi',\mu^2)\hat{w}(\mu,\varphi').
\end{equation}
If it is known only that $w\in H^{l}_\nu$, then formula (\ref{solution in weighted space}) provides the (unique) solution $y=A_{\nu,l}^{-1}w$ to (\ref{equation near c p}), (\ref{antiperiodicity}) which belongs to $y\in H^{l+2}_\nu$. Indeed, if follows from (\ref{OP resolvent at infinity}) and the Cauchy–Schwarz inequality that
$$|\mu|^2\|\hat{y}(\mu,\cdot)\|_{L_2([0,\alpha])}+|\mu|\|\partial_\varphi\hat{y}(\mu,\cdot)\|_{L_2([0,\alpha])}\le c\|\hat{w}(\mu,\cdot)\|_{L_2([0,\alpha])}.$$
Since equation (\ref{operator pencil}) implies $\partial_{\varphi}^{2(k+1)}\hat{y}=\partial_{\varphi}^{2k}\hat{w}+\mu^2\partial_{\varphi}^{2k}\hat{y}$, the iteration of the last formula yields
\begin{align*}
\|\partial_{\varphi}^{2(k+1)+q}\hat{y}(\mu,\cdot)\|_{L_2([0,\alpha])}\le |\mu|^2\|\partial_{\varphi}^{2k+q}\hat{y}(\mu,\cdot)\|_{L_2([0,\alpha])}+\|\partial_{\varphi}^{2k+q}\hat{w}(\mu,\cdot)\|_{L_2([0,\alpha])}\le\\
\le \sum_{m=0}^{2k}|\mu|^{2m+q}\|\partial_\varphi^{2(k-m)+q}\hat{w}(\mu,\cdot)\|_{L_2([0,\alpha])} \qquad (q=0,1, \ 2(k+1)+q\le l+2).
\end{align*}
Since $w\in H^{l}_\nu$ the right hand side of (\ref{Parseval equiv}) with $w,\hat{w}$ instead of $v,\hat{v}$ is finite. Thus, due to above estimates, the right hand side of (\ref{Parseval equiv}) with $y,\hat{y},l+2$ instead of $v,\hat{v},l$ obeys
\begin{equation}
\label{w spaces solution estim}
\|y\|_{H^{l+2}_\nu}\le c_\nu\|w\|_{H^{l}_{\nu}},
\end{equation}
i.e., $A_{\nu,l}^{-1}$ is continuous if $\nu$ does not coincide with any $\Im\mu_m$.

Now, suppose that $\nu,\nu'\ne\Im\mu_n$ for any $n$, $\nu'<\nu$, and $y\in H^{l+2}_\nu$ and $y'\in H^{l+2}_{\nu'}$ are two solutions to (\ref{equation near c p}), (\ref{antiperiodicity}) with the joint right hand side $w=A_{\nu,l} y=A_{\nu',l}y'$. Then $w\in H^{l}_{\nu''}$ for any $\nu''\in [\nu',\nu]$ and $\mu\mapsto\hat{w}(\mu,\cdot)$ is holomorphic in the strip $\Im\mu\in [\nu',\nu]$. Then formulas (\ref{solution in weighted space}), (\ref{operator pencil resolvent}), (\ref{OP resolvent poles}), (\ref{OP resolvent at infinity}) and the residue theorem imply
\begin{equation}
\label{ellyptic asymptotics in cone}
\begin{split}
y-y'=\frac{1}{\sqrt{2\pi}}\oint_{\mathcal{C}}d\mu \int\limits_{0}^{4\pi}d\varphi' e^{i\mu s}\mathfrak{R}_{4\pi}(\varphi,\varphi',\mu^2)\hat{w}(\mu,\varphi')&=\\
=\sqrt{2\pi}i\sum_{\Im \mu_m\in [\nu',\nu]}\int\limits_{0}^{4\pi}d\varphi'\underset{\mu=\mu_m}{\rm Res}\big(e^{i\mu s}\mathfrak{R}_{4\pi}(\varphi,\varphi',\mu^2)&\hat{w}(\mu,\varphi')\big)=\\
=\sum_{\Im \mu_m\in [\nu',\nu]}\sum_{\pm}c_{\pm,m}&(w)r^{\frac{2m-1}{4}}e^{\pm i\frac{2m-1}{4}\varphi},
\end{split}
\end{equation}
where $\mathcal{C}$ is the (positively oriented) boundary of the strip $\Im\mu\in [\nu',\nu]$, and the coefficients $c_{\pm,m}(w)$ are given by
\begin{equation}
\label{coeff formula in cone}
c_{\pm,m}(w)=\frac{-1}{4\pi(\frac{1}{2}-m)}\int\limits_{0}^{4\pi}d\varphi'\int\limits_{-\infty}^{+\infty}ds\, w(s,\varphi)\overline{r^{-\frac{2m-1}{4}}e^{\pm i\frac{2m-1}{4}\varphi'}}.
\end{equation}

Now, we are ready to prove Proposition \ref{prop elliptic asymptotics} using the facts above.

\

\paragraph{\it Asymptotics of $u$ near $P_k$: the case $\lambda=0$.} Condition (\ref{ell asym rhs condition}) can be rewritten as 
$$|\partial^{q}_{\varphi}(r\partial_{r})^{p}(r^2f^{(k)}(r,\varphi))|\le c_fr^{-\nu_0} \qquad (p+q\le l),$$
where 
$$\nu_0=\frac{1-2T}{4}<-1.$$
Also, the support $[\Delta,\chi_k]u$ is separated from $P_k$. Due to these facts and estimate (\ref{H l smooth incr}), the right hand side of (\ref{equation near c p}) (given by (\ref{ellip local})) satisfies 
\begin{equation}
\label{ellip rhs est}
\|w\|_{H^{l}_{\nu'}}\le c(c_f+\|u\|_{L_2(X;C)}) \qquad (\nu'>\nu_0).
\end{equation}
In view of (\ref{smoth inc weight}), the inclusion $u\in L_2(X;C)$ means that $y=\chi_k u^{(k)}\in H^{l+2}_1$ and, hence $A_{1,l}y=w$. For arbitrary $\nu'\in (\nu_0,\nu=1)$, put $y'=A_{\nu',l}^{-1}w$. Then formula (\ref{ellyptic asymptotics in cone}) provides the asymptotics of $u(z)$ as $z\to 0$. Rewriting (\ref{ellyptic asymptotics in cone}) in local coordinate $x=x_k$ and taking into account that the formulas
\begin{equation}
\label{anzats in z}
\begin{split}
f^{(k)}_{k,m,+}(r_k,\varphi_k)&=(\pm)(2r_k)^{\frac{2m-1}{4}}e^{\frac{i\varphi_k(2m-1)}{4}}, \\
f^{(k)}_{k,m,-}(r_k,\varphi_k)&=(\pm)\pi^{-1}\Big(\frac{1}{2}-m\Big)^{-1}(2r_k)^{\frac{1-2m}{4}}e^{\frac{i\varphi_k(2m-1)}{4}}
\end{split}
\end{equation}
are valid for the representatives %local trivializations 
of $f_{k,m,+}$ given by (\ref{anzats}), we arrive at expansion (\ref{elliptic asymptotics}) with $\lambda=0$, where the remainder $\tilde{u}$ obeys $\tilde{u}^{(k)}=y'$ for small $z=z_k$. Due to (\ref{w spaces solution estim}), (\ref{weighted Morrey}), the equality $y'=A_{\nu',l}^{-1}w$ implies
\begin{align*}
r^{\nu'}|\partial^q_\varphi (r\partial_r)^p \tilde{u}^{(k)}(r,\varphi)|\le c\|y'\|_{H^{l+2}_{\nu'}}\le c_\nu\|w\|_{H^{l}_{\nu'}} \qquad (p+q\le l).
\end{align*}
Rewriting the last condition in the local coordinate $x=x_k$ and taking into account (\ref{ellip rhs est}), we obtain estimate (\ref{remainder estimate ellip}).

\
 
\paragraph{\it Asymptotics of $u$ near $P_k$: the general case.} The general case is reduced to the case $\lambda=0$ by rewriting the equation $(\Delta-\lambda)u=f$ in the form $\Delta u=f+\lambda u$. However, replacing $f$ with $f+\lambda u$ in the above reasoning yields that the right hand side $w$ of (\ref{equation near c p}) contains the additional term $-4r^2\lambda u^{(k)}$, and, hence, estimate (\ref{ellip rhs est}) and formula (\ref{ellyptic asymptotics in cone}) are valid only for $\nu'\ge-1$. Therefore, the repeating of the reasoning above provides only the expansion
\begin{equation}
\label{nonhom zero approx}
u(x)=\sum_{(m,\pm)}c^{\lambda}_{k,m,\pm}(u)f_{k,m,\pm}(x)+\tilde{u}(x)
\end{equation}
where the sumation is taken over only terms decaying no faster than $O(x^{2})$ while the remainder obeys 
\begin{equation}
\label{rem zero approx}
|\partial^{p}_{\overline{x}}\partial^q_{x}\tilde{u}(x)|\le c_{\varepsilon}\big(c_f+(1+\lambda)\|u\|_{L_2(X;C)}\big)|x|^{3-p-q-\varepsilon} \qquad (p+q\le l, \ \varepsilon>0).
\end{equation}

To overcome this obstacle and obtain the next terms of the asymptotics of $u$ near $P_k$, we use the following trick. For each $f_{k,m,\pm}$ presented in the sum in (\ref{nonhom zero approx}), we construct the formal series
\begin{equation}
\label{ellip form series}
f_{k,m,\pm}(1+d(1,\pm i\mu_m)\lambda|z|^2+d(1,\pm i\mu_m)\lambda^2|z|^4+\dots)
\end{equation}
obeying the equation $(\Delta-\lambda)v=0$ near $P_k$; then the coefficients $d(n,\pm i\mu_m)$ in (\ref{ellip form series}) satisfy (\ref{chain of equations pencil}). Now, consider the section
$$u_1=u-\chi_k\sum_{(m,\pm)}c^{\lambda}_{k,m,\pm}(u)f^{\lambda,N}_{k,m,\pm}$$
instead of $u$, where $f_{k,m,\pm}^{\lambda,N}$ are truncated series (\ref{ellip form series}) with sufficiently large number, $N$, of terms. Note that
$$(\Delta -\lambda)f_{k,m,\pm}^{\lambda,N}=-\lambda^{N+1}f_{k,m,\pm}d(1,\pm i\mu_m)|z|^{4N}.$$
Also, 
$$|c^{\lambda}_{k,m,\pm}(u)|\le c_f(1+|\lambda|\|u\|_{L_2(X;C)})$$
due to (\ref{coeff formula in cone}). Thus, $(\Delta-\lambda)u_1$ satisfies condition (\ref{ell asym rhs condition}) with $c_{f}$ replaced by $c(c_f+|\lambda|^{N+2}\|u\|_{L_2(X;C)})$. Now, we repeat the above reasoning for $u_1$, $f_1=\Delta u_1$ instead of $u,f$. Since (\ref{rem zero approx}) and (\ref{ell asym rhs condition}) imply
$$\|w\|_{H^{l}_{\nu'}}\le c(c_f+(1+|\lambda|^{N+2})\|u\|_{L_2(X;C)}) \qquad (\nu'>\max\{-13/4,\nu_0\}),$$
we arrive at the asymptotics
$$u_1(x)=\sum_{(m,\pm)}c^{\lambda}_{k,m,\pm}(u_1)f_{k,m,\pm}(x)+\tilde{u}_1(x)$$
where the summation is taken over the terms decaying faster than $O(x^{2})$ but slower than $O(x^{\min\{T,7\}})$, while the remainder obeys 
$$|\partial^{p}_{\overline{x}}\partial^q_{x}\tilde{u}_1(x)|\le c_{\lambda,\varepsilon}\big(c_f+(1+\lambda)\|u\|_{L_2(X;C)}\big)|x|^{\min\{T,7\}-p-q-\varepsilon} \qquad (p+q\le l, \ \varepsilon>0),$$
where $c_{\lambda,\varepsilon}$ is a polynomial in $|\lambda|$.

If $T>7$, then we can repeat the above trick to find the next terms of the asymptotics of $u$ near $P_k$ and so on. As a result, we prove expansion (\ref{elliptic asymptotics}) and estimate (\ref{remainder estimate ellip}). Now, denote by 
\begin{equation*}
%\label{surface with holes}
X_{\epsilon}(p_1,\dots,p_n)
\end{equation*}
the domain obtained by deleting $\epsilon$-neighborhoods (in the metric $\rho^{-2}$) of points $p_1,\dots,p_n$ from $X$. Substitute $U=X_\epsilon(P_1,\dots,P_{2g-2})$, $f=u$, and $f'=\chi_k f_{k,m,\mp}$ in the Green formula (\ref{Green formula}) and pass to the limit $\epsilon\to 0$. Calculating the limit of the right hand side of (\ref{Green formula}) by the use of asymptotics (\ref{elliptic asymptotics}), we arrive at (\ref{elliptic asymptotics coeff}). Proposition \ref{prop elliptic asymptotics} is proved.

At the end of the paragraph, we formulate (in a simplified setting) an analogue of Proposition \ref{prop elliptic asymptotics} for surfaces with general conical singularities.
\begin{prop}
Let $m=\rho^{-2}(z)|dz|^2$ be a conformal metric on $X$ such that a neighborhood $U$ of some point $O$ of $X$ is isometric to a neighborhood of the cone of angle $\beta$, i.e., $m=|x^bdx|^2$ in some holomorphic coordinate $x$ obeying $x(P)=0$, where $b=\frac{\beta}{2\pi}-1$. Let also $h(z)|dz|^{-1}=\rho|dz|^{-1}$ be the corresponding metric on $C$. For simplicity, assume that $b$ is either positive or irrational. Let $u$ be a 
local square integrable section of $C$ in $U$ and let it be a generalized solution to the equation $(\Delta-\lambda)u=0$ in $U\backslash\{P\}$. Then $u$ is smooth in $U\backslash\{P\}$ and admits the convergent asymptotic expansion
\begin{equation}
\label{revisooor}
u(x)=\sum_{p,q,j}c^{\lambda}_{p,q}(u)d(p,q,j,b)Y_{p+j(b+1),q+j(b+1)}(x)
\end{equation}
near $P$. Here $Y_{p,q}(x)=x^p\overline{x}^q$, the coefficients $c^{\lambda}_{p,q}(u)$ are linear functionals of $u$. The sum in the right-hand side of {\rm(\ref{revisooor})} includes only the terms which are square integrable near $P$ and obey one of the following conditions: {\rm a)} $q=0$, $p$ is integer, {\rm b)} $p-b/2=0$, $q-b/2$ is integer. The numbers $d(p,q,j,b)$ are defined by the chain of equations 
\begin{equation}
\label{revisssoor 1}
d(p,q,0,b)=1, \qquad -(q+n(b+1))(p+n(b+1)-b/2)d(p,q,j,b)=d(p,q,j-1,b) \quad (n>0).
\end{equation}
Both sides of {\rm (\ref{revisooor})} can be differentiated in $\Re x,\,\Im x$ any number of times. Note that the assumption $b>0$ or $b\not\in\mathbb{Q}$ is needed only to to ensure the solvability of chain {\rm (\ref{revisssoor 1})}; in the general case series {\rm(\ref{revisooor})} may contain power-logarithmic terms {\rm(}for details, see Lemma 3.5.11, \cite{NP}{\rm)}.
\end{prop}

\subsubsection{Description of ${\rm Dom}\Delta_{min}^*$}
\label{par domain of adjoint}
In this paragraph, we prove Proposition \ref{prop Darboux basis} describing the structure of the space ${\rm Dom}\Delta_{min}^*/{\rm Dom}\Delta_{min}$.

Introduce the space $\mathcal{H}_\nu^l$ of sections of $C$ with finite norms
\begin{align}
\label{weighted spaces}
\|v\|_{\mathcal{H}_\nu^l}=\left(\sum_{k\le 2g-2}\|\chi_k v^{(k)}\|_{H^l_\nu}^2+\sum_{k>2g-2}\|\chi_k v^{(k)}\|^2_{H^l}\right)^{1/2}.
\end{align}
Note that the embedding $\mathcal{H}_{\nu'}^{l'}\subset\mathcal{H}_{\nu}^{l}$ is continuous for $l'\ge l$ and $\nu'\le\nu$. Also, we have the following equivalence of the norms 
\begin{equation*}
%\label{L 2 is H 0 1}
\|\cdot\|_{L_2(X;C)} \asymp \|\cdot\|_{\mathcal{H}_1^0}.
\end{equation*}
Note that $\mathcal{H}_\nu^l$ is complete. Moreover, as a corollary of Remark \ref{density in cylinder}, the set $C_c^\infty(\dot{X};C)$ is dense in $\mathcal{H}_\nu^l$ for each $l,\nu$. In view of definitions (\ref{weighted spaces}), (\ref{weighted norm}), the convergence $v_n\to v$ in $\mathcal{H}_{-1}^2$ implies that $v_n\to v$ and $\Delta v_n\to \Delta v$ in $L_2(X;C)$. Therefore, we have $\mathcal{H}_{-1}^2\subset {\rm Dom}\Delta_{min}$.

Conversely, let $v_n\in C_c^\infty(\dot{X};C)$ and $v_n\to v$ and $\Delta v_n\to f$ in $L_2(X;C)$. Then $v_n$ is a Cauchy sequence in $\mathcal{H}^2_{-1}$ due to inequality (\ref{H l smooth incr}) and estimate (\ref{w spaces solution estim}), where $y=\chi_k v^{(k)}_{n'}-\chi_k v^{(k)}_n$ and $w$ is given by (\ref{equation near c p}). Since $\mathcal{H}^2_{-1}$ is complete, we have $v_n\to v'$ in $\mathcal{H}^2_{-1}$. Hence, $v_n\to v'$ and $\Delta v_n\to \Delta v'$ in $L_2(X;C)$ and $v=v'\in \mathcal{H}^2_{-1}$, $f=\Delta v'$. Thus, ${\rm Dom}\Delta_{min}\subset\mathcal{H}_{-1}^2$. So, we have proved that
\begin{equation}
\label{emb to dom}
{\rm Dom}\Delta_{min}=\mathcal{H}^2_{-1}.
\end{equation}

Let $\Delta_{min}^*u=f$; then $u,f\in L_2(X;C)$ and $u$ is a generalized solution to $\Delta u=f$. So, estimates (\ref{H l smooth incr}) hold for $k>2g-2$ i.e. $u\in H^2(\Omega)$ for any domain $\Omega$ in $X$ whose closure does not contain conical points. Let us describe the asymptotics of $u$ near conical point $P_k$. In view of (\ref{H l smooth incr}), $\chi_k u^{(s)}$ (considered as a function of local coordinate (\ref{z coordinate}))  Let us describe the asymptotics of $u$ near each conical point $P_k$. Let $y,w$ be functions given by (\ref{ellip local}) with $z=z_k$, $r=r_k$, $\varphi=\varphi_k$. Then $w\in H^0_{-1}$, $y\in H^0_1$, and, moreover, $y\in H^2_1$ due to (\ref{smoth inc weight}). So, we have
\begin{equation}
\label{emb to dom adj}
{\rm Dom}\Delta_{min}^*\subset \mathcal{H}_1^2.
\end{equation}
Now formula (\ref{ellyptic asymptotics in cone}) with $l=2$, $\nu=1$, $\nu'=0$ yields the expansion
$$\chi_k u^{(k)}(r_k,\varphi_k)=\sum_{m=-1}^{2}\sum_{\pm}c_{k,m,\pm}(u)\chi_k f_{k,m,\pm}^{(k)}(r_k,\varphi_k)+y'(r_k,\varphi_k)$$
with remainder obeying $y'\in H^2_{-1}$. In other words, we have
\begin{equation}
\label{asymptotics near all vertices}
u=\sum_{k=1}^{2g-2}\sum_{m=-1}^{2}\sum_{\pm}c_{k,m,\pm}(u)\chi_k f_{k,m,\pm}+\tilde{u},
\end{equation}
where $\tilde{u}\in \mathcal{H}_{-1}^2\subset {\rm Dom}\Delta_{min}$ due to (\ref{emb to dom}). Thus, the elements (\ref{Darboux basis}) constitute a basis in ${\rm Dom}\Delta_{min}^*/{\rm Dom}\Delta_{min}$. 

Now let $f,f'\in{\rm Dom}\Delta^*$; then the expansions (\ref{asymptotics near all vertices}) are valid for them. Then from the Green formula (\ref{Green formula}) with $U=X_\epsilon(O_1,\dots,O_{2g-2})$, $\epsilon\to 0$, it follows that
\begin{equation}
\label{Green formula with asymptotics}
\begin{split}
(\Delta^* f,f')_{L_2(X;C)}-(f,\Delta^* f')_{L_2(X;C)}&=\\
=\sum_{k=1}^{2g-2}\sum_{m=-1}^{2}&\big(c_{k,m,+}(f)\overline{c_{k,m,-}(f')}-c_{k,m,-}(f)\overline{c_{k,m,+}(f')}\big).
\end{split}
\end{equation}
Comparing (\ref{Green formula with asymptotics}) with (\ref{symplectic form}), we obtain the equalities
$$\mathcal{G}(\pi(\chi_k f_{k,m,s}),\pi(\chi_k f_{k',m',-s'}))=\delta_{k,k'}\delta_{m,m'}\delta_{s,s'}$$
from which it follows that (\ref{Darboux basis}) is a Darboux basis on ${\rm Dom}\Delta_{min}^*/{\rm Dom}\Delta_{min}$. Proposition \ref{prop Darboux basis} is proved.

\subsubsection{Discreteness of spectra}
\label{par discretness spectra}
In this paragraph, we prove the following statement.
\begin{prop} 
\label{discreteness spectrum}
The spectrum of any self-adjoint extension of of $\Delta_{min}$ is discrete.
\end{prop}
The scheme of the proof is similar to that used in Chapter 4, \cite{NP}. Let $\Delta_1$ be a self-adjoint extension of of $\Delta_{min}$; then $(\Delta_1-i)^{-1}$ is a continuous operator in $L_2(X;C)\equiv \mathcal{H}_1^0$. Let $\nu_0\in (3/4,1)$. We first prove that $(\Delta_1-i)^{-1}$ acts continuously from $\mathcal{H}_1^0$ onto $\mathcal{H}_{\nu_0}^2$. Second, we prove that the embedding $\mathcal{H}_{\nu_0}^2\subset \mathcal{H}_{1}^0\equiv L_2(X;C)$ is compact. Due to these facts, $(\Delta_1-i)^{-1}$ is a compact operator in $L_2(X;C)$ and, hence, its spectrum is discrete. Thus, the spectrum of $\Delta_1$ is also discrete.

\

\paragraph{\it Continuity of $(\Delta_1-i)^{-1}$ in weighted spaces.} Let $f\in L_2(X;C)$ and $u=(\Delta_1-i)^{-1}f$. Since $\Delta_1$ is self-adjoint, we have 
$$\|u\|_{L_2(X;C)}\le \|f\|_{L_2(X;C)}.$$ In view of (\ref{H l smooth incr}), we also have 
$$\|\chi_k  u^{(k)}\|_{H^{2}}\le c(\|u\|_{L_2(X;C)}+\|f\|_{L_2(X;C)}), \qquad (k>2g-2).$$ 
Since $u\in {\rm Dom}\Delta_1\subset {\rm Dom}\Delta_{min}^*$, formula (\ref{emb to dom adj}) implies $y=\chi_k u^{(k)}\in H^2_1$ for $k\le 2g-2$. Also, $A_{1,2}y=w$, where 
$$w(r_k,\varphi_k)=-4r^2(\chi_k f(z_k)+iu(z)+[\Delta,\chi_k]u(z))$$ and $w\in H^0_{-1}$. Applying formula (\ref{ellyptic asymptotics in cone}) with $\nu=1$ and $\nu'=\nu_0\in (3/4,1)$ and taking into account that the strip $\Im\mu\in [\nu_0,1]$ contains no poles (\ref{OP resolvent poles}), we obtain $y=y'=A_{\nu_0,2}^{-1}w\in H^2_{\nu_0}$ and 
$$\|y\|_{H^2_{\nu_0}}\le c\|w\|_{H^0_2}\le c(\|u\|_{L_2(X;C)}+\|f\|_{L_2(X;C)})$$
due to (\ref{w spaces solution estim}). Combining the above estimates, we obtain the inequality
$$\|u\|_{\mathcal{H}_{\nu_0}^2}\le c\|f\|_{L_2(X;C)},$$
which means that $(\Delta_1-i)^{-1}$ acts continuously from $\mathcal{H}_1^0=L_2(X;C)$ onto $\mathcal{H}_{\nu_0}^2$.

\

\paragraph{\it Compactness of embeddings of weighted spaces.} Let $l'>l$ and $\nu'<\nu$. Now, we prove that the embedding $\mathcal{H}_{\nu'}^{l'}\subset\mathcal{H}_{\nu}^{l}$ is compact. Let $\mathfrak{B}_{\nu'}^{l'}$ is the unit ball in $\mathcal{H}_{\nu'}^{l'}$. To prove the claim, it is sufficient to show that, for any $\varepsilon>0$, there is a finite $\varepsilon$-net for $\mathfrak{B}_{\nu'}^{l'}$ in $\mathcal{H}_{\nu}^{l}$.

Let $\varepsilon>0$. First, each $v\in\mathfrak{B}_{\nu'}^{l'}$ can be represented as 
$$v=\sum_{k=2g-2}^{K}\chi_k v+\sum_{k=1}^{2g-2}\big(\chi_k (1-\chi_\epsilon) v+\chi_k\chi_\epsilon v\big).$$
In view of definition (\ref{weighted norm}), we have
$$\|\chi_k\chi_{k,\epsilon} v^{(z_k)}\|_{H_{\nu}^{l}}\le c\epsilon^{\nu-\nu'}\|\chi_k\chi_{k,\epsilon} v^{(z_k)}\|_{H_{\nu'}^{l'}}$$
so one can choose $\epsilon>0$ sufficiently small that $\|\sum_{k=1}^{2g-2}\chi_k\chi_{k,\epsilon} v\|_{\mathcal{H}_{\nu}^{l}}\le\varepsilon/3$. In what follows, we assume that $\epsilon$ is fixed in such a way.

In view of definition (\ref{weighted spaces}), for any $\nu'',l''$, the norm $\|u\|_{\mathcal{H}_{\nu''}^{l''}}$ of any section $u$ whose support is contained in $U_k$ and separated by a fixed distance from each conical point, is equivalent to $\|u^{(x_k)}\|_{H^{l''}(B)}$, where $B$ is a sufficiently large fixed ball in $\mathbb{C}$. Note that $u=\chi_k v$ or $u=\chi_k (1-\chi_\epsilon) v$ with $v\in\mathfrak{B}_{\nu'}^{l'}$ are examples of such sections, and there is a constant $R$ such that $\|u^{(k)}\|_{H^{l'}(B)}\le \|v\|_{\mathcal{H}_{\nu'}^{l'}}=R$ for them. 

Let $\mathfrak{B}^{l'}$ be the ball of radius $R$ in $H^{l'}(B)$. Due to the Rellich–Kondrachov theorem, the embedding $H^{l'}(B)\subset H^{l}(B)$ is compact. So, for arbitrary $\delta>0$, one can choose the finite $\delta$-net $a_1,\dots,a_n$ for $\mathfrak{B}^{l'}$ in $H^{l}(B)$. Then for any $u=\chi_k v$ or $u=\chi_k (1-\chi_\epsilon) v$ with $v\in\mathfrak{B}_{\nu'}^{l'}$, there exists $a_s$ such that 
\begin{equation}
\label{delta net local}
\|u^{(x_k)}-a_s\|_{H^{l}(B)}<\delta.
\end{equation}
Introduce the section $a_{s,k}$ supported in $U_k$ by the rule $a_{s,k}(x_k)=\tilde{\chi}_k a_s(x_k)$ if $k>2g-2$ and by the rule $a_{s,k}(x_k)=\tilde{\chi}_k (1-\chi_\epsilon/C)a_s(x_k)$ ($C$ is sufficiently large) if $k\le 2g-2$. Let $u=\chi_k v$ or $u=\chi_k (1-\chi_\epsilon) v$ with $v\in\mathfrak{B}_{\nu'}^{l'}$; then $\tilde{\chi}_k u=u$ or $\tilde{\chi}_k(1-\chi_\epsilon/C)u=u$. Now, estimate (\ref{delta net local}) implies
$$\|u-a_{s,k}\|_{\mathcal{H}_{\nu}^{l}}\le c\|u^{(x_k)}-a_s\|_{H^{l}(B)}\le c\delta.$$
Choose $\delta=\varepsilon/(13cg)$. As a corollary of the facts above, the sums $\sum_{k=1}^K a_{s_k,k}$ with $a_{s_k}=1,\dots,n$ provide the finite $\varepsilon$-net for $\mathfrak{B}_{\nu'}^{l'}$ in $\mathcal{H}_{\nu}^{l}$. By this, we have proved the compactness of the embedding $\mathcal{H}_{\nu'}^{l'}\subset\mathcal{H}_{\nu}^{l}$. 

\subsection{Self-adjoint extensions.}
In this subsection, basic properties of the self-adjoint extensions $\Delta_F$, $\Delta_S$, $\Delta_h$ of $\Delta_{min}$ are described.
\subsubsection{Friedrichs extension.}
\label{par Friedrichs}
In this paragraph, we prove Proposition \ref{prop Friedrichs} describing the domain of the Friedrichs extension $\Delta_F$ of $\Delta_{min}$. Also, we prove that $\Delta_F$ is positive defined and admits the representation 
\begin{equation}
\label{delta f via over D}
\Delta_F=\overline{\mathcal{D}}_\omega^*\overline{\mathcal{D}}_\omega,
\end{equation}
where $\overline{\mathcal{D}}_\omega$ is the closure (in $L_2(X;C)$) of the operator
\begin{equation}
\label{over D}
\overline{\mathcal{D}}_\omega u=\overline{\omega}^{-1}\overline{\partial}u \qquad (u\in C_c^\infty(\dot{X};C)).
\end{equation}

\

\paragraph{\it Quadratic form.} First, let us recall the construction of the Friedrichs extension. Introduce the quadratic form $a_{\mu,\cdot}[f]:=((\Delta_{min}-\mu)f,f)$ with the domain ${\rm Dom}a_{\mu,\cdot}={\rm Dom}\Delta_{min}$. In view of the Green formula 
\begin{align}
\label{half Green formula}
(\Delta f,f)_{L_2(U;C)}=\|\overline{\partial} f\|^2_{L_2(U;C\otimes\overline{K})}+\frac{1}{2i}\int_{\partial U}\overline{f}h\overline{\partial}fd\overline{z}
\end{align}
with $U=X_{\epsilon}(P_1,\dots,P_{2g-2})$, $\epsilon\to 0$, we have
$$a_{\mu,\cdot}[f]=\|\overline{\partial} f\|^2_{L_2(X;C\otimes\overline{K})}-\mu \|f\|^2_{L_2(X;C)},$$
and $a_{\mu,\cdot}$ is positive definite for $\mu<0$. Since the multiplication operator $\overline{\omega}^{-1}: \ L_2(X;C\otimes\overline{K})\mapsto L_2(X;C)$ is unitary, one can represent $a_{\mu,\cdot}$ as
$$a_{\mu,\cdot}[f]=\|\overline{\mathcal{D}}_\omega f\|^2_{L_2(X;C)}-\mu \|f\|^2_{L_2(X;C)}.$$
The form $a_{\mu,\cdot}$ admits closure denoted by $a_{\mu,F}$. The  Friedrichs extension $\Delta_F$ of $\Delta_{min}$ is defined as $\Delta_F=A+\mu$, where $A$ is the (unique) self-adjoint operator corresponding to the form $a_{\mu,\cdot}$ (i.e., $a_{\mu,F}[f]=(Af,f)_{L_2(U;C)}$ on ${\rm Dom}A\subset {\rm Dom}a_{\mu,F}$). Note that $\Delta_F$ is independent of $\mu<0$ and the exact lower bounds of the operators $\Delta_{min},\Delta_F$ coincide.

\

\paragraph{\it Domain of $\Delta_F$.} Let us prove Proposition \ref{prop Friedrichs}. Denote by $\Delta_{F'}$ the self-adjoint extension of $\Delta_{min}$ defined on the functions with asymptotics (\ref{asymptotics Friedrichs}). To prove that $\Delta_{F'}=\Delta_F$, it is sufficient to show that ${\rm Dom}\Delta_{F'}\subset {\rm Dom}a_{\mu,F}$ and $a_{\mu,F}[f]=((\Delta_{F'}-\mu)f,f)_{L_2(U;C)}$ on ${\rm Dom}\Delta_{F'}$. 

Recall that each $u\in {\rm Dom}\Delta_{F'}$ is the sum $u=u_0+\tilde{u}$, where $u_0$ is a linear combination of the sections
\begin{equation}
\label{Fried doma}
\chi_k f_{k,m,s} \qquad (k=1,\dots,2g-2, \quad (m,s)=(-1,-),(0,-),(1,+),(2,+))
\end{equation}
while $\tilde{u}\in \mathcal{H}^2_{-1}$ due to (\ref{emb to dom}). By definitions of $\chi_k f_{k,m,s}$, we have 
\begin{equation}
\label{Fried est near cp}
\partial_{\varphi_k}^p (r_k\partial_{r_k})^qu_0^{(k)}=O(r_k^{1/4-p-q}), \qquad \Delta u_0=0
\end{equation}
near each $P_k$. Hence, $u_0,u\in \mathcal{H}^2_{0}$ in view of definitions (\ref{weighted spaces}), (\ref{weighted norm}). Since $C_c^\infty(\dot{X};C)$ is dense in $\mathcal{H}^2_{-1}$, there is a sequence $\{\tilde{u}_n\}\subset C_c^\infty(\dot{X};C)$ converging to $\tilde{u}$ in $\mathcal{H}^2_{-1}$. Denote $u_n=u_0+\tilde{u}_n$. Then $u_n\to u$ in $\mathcal{H}^2_{0}$ and $L_2(X;C)$. Also, since $\Delta u_0=0$ near conical points, we have $\Delta u_n\to \Delta_{F'}u$ in $L_2(X;C)$.

Next, let us substitute $U=X_{\epsilon}(P_1,\dots,P_{2g-2})$, $\epsilon\to 0$ and $f=u_n$ into Green formula (\ref{half Green formula}). In view of (\ref{Fried est near cp}), we obtain $(\Delta u_n,u_n)_{L_2(X;C)}=\|\overline{\mathcal{D}}_\omega u_n\|^2_{L_2(X;C)}$. Hence, 
\begin{equation}
\label{half Green for extensions}
(\Delta u,u)_{L_2(X;C)}=\|\overline{\mathcal{D}}_\omega u\|^2_{L_2(X;C)}.
\end{equation}
is valid due to the above convergences. 

Finally, due to definitions (\ref{weighted spaces}), (\ref{weighted norm}), the convergence $u_n\to u$ in $\mathcal{H}^2_{0}$ implies that $\overline{\mathcal{D}}_\omega u_n\to \overline{\mathcal{D}}_\omega u$ in $L_2(X;C)$ and $a_{\mu,F}[u_n]\to \|\overline{\mathcal{D}}_\omega u\|^2_{L_2(X;C)}-\mu \|u\|^2_{L_2(X;C)}$. Hence, 
\begin{equation}
\label{fried proved}
u\in {\rm Dom}a_{\mu,F}, \quad a_{\mu,F}[u]=\|\overline{\mathcal{D}}_\omega u\|^2_{L_2(X;C\otimes\overline{K})}-\mu \|u\|^2_{L_2(X;C)}=((\Delta_{F'}-\mu)u,u)_{L_2(X;C)}
\end{equation}
for any $u\in {\rm Dom}\Delta_{F'}$ and $\mu\in\mathbb{C}$. Hence, $\Delta_{F'}=\Delta_F$ and Proposition \ref{prop Friedrichs} is proved. 

In addition, formula (\ref{fried proved}) implies that $\Delta_F$ is non-negative and each $u\in {\rm Ker}\Delta_F$ is holomorphic outside conical points. Also, each conical point $P_k$ is a zero of $u$ due to (\ref{asymptotics Friedrichs}). Hence, $u$ has at least $2g-2$ zeroes on $X$. Since the degree of the divisor class of $C$ is equal to $g-1\ge 1$, the latter is possible only if $u=0$ on $X$. So, ${\rm Ker}\Delta_F=0$ and $\Delta_F$ is positive definite. Note that this fact is independent on the assumption $h^0(C)=0$.

\

\paragraph{\it Connection between $\Delta_F$ and $\overline{\mathcal{D}}_\omega$, $\overline{\mathcal{D}}_\omega^*$.} Since
\begin{align*}
%\label{1st order operator adjoint 1}
\begin{split}
(\overline{\mathcal{D}}_\omega u,v)_{L_2(X,C)}=&\int_{X}\overline{\omega}^{-1}\overline{\partial}u\cdot h\overline{v}\rho^{-2}dzd\overline{z}=\\
=&\int_X\overline{\partial}(\overline{\omega}^{-1}uh\overline{v}\rho^{-2})dzd\overline{z}-\int_X u\overline{\partial}(\overline{\omega}^{-1}h\overline{v}\rho^{-2})dzd\overline{z}=\\
=0+&\int_X uh \overline{\big[-(\rho^{-2}h)^{-1}\partial (\rho^{-2}h\omega^{-1}v)\big]}dS \qquad (u,v\in C_c^\infty(\dot{X};C)),
\end{split}
\end{align*}
the operator $\overline{\mathcal{D}}_\omega^*$ satisfies
\begin{equation}
\label{1st order operator adjoint}
\overline{\mathcal{D}}_\omega^* u=-(\rho^{-2}h)^{-1}\partial \big(\rho^{-2}h\omega^{-1}u\big)=-|\omega|^{-1}\partial\big(|\omega|\omega^{-1}u\big) \quad (u\in C_c^\infty(\dot{X};C)).
\end{equation}
In particular, $\overline{\mathcal{D}}_\omega^*$ is densely defined and operator (\ref{over D}) admits closure (denoted by $\overline{\mathcal{D}}_\omega$) while the operators $\overline{\mathcal{D}}_\omega^*\overline{\mathcal{D}}_\omega$ and $\overline{\mathcal{D}}_\omega\overline{\mathcal{D}}_\omega^*$ are self-adjoint. In addition,
$$\overline{\mathcal{D}}_\omega^*\overline{\mathcal{D}}_\omega u=-|\omega|^{-1}\partial\big(|\omega|\omega^{-1}\overline{\omega}^{-1}\overline{\partial}u\big)=\Delta_{min} u\qquad \text{ for } u\in C_c^\infty(\dot{X};C).$$ 
Taking the closure of both parts yields $\Delta_{min}\subset\overline{\mathcal{D}}_\omega^*\overline{\mathcal{D}}_\omega$.

Next, let $u\in {\rm Dom}\Delta_F$. As shown after (\ref{Fried doma}), there us a sequence $\{u_n\}\subset C_c^\infty(\dot{X};C)$ converging to $u$ in $\mathcal{H}^2_0$ and in $L_2(X;C)$ and such that $\overline{\mathcal{D}}_\omega^*\overline{\mathcal{D}}_\omega u_n=\Delta u_n\to \Delta_{F}u$ in $L_2(X;C)$. Hence, $u\in {\rm Dom}\overline{\mathcal{D}}_\omega^*\overline{\mathcal{D}}_\omega$ and $\overline{\mathcal{D}}_\omega^*\overline{\mathcal{D}}_\omega u=\Delta_F u$. Thus, $\Delta_F\subset \overline{\mathcal{D}}_\omega^*\overline{\mathcal{D}}_\omega$ and (since both operators are self-adjoint) $\overline{\mathcal{D}}_\omega^*\overline{\mathcal{D}}_\omega=\Delta_F$. So, formula (\ref{delta f via over D}) is proved. 

\subsubsection{Szeg\"o extension.}
\label{par Szego}
In this paragraph, we prove Proposition \ref{prop Szego ext}. First, recall that, under the assumption $h^0(C)=0$, the Szeg\"o kernel $(z,z')\mapsto \mathcal{S}(z,z')$ (the integral kernel of $(-\frac{1}{\pi}\overline{\partial})^{-1}$) is a section of $C$ in both variables, holomorphic outside the diagonal $z=z'$, antisymmetric $\mathcal{S}(z,z')=-\mathcal{S}(z',z)$, and admitting the asymptotics 
$$\mathcal{S}(z,z')=\frac{1}{z-z'}+O(1)$$
near the diagonal $z=z'$ (see p.29, \cite{Fay}). Here $z,z'$ are understood as different values of the same (arbitrary) holomorphic local coordinate on $X$. In what follows, we assume that $z,z'$ are given by (\ref{z coordinate}).

Denote by $\Delta_{S}$ the self-adjoint extension of $\Delta_{min}$ defined on the functions with asymptotics (\ref{Szego domain}). Let $\mathcal{G}^{\mathcal{S}}$ be the Green function of $\Delta_{S}$. Recall that $(z,z')\mapsto\mathcal{G}^{\mathcal{S}}(z,z')$ is a section of $C$ ($\overline{C}$) in $z$ ($z'$) obeying $\Delta \mathcal{G}^{\mathcal{S}}(\cdot,z')=0$ outside the diagonal $z=z'$ and conical points. Also, it is symmetric $\mathcal{G}^{\mathcal{S}}(z,z')=\overline{\mathcal{G}^{\mathcal{S}}(z',z)}$ and obeys
\begin{equation}
\label{Green f near diag}
\mathcal{G}^{\mathcal{S}}(z,z')=-\frac{2}{\pi}{\rm log}(|z'-z|)+\dots
\end{equation} 
outside conical points, where dots denote a smooth function of local coordinates $z,z'$. In view of Proposition \ref{prop elliptic asymptotics} and formula (\ref{Szego domain}), we have
\begin{equation}
\label{Szego Green asymp 1}
\begin{split}
\mathcal{G}(x_k,z')=c_{k,0,+}(z')+&c_{k,1,+}(z')x_k+c_{k,2,+}(z')x_k^2+\\
+&(-2/3\pi)c_{k,-1,-}(z')\overline{x_k}|x_k|+O(|x_k|^3)
\end{split}
\end{equation}
near each $P_k$; asymptotics (\ref{Szego Green asymp 1}) admits differentiation which means that the derivatives of the remainder $\tilde{\mathcal{G}}(x_k,z')=O(|x_k|^3)$ satisfy
$$\partial_{z'}^l\tilde{\mathcal{G}}(x_k,z'),\partial_{\overline{z'}}^l\tilde{\mathcal{G}}(x_k,z')=O(|x_k|^{3-l}), \quad \partial_{x_k}^l\tilde{\mathcal{G}}(x_k,z'),\partial_{\overline{x_k}}^l\tilde{\mathcal{G}}(x_k,z')=O(|x_k|^{3-l}) \qquad (l=1,2,\dots).$$

Denote 
$$\tilde{\mathcal{S}}(z,z')=-\pi h(z')\partial_{z'}\mathcal{G}^{\mathcal{S}}(z,z').$$
Then 
$$-\frac{1}{\pi}\partial_{\overline{z'}}\tilde{\mathcal{S}}(z,z')=\overline{\partial_{z'}h(z')\partial_{\overline{z'}}\mathcal{G}^{\mathcal{S}}(z',z)}=-\overline{\rho^{-2}(z')h(z')\Delta_{z'}\mathcal{G}^{\mathcal{S}}(z',z)}=0,$$
i.e., $z'\mapsto \tilde{\mathcal{S}}(z,z')$ is holomorphic in $z'$ outside the diagonal $z'=z$ and conical points. Next, formula (\ref{Green f near diag}) implies 
$$\tilde{\mathcal{S}}(z,z')=-\pi\partial_{z'}(-\frac{2}{\pi}{\rm log}(|z'-z|)+O(1))=\frac{1}{z'-z}+O(1).$$
Finally, due to the symmetry $\mathcal{G}^{\mathcal{S}}(z,z')=\overline{\mathcal{G}^{\mathcal{S}}(z',z)}$, formula (\ref{Szego Green asymp 1}) can be rewritten as
\begin{align*}
\mathcal{G}^{\mathcal{S}}(z',x_k)=\overline{c^0_{k,0,+}(z')}+\overline{c^0_{k,1,+}(z')x_k}+&\overline{c^0_{k,2,+}(z')x_k^2}+\\
+&(-2/3\pi)\overline{c^0_{k,-1,-}(z')}x_k|x_k|+O(|x_k|^3).
\end{align*}
Differentiation of both sides yields
$$\tilde{\mathcal{S}}(z',x_k)=(2/3)\overline{c^0_{k,-1,-}(z')}|x_k|^{-1}\partial_{x_k}(|x_k|x_k)+O(|x_k|)=\overline{c^0_{k,-1,-}(z')}+O(|x_k|).$$
Therefore $\tilde{\mathcal{S}}(z',\cdot)$ is holomorphic at each conical point. In view of the facts above, $[\tilde{\mathcal{S}}-\mathcal{S}](z',\cdot)$ is a holomorphic section of $C$ for each $z'$. Since $h^0(C)=0$, we have $\tilde{\mathcal{S}}=\mathcal{S}$ and formula (\ref{Green to Szego}) is valid.

Since $\mathcal{S}$ is the integral kernel of $(-\frac{1}{\pi}\overline{\partial})^{-1}$, the formula
\begin{equation}
\label{Szego is kernel}
f(z)=\int_{X}\mathcal{S}(z,z')\big[-\frac{1}{\pi}\overline{\partial}f\big](z')\frac{d\overline{z'}dz'}{2i}
\end{equation}
is valid on smooth sections $f$ of $C$. Namely, (\ref{Szego is kernel}) follows from the Stokes theorem (applied to the domain obtained from $X$ by removing small neighborhoods of the singularities of the integrand) and the near-diagonal asymptotics of $\mathcal{S}$. The same argument shows that equality (\ref{Szego is kernel}) is valid for sections $f$ with logarithmic singularities. Thus, substituting $f=\mathcal{G}^{\mathcal{S}}(\cdot,z'')$ into (\ref{Szego is kernel}) and taking into account that $\rho=h$, we obtain
\begin{align*}
\mathcal{G}^{\mathcal{S}}(z,z'')=\int\limits_{X}\mathcal{S}(z,z')\big[-\frac{1}{\pi}\overline{\partial_{z'}\mathcal{G}^{\mathcal{S}}(z'',z')}\big]\frac{d\overline{z'}dz'}{2i}=\\
=-\frac{1}{\pi^2}\int\limits_{X}\mathcal{S}(z,z')h(z')\overline{\mathcal{S}(z',z'')}dS(z')
\end{align*}
Therefore, formula (\ref{Szego to Green}) is valid. Thus, $\Delta_{S}$ is the Szeg\"o extension of $\Delta_{min}$. Proposition \ref{prop Szego ext} is proved.

\

\paragraph{\it Positivity of $\Delta_S$.} In view of (\ref{Szego domain}), each $u\in {\rm Dom}\Delta_S$ admits representation $u=u_0+\tilde{u}$, where $\partial_{\overline{x_k}}u_0(x_k)=-\pi c_{k,-1,-}|x_k|$ near $P_k$ while $\tilde{u}\in \mathcal{H}^2_{-1}$ due to (\ref{emb to dom}). Then formula (\ref{half Green formula}) yields (\ref{half Green for extensions}) with $\Delta=\Delta_S$. In particular, $\Delta_S$ is non-negative, and each $u\in {\rm Ker}\Delta_S$ is holomorphic outside conical points. Due to (\ref{Szego domain}), $u(x_k)$ is bounded whence $u$ is holomorphic at conical points. Since $h^0(C)=0$, we obtain $u=0$. So, the Szeg\"o extension $\Delta_S$ is positive.

\subsubsection{Holomorphic extension.}
\label{par holom}
In this paragraph, we prove Proposition \ref{prop isospectrality of F and h} on almost-isospectrality of holomorphic and Friedrichs extensions of $\Delta_{min}$. In addition, we describe the Bergman kernel of $\Delta_h$. 

\

\paragraph{\it Almost-isospectrality of $\Delta_F$ and $\Delta_h$.} Let us prove Proposition \ref{prop isospectrality of F and h}. To this end, we consider the self-adjoint operator $\overline{\mathcal{D}}_\omega\overline{\mathcal{D}}_\omega^*$ ($\overline{\mathcal{D}}_\omega$ is defined after (\ref{delta f via over D})). In view of (\ref{over D}) and (\ref{1st order operator adjoint}), we have
\begin{equation}
\label{D D adjoint}
\overline{\mathcal{D}}_\omega\overline{\mathcal{D}}_\omega^*=-\overline{\omega}^{-1}\overline{\partial}\big(|\omega|^{-1}\partial\big(|\omega|\omega^{-1}\dots\big)\big)=U_\omega^{-1}\overline{\iota}^{-1} \Delta_0 \overline{\iota} U_\omega \quad \text{ on } C_c^\infty(\dot{X};C),
\end{equation}
where 
$$U_\omega: \ u\mapsto |\omega|\omega^{-1}u$$
is the unitary operator acting from $L_2(X;C)$ onto $L_2(X;\overline{C})$ and
$$\overline{\iota}: u\mapsto \overline{u}$$
is anti-linear isometry from $L_2(X;\overline{C})$ onto $L_2(X;C)$. Since $\overline{\mathcal{D}}_\omega\overline{\mathcal{D}}_\omega^*$ is self-adjoint, so is the operator $\overline{\Delta}_{\odot}=U_\omega\overline{\mathcal{D}}_\omega\overline{\mathcal{D}}_\omega^*U_\omega^{-1}$ in $L_2(X;\overline{C})$. Since $\overline{\iota}$ is anti-linear, the formulas $\Delta_\odot:=\overline{\iota} \overline{\Delta}_{\odot}\overline{\iota}^{-1}$ and ${\rm dom}\Delta_\odot=\overline{\iota} U_\omega{\rm dom}(\overline{\mathcal{D}}_\omega\overline{\mathcal{D}}_\omega^*)$ define the linear unbounded operator $\Delta_\odot$ in $L_2(X;C)$. The equality
\begin{align*}
(u,\overline{\Delta}_{\odot}^*v)_{L_2(X;C)}=(\Delta_\odot u,v)_{L_2(X;C)}=\int_{X}\overline{\overline{\Delta}_{\odot}\overline{u}}h\overline{v}dS=\overline{(\overline{\Delta}_{\odot}\overline{u},\overline{v})_{L_2(X;\overline{C})}}=\\=\overline{(\overline{u},\overline{\Delta}_{\odot}\overline{v})_{L_2(X;\overline{C})}}=\int_{X}uh\overline{\overline{\Delta}_{\odot}\overline{v}}dS=(u,\overline{\Delta}_{\odot}v)_{L_2(X;C)}
\end{align*}
shows that $\Delta_\odot$ is self-adjoint. Due to (\ref{D D adjoint}), $\Delta_\odot$ is a self-adjoint extension of $\Delta_{min}$. 

According to Proposition \ref{discreteness spectrum}, the spectrum of $\Delta_\odot$ is discrete. We have
\begin{align*}
\Delta_F u=\lambda u \ \Rightarrow \ \overline{\mathcal{D}}_\omega\overline{\mathcal{D}}_\omega^*(\overline{\mathcal{D}}_\omega u)&=\overline{\mathcal{D}}_\omega \Delta_F u=\lambda(\overline{\mathcal{D}}_\omega u) \Leftrightarrow \\ &\Leftrightarrow \ \Delta_{\odot}(\overline{\iota} U_\omega\overline{\mathcal{D}}_\omega u)=\lambda(\overline{\iota} U_\omega\overline{\mathcal{D}}_\omega u)
\end{align*}
Suppose that $\lambda\ne 0$; since $\overline{\iota} U_\omega\overline{\mathcal{D}}_\omega u=0$ implies $\overline{\mathcal{D}}_\omega u=0$ and $\Delta_F u=\overline{\mathcal{D}}_\omega^*\overline{\mathcal{D}}_\omega u=0$, the last formula means that ${\rm dim}\,{\rm Ker}(\Delta_{\odot}-\lambda)\ge {\rm dim}\,{\rm Ker}(\Delta_F-\lambda)$. Similarly, we have 
\begin{align*}
\Delta_{\odot}v=\lambda v \ \Leftrightarrow \ \overline{\mathcal{D}}_\omega\overline{\mathcal{D}}_\omega^*(U_\omega^{-1}\overline{\iota}^{-1}v)&=\lambda(U_\omega^{-1}\overline{\iota}^{-1}v)\Rightarrow \\ 
\Rightarrow& \Delta_F(\overline{\mathcal{D}}_\omega^* U_\omega^{-1}\overline{\iota}^{-1}v)=\lambda(\overline{\mathcal{D}}_\omega^* U_\omega^{-1}\overline{\iota}^{-1}v).
\end{align*}
Since $\overline{\mathcal{D}}_\omega^* U_\omega^{-1}\overline{\iota}^{-1}v=0$ implies $\Delta_{\odot}v=(\overline{\iota}U_\omega\overline{\mathcal{D}}_\omega\overline{\mathcal{D}}_\omega^*U_\omega^{-1}\overline{\iota}^{-1})v=0$, the last formula for non-zero $\lambda$ means that ${\rm dim}\,{\rm Ker}(\Delta_F-\lambda)\ge {\rm dim}\,{\rm Ker}(\Delta_{\odot}-\lambda)$. Thus, we have arrived to the equality ${\rm dim}\,{\rm Ker}(\Delta_F-\lambda)={\rm dim}\,{\rm Ker}(\Delta_{\odot}-\lambda)$ for all $\lambda\ne 0$ which means that the operators $\Delta_F$ and $\Delta_{\odot}$ are almost isospectral and the anti-linear map
$$u\mapsto \overline{\iota} U_\omega\overline{\mathcal{D}}_\omega u$$ 
provides the correspondence between eigensections of $\Delta_F$ and eigensections of $\Delta_\odot$ with non-zero eigenvalues.

To prove Proposition \ref{prop isospectrality of F and h} it remains to show that $\Delta_\odot=\Delta_h$. Let $(\lambda,u)$ be any eigenpair (i.e., a pair of eigenvalue and eigensection) of $\Delta_F$ obeying $\lambda\ne 0$. Then $u$ is smooth outside of the conical points. In view of Proposition \ref{prop elliptic asymptotics} and formula (\ref{asymptotics Friedrichs}), the asymptotics 
\begin{align*}
u(x_k)=&(-2/\pi)c_{k,0,-}|x_k|+c_{k,1,+}x_k+c_{k,2,+}x_k^2+\\+&(-2/3\pi)c_{k,-1,-}\overline{x_k}|x_k|+c_{k,3,+}x_k^3+(-2/5\pi)c_{k,-2,-}\overline{x_k}^2|x_k|+\\
+&c_{k,4,+}x_k^4(-2/7\pi)c_{k,-3,-}\overline{x_k}^3|x_k|+O(|x_k|^5), \qquad x_k\to 0
\end{align*}
holds and admits differentiation (i.e., the remainder $\tilde{u}=O(|x_k|^5)$ obeys $\partial^l_{x_k}\tilde{u},\partial^l_{\overline{x_k}}\tilde{u}=O(|x_k|^{5-l})$). Hence,
\begin{align*}
-\pi[\overline{\mathcal{D}}_\omega u](x_k)=-\pi\overline{x_k}^{-1}\partial_{\overline{x_k}}u(x_k)=c_{k,0,-}\overline{x_k}^{-2}|x_k|+c_{k,-1,-}\overline{x_k}^{-1}|x_k|+\\
+c_{k,-2,-}|x_k|+c_{k,-3,-}\overline{x_k}|x_k|+O(|x_k|^3), \qquad x_k\to 0
\end{align*}
and
\begin{align*}
-\pi[\overline{\iota} U_\omega\overline{\mathcal{D}}_\omega u](x_k)=&\overline{x_k^{-1}|x_k|\big(-\pi\overline{x_k}^{-1}\partial_{\overline{x_k}}u(x_k))}=\\
=c_{k,0,-}x_k^{-1}+&c_{k,-1,-}+c_{k,-2,-}x_k+c_{k,-3,-}x_k^2+O(|x_k|^3), \qquad x_k\to 0.
\end{align*}
Thus, $\overline{\iota} U_\omega\overline{\mathcal{D}}_\omega u\in {\rm Dom}\Delta_h$ due to (\ref{holomorphic asymp}). Therefore, each eigenpair $(\overline{\iota} U_\omega\overline{\mathcal{D}}_\omega \lambda\ne 0,u)$ of $\Delta_\odot$ is an eigenpair of $\Delta_h$. In other words, we have $\Delta_\odot(I-P)=\Delta_h(I-P)$, where $P$ is a projection on the (finite-dimensional) subspace ${\rm Ker}\Delta_\odot$. At the same time, $\Delta_\odot=\Delta_h$ on set $C_c^\infty(\dot{X};C)$. Since $C_c^\infty(\dot{X};C)$ is dense in $L_2(X;C)$, we obtain $\Delta_\odot\equiv\Delta_h$. Proposition \ref{prop isospectrality of F and h} is proved.

\

\paragraph{\it Bergman kernel of $\Delta_h$.} In view of (\ref{holomorphic asymp}), each $u\in {\rm Dom}\Delta_h$ admits representation $u=u_0+\tilde{u}$, where $\partial_{\overline{x_k}}u_0(x_k)=0$ near $P_k$ while $\tilde{u}\in \mathcal{H}^2_{-1}$ due to (\ref{emb to dom}). Then formula (\ref{half Green formula}) yields (\ref{half Green for extensions}) with $\Delta=\Delta_h$. In particular, $\Delta_h$ is non-negative and each $u\in {\rm Ker}\Delta_h$ is holomorphic outside conical points. In view of (\ref{holomorphic asymp}), $u\in {\rm Ker}\Delta_h$ may have simple poles at conical points. Due to the assumption $h^0(C)=0$, the Szeg\"o kernels $\mathcal{S}(P_j,\cdot)$ ($j=1,\dots,2g-2$) provide a basis in ${\rm Ker}\Delta_h$; in particular, ${\rm dim}{\rm Ker}\Delta_h=2g-2$. Then the Bergman kernel of $\Delta_h$ is given by
\begin{equation*}
%\label{Bergman kernel holom ext}
\mathscr{B}^h(z,z')=-\sum_{kj}\mathcal{S}(z,P_k)\mathfrak{S}^{-1}_{kj}\overline{\mathcal{S}(P_j,z')},
\end{equation*}
where $\mathfrak{S}$ is given by (\ref{Bergman matrix}). 

\subsection{Heat kernel}
\label{ssec heat}
In this subsection, we construct the parametrix $\mathscr{H}_0(\cdot,\cdot,t)$ for the fundamental solution $\mathscr{H}(\cdot,\cdot,t)$ to heat equation $(\partial_t+\Delta_F)u=0$. Next, we derive the short-time asymptotics of $\mathscr{H}(\cdot,\cdot,t)$. In particular, we obtain the short-time asymptotics of the heat trace $t\mapsto K(t|\Delta_F)$ and prove that the zeta function $s\mapsto\zeta(s|\Delta_F)$ admits analytic continuation into the neighborhood of $s=0$. Due to this fact, the determinant of $\Delta_F$ is well-defined.

\subsubsection{Parametrix for heat kernel}

\

\paragraph{\it Local solution outside conical points.} Outside conical points, coordinates (\ref{z coordinate}) are well-defined and $\Delta$ is of the form (\ref{elipticity}). Hence, the local fundamental solution to the heat  equation is given by 
$$\mathscr{H}_{\mathbb{C}}(z,z',t)=\frac{1}{\pi t}e^{-|z-z'|^2/t},$$
where $z,z'$ are values of coordinate (\ref{z coordinate}).

\

\paragraph{\it Local solution near conical points.} Near $P_k$, the heat equation $(\partial_t+\Delta_F)u=0$ takes the form 
\begin{equation}
\label{heat eq in cone}
\big(4\partial_t-r^{-2}((r\partial_{r})^2+\partial_\varphi^2)\big)v(r,\varphi)=0,
\end{equation}
where $r=r_k$, $\varphi=\varphi_k$. In addition, the representative %local trivialization 
$v=u^{(k)}$ of the section $u$ of $C$ in coordinate $z_k$ obeys (\ref{antiperiodicity}), i.e., it is $4\pi$-antiperiodic in $\varphi=\varphi_k$.

The fundamental solution to equation (\ref{heat eq in cone}), which is $\alpha-$periodic in $\varphi,\varphi'$, is given by 
\begin{equation}
\label{Dowker heat kernel}
\mathscr{H}_{\circlearrowleft}(r,\varphi,r',\varphi',t | \alpha)=\frac{1}{2\pi\alpha it}\int\limits_{\mathcal{C}\cup(-\mathcal{C})}{\rm exp}\Big(-\frac{r^2-2rr'{\rm cos}\vartheta+r'^{2}}{t}\Big){\rm cot}\theta d\vartheta.
\end{equation}
(see \cite{Carslaw}, see also \cite{Dowker,KokKorot}). Here $\theta:=\pi\alpha^{-1}(\vartheta+\varphi-\varphi')$ and $\mathcal{C}$ is the contour in the semi-strip $\Re\vartheta\in (-\pi,\pi)$, $\Im\vartheta>0$ running from $\vartheta=i\infty+\pi$ to $\vartheta=i\infty-\pi$. We take the function
\begin{equation*}
\mathscr{H}_{\curvearrowleft}(r,\varphi,r',\varphi',t | \alpha)=\mathscr{H}_{\circlearrowleft}(r,\varphi,r',\varphi',t | 2\alpha)-\mathscr{H}_{\circlearrowleft}(r,\varphi+\alpha,r',\varphi',t | 2\alpha)
\end{equation*}
as an $\alpha$-anti-periodic fundamental solution to (\ref{heat eq in cone}).

Let us find the asymptotics of $\mathscr{H}_{\curvearrowleft}(r,\varphi,r',\varphi',t | \alpha)$ as $r\to 0$. From (\ref{Dowker heat kernel}), it follows that
\begin{equation}
\label{antiperiodic heat kernel}
\begin{split}
\mathscr{H}_{\curvearrowleft}(r,\varphi,r',\varphi',t | \alpha)&=\\
=\frac{1}{4\pi\alpha it}&\int\limits_{\mathcal{C}\cup(-\mathcal{C})}{\rm exp}\Big(-\frac{r^2-2rr'{\rm cos}\vartheta+r'^{2}}{t}\Big)({\rm cot}\frac{\theta}{2}+{\rm tan}\frac{\theta}{2})d\vartheta\\
=\frac{1}{\pi\alpha t}&\int\limits_{\mathcal{C}\cup(-\mathcal{C})}{\rm exp}\Big(-\frac{r^2-2rr'{\rm cos}\vartheta+r'^{2}}{t}\Big)\frac{a}{a^2-1}d\vartheta,
\end{split}
\end{equation}
where $a={\rm exp}(i\theta)$. Since
\begin{align*}
\vartheta&\in\mathcal{C} \quad \Rightarrow \quad |a|<1, \quad \frac{a}{a^2-1}=-\sum_{k\ge 0}a^{2k+1},\\
\vartheta&\in(-\mathcal{C}) \quad \Rightarrow \quad |a|>1, \quad \frac{a}{a^2-1}=\sum_{k\le 0}a^{2k-1},
\end{align*}
we have
\begin{align*}
-\pi\alpha t &e^{(r^2+r'^{2})/t}\mathscr{H}_{\curvearrowleft}(r,\varphi,r',\varphi',t | \alpha)=\sum_\pm (\pm)\sum_{\pm k\ge 0}\int\limits_{\pm\mathcal{C}}e^{\frac{2rr'{\rm cos}(\vartheta)}{t}}a^{2k\pm 1} d\vartheta=\\
=&\sum_\pm \sum_{\pm k\ge 0}e^{\frac{\pm i\pi(\varphi-\varphi')(2(\pm k)+1)}{\alpha}}\int\limits_{\pm\vartheta\in\mathcal{C}}e^{\frac{2rr'{\rm cos}(\pm\vartheta)}{t}}e^{\frac{i\pi(\pm\vartheta)(2(\pm k)+1)}{\alpha}}d(\pm\vartheta)=\\
=&2\sum_{k\ge 0}{\rm cos}\Big(\frac{\pi(\varphi-\varphi')(2k+1)}{\alpha}\Big)\int\limits_{\mathcal{C}}e^{\frac{2rr'{\rm cos}(\vartheta)}{t}}e^{\frac{i\pi\vartheta(2k+1)}{\alpha}}d\vartheta=[\vartheta=i\gamma]=\\
=&2i\sum_{k\ge 0}{\rm cos}\Big(\frac{\pi(\varphi-\varphi')(2k+1)}{\alpha}\Big)\int\limits_{-i\mathcal{C}}{\rm exp}\Big(\frac{2rr'}{t}{\rm cosh}\gamma-\frac{\pi(2k+1)}{\alpha}\gamma\Big)d\gamma.
\end{align*}
Taking into account the expression
$$\frac{1}{2\pi i}\int\limits_{-i\mathcal{C}}{\rm exp}\big(z{\rm cosh}s-\nu s\big)ds=I_\nu(z),$$
for the modified Bessel functions of the first kind, we finally obtain
\begin{equation}
\label{heat kernel asymptotics near vertex}
\begin{split}
\mathscr{H}_{\curvearrowleft}(r,\varphi&,r',\varphi',t | \alpha)=\\
=\frac{1}{\alpha t}&{\rm exp}\big(-\frac{r^2+r'^{2}}{t}\big)\sum_{k\ge 0}{\rm cos}\Big(\frac{\pi(\varphi-\varphi')(2k+1)}{\alpha}\Big)I_{\pi(2k+1)/\alpha}(\frac{2rr'}{t})=O(r^{\frac{\pi}{\alpha}})
\end{split}
\end{equation}
as $r\to 0$. Let us show that $\mathscr{H}_{\curvearrowleft}(r,\varphi,r',\varphi',t | \alpha)$ and all its derivatives decay exponentially as $t\to +0$ if $(r,\varphi{\rm \ mod \ }\alpha)$ is separated from $(r',\varphi'{\rm \ mod \ }\alpha)$. First, note that one can replace the integration contour $\mathcal{C}\cup(-\mathcal{C})$ in (\ref{antiperiodic heat kernel}) with the union of two lines $\pm l:=\{\vartheta=\pm(\pi-i\acute{\vartheta})\}_{\acute{\vartheta}\in\mathbb{R}}$ and small anti-clockwise circles $\circ(\vartheta_*)$ centered at the zeroes $\vartheta_*$ of functions $a\pm 1$. The exponent in the right hand side of (\ref{antiperiodic heat kernel}) decays as $t\to +0$ uniformly in $\vartheta\in (\pm l)$ since ${\rm cos}\vartheta=-{\rm cosh}\acute{\vartheta}<-1$ on $\pm l$. In addition, $|a/(a^2-1)|=O({\rm exp}(-\pi\alpha^{-1}|\acute{\vartheta}|))$ holds on the lines $\pm l$. Thus, the integrals over $\pm l$ in the right hand side of (\ref{antiperiodic heat kernel}) converge and decay exponentially as $t\to +0$. Similarly, the exponent in the right hand side of (\ref{antiperiodic heat kernel}) decays uniformly in $\vartheta\in \circ(\vartheta_*)$, where $\vartheta_*\ne 0$. Finally, the integral over $\circ(0)$ is present in the right hand side of (\ref{antiperiodic heat kernel}) if and only if $\varphi=\varphi'\,({\rm mod}\alpha)$. In the last case, this integral is equal to $e^{-(r-r')^{2}/4t}/\pi t$ due to the residue theorem. This implies the needed decay of $\mathscr{H}_{\curvearrowleft}(r,\varphi,r',\varphi',t | \alpha)$ as $t\to +0$ for $(r,\varphi{\rm \ mod \ }\alpha)$ separated from $(r',\varphi'{\rm \ mod \ }\alpha)$.

If $r=r'$, $\varphi=\varphi'$, then (\ref{antiperiodic heat kernel}) takes the form
$$\mathscr{H}_{\curvearrowleft}(r,\varphi,r,\varphi,t | \alpha)=\frac{1}{2\pi i\alpha t}\int\limits_l {\rm exp}\big(-\frac{4r^2{\rm sin}^2(\vartheta/2)}{t}\big)\frac{d\vartheta}{{\rm sin}(\pi\vartheta/\alpha)},$$
where the integration is over the union $l$ of $\pm l$ and $\circ(\vartheta_*)$-s. Therefore, 
\begin{equation}
\label{integration of conic heat kernel}
\begin{split}
\int\limits_0^R rdr\int\limits_0^\alpha d\varphi \mathscr{H}_{\curvearrowleft}(r,\varphi&,r,\varphi,t | \alpha)=\\
=\frac{1}{16\pi i}\int\limits_l &\frac{d\vartheta}{{\rm sin}^2(\vartheta/2){\rm sin}(\pi\vartheta/\alpha)}{\rm exp}\big(-\frac{4r^2{\rm sin}^2(\vartheta/2)}{t}\big)\Big|^0_{r=R}.
\end{split}
\end{equation}
Here the term with $r=0$ equals zero while the term with $r=R$ coincides with
\begin{align*}
\frac{-1}{16\pi i}\int\limits_{\circ(0)} \frac{d\vartheta}{{\rm sin}^2(\vartheta/2){\rm sin}(\pi\vartheta/\alpha)}&{\rm exp}\big(-\frac{4R^2{\rm sin}^2(\vartheta/2)}{t}\big)=\\
=-\frac{1}{8}\underset{\vartheta=0}{{\rm Res}}\Big(&\frac{{\rm exp}\big(-\frac{4R^2{\rm sin}^2(\vartheta/2)}{t}\big)}{{\rm sin}^2(\vartheta/2){\rm sin}(\pi\vartheta/\alpha)}\Big)=-\frac{1}{8}\Big(\frac{\alpha}{3\pi}+\frac{2\pi}{3\alpha}-\frac{\alpha}{\pi}\frac{4R^2}{t}\Big)
\end{align*}
up to a remainder that exponentially decays as $t\to +0$, $R\ge {\rm const}>0$. Thus, (\ref{integration of conic heat kernel}) implies 
\begin{align}
\label{heat kernel trace asymptotics conic}
\begin{split}
\int_0^R \chi(r)rdr\int_0^\alpha d\varphi \mathscr{H}_{\curvearrowleft}(r,\varphi,r,\varphi,t)&=\\
=\frac{1}{\pi t}\int_0^R \chi(r)rdr\int_0^\alpha &d\varphi-\frac{1}{8}\Big(\frac{\alpha}{3\pi}+\frac{2\pi}{3\alpha}\Big)+O(e^{-\epsilon /t}),
\end{split}
\end{align}
where $\chi$ is a smooth cut-off function on $\mathbb{R}$ equal to $1$ near zero and $\epsilon$ is some positive number.

\

\paragraph{\it Construction of parametrix.}
For $t>0$, let $(z,z')\mapsto\mathscr{H}_{0,k}(z,z',t)$ be the section of $C$ ($\overline{C}$) in $z$ ($z'$) which vanishes outside $U_k\times U_k$ and is given by
\begin{align}
\label{parametrix parts}
\mathscr{H}_{0,k}(z_k,z'_k,t)=\left\{\begin{array}{lr}
\tilde{\chi}_k(z_k)\chi_k(z'_k)\mathscr{H}_{\mathbb{C}}(z_k,z'_k,t) &\qquad (k>2g-2),\\
\tilde{\chi}_k(z_k)\chi_k(z'_k)\mathscr{H}_{\curvearrowleft}(r_k,\varphi_k,r'_k,\varphi'_k,t | 4\pi)
 &\qquad (k\le 2g-2)
\end{array}\right.
\end{align}
in $U_k\times U_k$ (here $z_k$, $z'_k$ denote different values of the same coordinate). As a parametrix for the fundamental solution to $(\partial_t+\Delta_F)u=0$ we take the sum 
\begin{equation*}
%\label{heat kernel parametrix}
\mathscr{H}_0=\sum_k\mathscr{H}_{0,k}.
\end{equation*}

For $t>0$ define 
$$\tilde{f}(\cdot,z',t):=(\partial_t+\Delta_{z})\mathscr{H}_0(\cdot,z',t).$$
Since each $\mathscr{H}_{0,k}$ in (\ref{parametrix parts}) is a fundamental solution to the heat equation in $U_k\times U_k$, one has $\tilde{f}=\sum_k\tilde{f_k}$, where $\tilde{f_k}(\cdot,\cdot,t)$ vanishes outside $U_k\times U_k$, and is given by
\begin{align*}
\tilde{f_k}(z_k,z'_k,t)=\left\{\begin{array}{lr}
-\chi_k(z'_k)[\Delta,\tilde{\chi}_k](z_k)\mathscr{H}_{\mathbb{C}}(z_k,z'_k,t) &\qquad (k>2g-2),\\
-\chi_k(z'_k)[\Delta,\tilde{\chi}_k](z_k)\mathscr{H}_{\curvearrowleft}(r_k,\varphi_k,r'_k,\varphi'_k,t | 4\pi)
 &\qquad (k\le 2g-2)
\end{array}\right.
\end{align*}
in $U_k\times U_k$. Since $\tilde{\chi}_k\chi_k=\chi_k$, the coordinate $z_k$ is separated from $z'_k$ on the support of $\chi_k(z'_k)[\Delta,\tilde{\chi}_k](z_k)$, and, hence, $\tilde{f_k}(z_k,z'_k,t)$ and all its derivatives decay exponentially as $t\to +0$. Therefore, $\tilde{f}(z,z',t)$ and all its derivatives (considered as functions of coordinates $z_k,z_k'$) decay exponentially as $t\to +0$.

Next, since each $\mathscr{H}_{0,k}$ in (\ref{parametrix parts}) is a fundamental solution to the heat equation in $U_k\times U_k$ while the cut-off functions obey $\tilde{\chi}_k\chi_k=\chi_k$, $\sum_k \chi_k=1$, one has 
\begin{equation}
\label{parametrix initial condition}
\lim_{t\to +0}\mathscr{H}_0(z,z',t)=\delta(z-z')
\end{equation}
(here the limit is understood in the sense of distributions).

Finally, in view of formula (\ref{heat kernel asymptotics near vertex}) with $\alpha=4\pi$, the asymptotics of $\mathscr{H}_0(\cdot,z',t)$ near conical $P_k$ is of the form (\ref{asymptotics Friedrichs}). Therefore, 
\begin{equation}
\label{parametrix domain}
\mathscr{H}_0(\cdot,z',t)\in{\rm Dom}\Delta_F \qquad (t>0, \ z'\in\dot{X}).
\end{equation}

\subsubsection{Short-time asymptotics of heat kernel}
Since the operator $\Delta_F$ is positive, the $C_0$-semigroup of contractions $t\mapsto {\rm exp}\big(-t\Delta_{F}\big)$ is well-defined. For each $z'\in \dot{X}$, the (unique) solution $\mathscr{H}_1(\cdot,z',\cdot)\in C^{1}(\mathbb{R};L_2(X;C))\cap C(\mathbb{R};{\rm Dom}\Delta_F)$ (where ${\rm Dom}\Delta_F$ is endowed with the graph norm) to the Cauchy problem 
\begin{equation}
\label{h eq h ker op sol}
\begin{split}
(\partial_t+\Delta_{F})\mathscr{H}_1(\cdot,z',t)&=-\tilde{f}(\cdot,z',t'), \quad t>0,\\
\mathscr{H}_1(\cdot,z',0)&=0
\end{split}
\end{equation}
is given by
$$\mathscr{H}_1(\cdot,z',t)=-\int_{0}^{t}{\rm exp}\big((t'-t)\Delta_{F}\big)\tilde{f}(\cdot,z',t')dt'.$$
Similarly, the solution $\mathscr{H}_{1,l}(\cdot,z',\cdot)\in C(\mathbb{R};{\rm Dom}\Delta_F)\cap C^{1}(\mathbb{R};L_2(X;C))$ to
\begin{equation*}
\begin{split}
(\partial_t+\Delta_{F})\mathscr{H}_{1,l}(\cdot,z',t)&=-\partial_{t}^l\tilde{f}(\cdot,z',t), \quad t>0,\\
\mathscr{H}_{1,l}(\cdot,z',0)&=0
\end{split}
\end{equation*}
is given by
$$\mathscr{H}_{1,l}(z,z',t):=\int_0^{t}{\rm exp}\big((t'-t)\Delta_{F}\big)[-\partial_{t'}^l\tilde{f}(\cdot,z',t')]dt'.$$
The solutions $\mathscr{H}_1$ and $\mathscr{H}_{1,l}$ with $l>0$ are related via the equality
\begin{equation}
\label{parametrix remainder}
\mathscr{H}_{1}(z,z',t)=\int\limits_0^t dt_{l-1}\int\limits_0^{t_{l-1}} dt_{l-2}\dots\int\limits_0^{t_1} dt_{1}\mathscr{H}_{1,l}(z,z',t_1).
\end{equation}
Indeed, since $\mathscr{H}_{1,l}(\cdot,z',\cdot)\in C(\mathbb{R};{\rm Dom}\Delta_F)\cap C^{1}(\mathbb{R};L_2(X;C))$, each integration in the right hand side of (\ref{parametrix remainder}) commutes with the action of $\Delta_F$. So, the right hand side of (\ref{parametrix remainder}) is also a solution to (\ref{h eq h ker op sol}); hence, it coincides with $\mathscr{H}^{(1)}(\cdot,z',t)$.

In view of (\ref{parametrix remainder}), one has $\mathscr{H}_1(\cdot,z',\cdot)\in C^{\infty}(\mathbb{R};{\rm Dom}\Delta_F)\cap C^{\infty}(\mathbb{R};L_2(X;C))$ and 
\begin{equation}
\label{parametrix rem diff}
(-\Delta_F)^l \mathscr{H}_1(\cdot,z',t)= \mathscr{H}_{1,l}(\cdot,z',t)+\sum_{s=1}^{l}(-\Delta_F)^{l-s}\partial_t^{s-1}\tilde{f}(\cdot,z',t).
\end{equation}
Recall that the $L_2(X;C)$-norms of $\Delta_F^{k}\partial_t^{s}\tilde{f}(\cdot,z',t)$ (where $k,s=0,1,\dots$) decay exponentially as $t\to +0$ and uniformly with respect to $z'$. Since ${\rm exp}\big((t'-t)\Delta_{F}\big)$ is a contraction in $L_2(X;C)$ for $t'<t$, the same is true for the $L_2(X;C)$-norm of the right hand side of (\ref{parametrix rem diff}). Thus, the $L_2(X;C)$-norms of $\Delta_F^l\mathscr{H}_{1}(\cdot,z',t)$ decay exponentially and uniformly with respect to $z'$ as $t\to +0$. Hence, due to estimates (\ref{C l smoth incr}), $\mathscr{H}_{1}(z_k,z',t)$ ($k> 2g-2$) decays exponentially as $t\to +0$ uniformly with respect to $z_k,z'$. Moreover, due to Proposition \ref{prop elliptic asymptotics}, the functions $\mathscr{H}_{1}(x_k,z',t)$ ($k\le 2g-2$) decay exponentially as $t\to +0$ uniformly with respect to $x_k,z'$.

Put
\begin{equation}
\label{heat kernel expansion}
\mathscr{H}:=\mathscr{H}_0+\mathscr{H}_1.
\end{equation}
In view of (\ref{parametrix domain}) and (\ref{h eq h ker op sol}), we have $(\partial_t+\Delta_F)\mathscr{H}(\cdot,z',t)=0$ for $t>0$ and $z\in\dot{X}$. Due to (\ref{parametrix initial condition}) and the exponential decrease of $\mathscr{H}_{1}$ as $t\to +0$, we have $\lim_{t\to +0}\mathscr{H}(z,z',t)=\delta(z-z')$ in the sense of distributions.

Let $0<\lambda^F_1\le \lambda^F_2\le\dots$ be the eigenvalues of $\Delta_F$ counted with their multiplicities and $\{u^F_k\}$ be the orthonormal basis of the corresponding eigensections. For $z'\in\dot{X}$ and $t>0$, put $c_k(t,z'):=(\mathscr{H}(\cdot,z',t),u^F_k)_{L_2(X;C)}$.
Then 
\begin{align*}
\lambda^F_k c_k(t,z')=(\mathscr{H}(\cdot,z',t),\Delta_F u^F_k)_{L_2(X;C)}=&(\Delta_F \mathscr{H}(\cdot,z',t),u^F_k)_{L_2(X;C)}=\\=-(\partial_t\mathscr{H}&(\cdot,z',t),u^F_k)_{L_2(X;C)}=-\partial_t c_k(t,z').
\end{align*}
Since $\lim_{t\to +0}\mathscr{H}(z,z',t)=\delta(z-z')$, we have 
$$\lim_{t\to +0}c_k(t,z')=\overline{u^F_k(z')}.$$ 
Therefore, $c_k(t,z')={\rm exp}(-\lambda^F_k t)\overline{u^F_k(z')}$ and 
\begin{equation}
\label{heat kernel}
\mathscr{H}(z,z',t)=\sum_{k=1}^{\infty}e^{-\lambda^F_k t}u^F_k(z)\overline{u^F_k(z')}.
\end{equation}
Thus, $\mathscr{H}$ is the integral kernel of the operator exponent $e^{-t\Delta_F}$ (the heat kernel corresponding to the Friedrich extension of Laplacian). Due to (\ref{heat kernel expansion}) and the exponential decrease of $\mathscr{H}_{1}$ as $t\to +0$, the main term of asymptotics of $\mathscr{H}(z,z',t)$ as $t\to +0$ is given by (\ref{heat kernel}), (\ref{parametrix parts}). Hence, formula (\ref{heat kernel trace asymptotics conic}) yields
\begin{equation}
\label{heat trace asymptotics}
\int_X \mathscr{H}(z,z,t)h(z)dS(z)=\frac{{\rm Area}(X,|\omega|^2)}{\pi t}-\frac{3(g-1)}{8}+O(e^{-\epsilon /t}), \quad t\to +0,
\end{equation}
where $\epsilon>0$.

\subsubsection{Zeta function of $\Delta_F$}
In view of (\ref{heat kernel}), the left hand side of (\ref{heat trace asymptotics}) coincides with $K(t|\Delta_F):={\rm Tr}\,e^{-t\Delta_F}=\sum_k e^{-\lambda^F_k t}$. Hence, 
\begin{equation}
\label{spec partition asymp}
K(t|\Delta_F)=\frac{{\rm Area}(X,|\omega|^2)}{\pi t}-\frac{3(g-1)}{8}+\tilde{K}(t),
\end{equation} 
where $\tilde{K}(t)=O(e^{-\epsilon /t})$ as $t\to +0$. Also, $K(t|\Delta_F)=O(e^{-\lambda^F_1 t})$ as $t\to +\infty$. Since the zeta-function $\zeta(s|A)={\rm Tr}A^{-s}$ of the self-adjoint operator $A$ is related to $K(t|A):={\rm Tr}\,e^{-tA}$ via
\begin{equation}
\label{zeta function via spec part func}
\zeta(s|A)=\frac{1}{\Gamma(s)}\int_{0}^{+\infty}t^{s-1}K(t|A)dt,
\end{equation}
formula (\ref{spec partition asymp}) yields
\begin{align*}
\zeta(s|\Delta_F)=\frac{{\rm Area}(X,|\omega|^2)}{\pi \Gamma(s)}\int_{0}^{1}t^{s-2}dt-\frac{3(g-1)}{8\Gamma(s)}\int_{0}^{1}t^{s-1}dt&+\\
+\frac{1}{\Gamma(s)}\Big(\int_{0}^{1}t^{s-1}\tilde{K}(t|\Delta_F)dt+\int_{1}^{+\infty}t^{s-1}K(t|\Delta_F)dt\Big)&=\frac{{\rm Area}(X,|\omega|^2)}{\pi(s-1)\Gamma(s)}+\tilde{\zeta}(s),
\end{align*}
where $\tilde{\zeta}$ is an entire function. Thus, $s\mapsto\zeta(s|\Delta_F)$ is continued to a meromorphic function with the unique (simple) pole at $s=1$. As a corollary, the determinant ${\rm det}\Delta_F={\rm exp}\big(-\partial_s\zeta(s|\Delta_F)\big|_{s=0}\big)$ of $\Delta_F$ is well-defined.

\section{Comparison formulas}
\label{sec comparing det}
In this section, we prove comparison formula (\ref{comparing determinants F S}) for determinants of $\Delta_F$ and $\Delta_S$.

\subsection{Comparing resolvents}
In this subsection, we introduce the scattering matrix $T(\lambda)$ associated with $\Delta_F$ and $\Delta_S$ and describe its properties. Then, we derive the formula for the difference of resolvent kernels $\mathscr{R}^F_{\lambda}$ and $\mathscr{R}^S_{\lambda}$ of $\Delta_F$ and $\Delta_S$ in terms of $T(\lambda)$.

\subsubsection{Scattering matrix}
Introduce the solutions
\begin{equation}
\label{special solutions}
\mathscr{V}^\lambda_{k,0,-}=\chi_k f_{k,0,-}-(\Delta_S-\lambda)^{-1}\tilde{\mathscr{W}}^\lambda_{k,0,-}
\end{equation}
to the equation $(\Delta-\lambda)\mathscr{V}=0$, where $\lambda$ is not an eigenvalue of $\Delta_S$ and 
$$\tilde{\mathscr{W}}^\lambda_{k,0,-}=(\Delta-\lambda)\big(\chi_k f_{k,0,-}\big)=[\Delta,\chi_k]f_{k,0,-}-\lambda \chi_k f_{k,0,-}.$$
In view of (\ref{Szego domain}) and (\ref{asymptotics Friedrichs}), asymptotic expansion (\ref{elliptic asymptotics}) for $u=(\Delta_S-\lambda)^{-1}\tilde{\mathscr{W}}^\lambda_{k,0,-}$ provided by Proposition \ref{prop elliptic asymptotics} yields
\begin{align}
\label{special solution asymptotics}
\mathscr{V}^\lambda_{k,0,-}-\Big(\chi_k f_{k,0,-}+\sum_{j=1}^{2g-2}T_{jk}(\lambda)\chi_j f_{j,0,+}\Big)\in{\rm Dom}\Delta_S\cap{\rm Dom}\Delta_F,
\end{align}
where $T_{jk}(\lambda)$ are some coefficients. The $(2g-2)\times(2g-2)$-matrix $T(\lambda)$ with entries $T_{jk}(\lambda)$ is called the {\it scattering matrix}. Below, we describe the properties of $T(\lambda)$. Note that neither solutions (\ref{special solutions}) nor the scattering matrix depend on the choices of the cut-off functions $\chi_1,\dots,\chi_{2g-2}$. Since $\mathscr{R}^S_\lambda$ is meromorphic with respect to $\lambda$ and has poles at the eigenvalues of $\Delta_S$, the scattering matrix $\lambda\mapsto T(\lambda)$ is also meromorphic and each pole of $T$ is an eigenvalue of $\Delta_S$.

\

Let us express the scattering matrix $T(\lambda)$ in terms of $\mathscr{R}^S_{\lambda}$.
\begin{lemma}
We have
\begin{align}
\label{special solutions via resolvent kernel}
\mathscr{V}^\lambda_{k,0,-}&=\mathscr{R}^S_\lambda(\cdot,P_k),\\
\label{T-matrix}
T_{jk}(\lambda)&=\mathscr{R}^S_\lambda(P_j,P_k)=\overline{T_{kj}(\overline{\lambda})},\\
\label{derivative of T-matrix}
\partial_\lambda T_{jk}(\lambda)&=(\mathscr{R}^S_\lambda(\cdot,P_k),\mathscr{R}^S_{\overline{\lambda}}(\cdot,P_j))_{L_2(X;C)}.
\end{align}
\end{lemma}
\begin{proof}
Let us prove (\ref{special solutions via resolvent kernel}) and (\ref{T-matrix}). If $z'$ is not a conical point, then the resolvent kernel $\mathscr{R}^\star_{\lambda}$ of $\Delta_\star$ ($\star=F,S$) can be represented as
\begin{equation}
\label{resolvent kernel representation}
\mathscr{R}^S_\lambda(\cdot,z')=\psi_{z'}\mathfrak{R}_\lambda^{0}(\cdot,z')+\big[(\Delta_S-\lambda)^{-1}(\Delta-\lambda)(\psi_{z'}\mathfrak{R}_\lambda^{0})\big](\cdot,z'),
\end{equation}
where $\psi_{z'}$ is a smooth cut-off function on $X$ equal to one near $z'$, the support of $\psi_{z'}$ is sufficiently small, and $\mathfrak{R}_\lambda^{0}$ is a local fundamental solution to $(\Delta-\lambda)v=0$ given by
\begin{equation}
\label{log singularity}
\mathscr{R}_\lambda^{0}(z,z'):=
\left\{\begin{array}{ll}
(2/\pi)K_0(-2i|z-z'|\sqrt{\lambda}) & (\lambda\ne 0)\\
(-2/\pi){\rm log}(|z-z'|) & (\lambda=0),
\end{array}\right.
\end{equation}
in coordinates (\ref{z coordinate}), where $K_0$ is the Macdonald function. 

Green formula (\ref{Green formula}) with $U=U_\epsilon=X_\epsilon(P_1,\dots,P_{2g-2},z)$, $L=C$, and $f=\chi_k f_{k,0,-}$, $f'=\mathscr{R}^S_{\overline{\lambda}}(\cdot,z)$ and the symmetry $\mathscr{R}^S_{\lambda}(z,z')=\overline{\mathscr{R}^S_{\overline{\lambda}}(z',z)}$ lead to 
\begin{align}
\label{referadd 1}
\begin{split}
((\Delta-\lambda)\big(&\chi_k f_{k,0,-}\big),\mathscr{R}^S_{\overline{\lambda}}(\cdot,z))_{L_2(X;C)}-(\chi_k f_{k,0,-},(\Delta-\overline{\lambda})\mathscr{R}^S_{\overline{\lambda}}(\cdot,z))_{L_2(X;C)}=\\
=&\lim_{\epsilon\to 0}\frac{1}{2i}\Big(\int\limits_{{\rm dist}(z',z)=\epsilon}+\int\limits_{{\rm dist}(z',P_k)=\epsilon}\Big)\big[\chi_k f_{k,0,-}\,\partial_{z'}\mathscr{R}^S_{\lambda}(z,z')dz'+\mathscr{R}^S_{\lambda}(z,z')\,\partial_{\overline{z'}}(\chi_k f_{k,0,-})d\overline{z'}\big].
\end{split}
\end{align}
Note that the first term in the right hand side of (\ref{resolvent kernel representation}) vanishes near conical points, while the second term is smooth in $\dot{X}$ and belongs to ${\rm Dom}\Delta_\star$. Hence, Proposition \ref{prop elliptic asymptotics} provides asymptotic expansion (\ref{elliptic asymptotics}) of $u=\mathscr{R}_\lambda^{S}(\cdot,z')$ near each $P_k$, where $c_{k,0,+}=\mathscr{R}_\lambda^{S}(P_k,z')$. At the same time, $\partial_{z'}\mathscr{R}^S_{\lambda}(z,z')=-1/\pi(z'-z)+O(1)$ as $z'\to z$, so the right-hand side of (\ref{referadd 1}) is equal to 
$$[\chi_k f_{k,0,-}](z)-\overline{\mathscr{R}^S_{\overline{\lambda}}(P_k,z)}$$
Due to the equalities $(\Delta-\overline{\lambda})\mathscr{R}^S_{\overline{\lambda}}(\cdot,z)=0$ and $(\Delta-\lambda)\big(\chi_k f_{k,0,-}\big)=\tilde{\mathscr{W}}^\lambda_{k,0,-}$, the left-hand side of (\ref{referadd 1}) is equal to 
$$(\tilde{\mathscr{W}}^\lambda_{k,0,-},\mathscr{R}^S_{\overline{\lambda}}(\cdot,z))_{L_2(U_\epsilon;C)}=\int\limits_X \tilde{\mathscr{W}}^\lambda_{k,0,-}(z')\mathscr{R}^S_{\lambda}(z,z')h(z')dS(z')=(\Delta_S-\lambda)^{-1}\tilde{\mathscr{W}}^\lambda_{k,0,-}(z).$$
Thus, (\ref{referadd 1}) takes the form
$$(\Delta_S-\lambda)^{-1}\tilde{\mathscr{W}}^\lambda_{k,0,-}(z)=[\chi_k f_{k,0,-}](z)-\mathscr{R}^S_{\lambda}(z,P_k).$$
Comparing the last equality with (\ref{special solutions}) leads to (\ref{special solutions via resolvent kernel}). Now, substituting (\ref{special solutions via resolvent kernel}) with $z=x_k$ into (\ref{special solution asymptotics}) and passing to the limit as $x_k\to 0$ yields (\ref{T-matrix}).

Now, let us prove (\ref{derivative of T-matrix}). In view of Proposition \ref{prop elliptic asymptotics}, one can differentiate the equation $(\Delta-\lambda)\mathscr{V}^\lambda_{k,0,-}=0$ and asymptotics (\ref{special solution asymptotics}) with respect to $\lambda$. As a result, one obtains $(\Delta-\lambda)\partial_\lambda\mathscr{V}^\lambda_{k,0,-}=\mathscr{V}^\lambda_{k,0,-}$ and 
$$\partial_\lambda\mathscr{V}^\lambda_{k,0,-}-\sum_{j=1}^{2g-2}\partial_\lambda T_{jk}(\lambda)\chi_j f_{j,0,+}\in{\rm Dom}\Delta_S\cap{\rm Dom}\Delta_F.$$
In view of Green formula (\ref{Green formula with asymptotics}) and the equality $(\Delta-\overline{\lambda})\mathscr{V}^{\overline{\lambda}}_{j,0,-}=0$, these relations imply 
\begin{align*}
(\mathscr{V}^\lambda_{k,0,-},\mathscr{V}^{\overline{\lambda}}_{j,0,-}&)_{L_2(X;C)}=\\
=((\Delta-\lambda)&\partial_\lambda\mathscr{V}^\lambda_{k,0,-},\mathscr{V}^{\overline{\lambda}}_{j,0,-})_{L_2(X;C)}-(\mathscr{V}^\lambda_{k,0,-},(\Delta-\overline{\lambda})\mathscr{V}^{\overline{\lambda}}_{j,0,-})_{L_2(X;C)}=\partial_\lambda T_{jk}(\lambda).
\end{align*}
In view of (\ref{special solutions via resolvent kernel}), the last formula is equivalent to (\ref{derivative of T-matrix}).
\end{proof}

Now, let us derive the asymptotics of scatterring matrix at infinity.
\begin{lemma}
The following asymptotics holds
\begin{align}
\label{asymptotics of detT}
{\rm det}T(\lambda)=t_\infty(-\lambda)^{p_\infty}+O(e^{-\epsilon|\lambda|}), \qquad \Re\lambda\to -\infty,
\end{align}
where
\begin{equation}
\label{asymptotics of det T constants}
t_\infty:=\Gamma(3/4)^{4-4g}, \qquad p_\infty:=\frac{1-g}{2}.
\end{equation}
\end{lemma}
\begin{proof}
Let us represent $\mathscr{V}^\lambda_{k,0,-}$ in the form
\begin{equation}
\label{special solution large lambda}
\mathscr{V}^\lambda_{k,0,-}=\chi_k f^\lambda_{k,0,-}-(\Delta_S-\lambda)^{-1}[\Delta,\chi_k]f^\lambda_{k,0,-},
\end{equation}
where $f^\lambda_{k,0,-}(x_k)=\mathfrak{D}_{\lambda}(|x_k|)$ is a local solution $f^\lambda_{k,0,-}$ to $(\Delta-\lambda)u=0$ near $P_k$. Since
$$\lambda\mathfrak{D}_{\lambda}(|x_k|)=\Delta f^\lambda_{k,0,-}=-|x_k|^{-1}\frac{\partial}{\partial x_k}|x_k|^{-1}\frac{\partial}{\partial \overline{x_k}}\mathfrak{D}_{\lambda}(|x_k|)=-\frac{1}{4|x_k|^2}\mathfrak{D}_{\lambda}''(|x_k|),$$
we have 
$$\mathfrak{D}_{\lambda}''(\rho)+4\lambda\rho^2\mathfrak{D}_{\lambda}(\rho)=0, \qquad \rho=|x_k|.$$ 
The solution exponentially decaying as $|\lambda|r^{4}\to +\infty$, $|{\rm arg}\lambda|>\epsilon>0$ is given by 
$$\mathfrak{D}_{\lambda}(r)=c(\lambda)D_{-1/2}(2r(-\lambda)^{1/4}),$$
where $D_{-1/2}$ is the parabolic cylinder function. In view of the Taylor expansion  
$$D_{-\frac{1}{2}}(z)=\frac{\pi^{3/2}2^{-1/4}}{\Gamma(3/4)}-\frac{\pi^{3/2}2^{1/4}}{\Gamma(1/4)}z+O(z^4), \qquad z\to 0,$$
condition (\ref{special solution asymptotics}) implies $c(\lambda)=(-2/\lambda)^{1/4}\pi^{-3/2}\Gamma(3/4)^{-1}$. Therefore,
\begin{equation}
\label{asymptotics of the special solution near zero}
f^\lambda_{k,0,-}(x_k)=f_{k,0,-}(x_k)+\Gamma(3/4)^{-2}(-\lambda)^{-1/4}+O(|x_k|^4).
\end{equation}
Note that $[\Delta,\chi_k]f^\lambda_{k,0,-}$ and all its derivatives with respect to $x_k$ decay exponentially as $\Re\lambda\to -\infty$. In view of Proposition \ref{elliptic asymptotics}, the coefficients in asymptotics (\ref{Szego domain}) of $u=(\Delta_S-\lambda)^{-1}[\Delta,\chi_k]f^\lambda_{k,0,-}$ decay exponentially as $\Re\lambda\to -\infty$. Thus, formulas (\ref{special solution asymptotics}), (\ref{special solution large lambda}), and (\ref{asymptotics of the special solution near zero}) imply
$$T_{jk}(\lambda)=\delta_{jk}\Gamma(3/4)^{-2}(-\lambda)^{-1/4}+O(e^{-\epsilon|\lambda|}),$$
where $\epsilon>0$. As a corollary, we obtain (\ref{asymptotics of detT}), (\ref{asymptotics of det T constants}).
\end{proof}

\

Finally, we derive an expression for $T(0)$. Since (\ref{Szego to Green}) is valid for the Green function $\mathcal{G}^{\mathcal{S}}$ of the Szeg\"o extension $\Delta_S$, formula (\ref{T-matrix}) implies 
\begin{equation}
\label{S and T matrices 1}
T_{jk}(0)=\mathcal{G}^{\mathcal{S}}(P_j,P_k)=\pi^{-2}\mathfrak{S}_{jk},
\end{equation}
where $\mathfrak{S}$ is given by (\ref{Bergman matrix}). Hence, in view of (\ref{Szego kernel expression}), we obtain (\ref{expression for T(0)}).

\subsubsection{Comparing resolvents of $\Delta_F$ and $\Delta_S$}
\begin{lemma}
If $\lambda$ does not belong to the spectra of $\Delta_S$, $\Delta_F$, then the operator 
$$D(\lambda):=(\Delta_S-\lambda)^{-1}-(\Delta_F-\lambda)^{-1}$$ 
is well defined, it is of finite rank, and its trace is given by
\begin{equation}
\label{comparing traces of resolvents F and S}
\begin{split}
{\rm Tr}D(\lambda)=\partial_\lambda{\rm log}{\rm det}T(\lambda).
\end{split}
\end{equation}
\end{lemma} 
\begin{proof}
Let $y^S:=(\Delta_S-\lambda)^{-1}f$, where $f\in C_c^\infty(\dot{X};C)$ and $\lambda$ is not an eigenvalue of $\Delta_S$. In view of Proposition \ref{prop elliptic asymptotics} and formula (\ref{Szego domain}), the expansion 
\begin{equation}
\label{asymptotics of solution modulo dom S cap F}
y^S-\sum_j c_{j,0,+}(f)\chi_j f_{j,0,+}\in{\rm Dom}\Delta_S\cap{\rm Dom}\Delta_F
\end{equation}
is valid, where the coefficients can be found from Green formula (\ref{Green formula with asymptotics}) with $f=\mathscr{V}^{\overline{\lambda}}_{j,0,-}$, $f'=y^{S}$;
\begin{equation}
\label{resolvent diff projection coeff}
c_{j,0,+}(f)=(f,\mathscr{V}^{\overline{\lambda}}_{j,0,-})_{L_2(X;C)}=[\text{(\ref{special solutions via resolvent kernel})}]=\int_{X}f(z')\mathscr{R}^S_\lambda(P_j,z')h(z')dS(z').
\end{equation}
We search for a solution $y^F$ to $(\Delta_F-\lambda)y^F=f$ in the form
\begin{equation}
\label{resolvent diff solution form}
y^{F}=y^{S}-\sum_k c_k \mathscr{V}^\lambda_{k,0,-}.
\end{equation}
In view of (\ref{asymptotics of solution modulo dom S cap F}) and (\ref{special solution asymptotics}), the right hand side belongs to ${\rm Dom}\Delta_F$ if and only if
\begin{equation}
\label{resolvent diff scat matr}
c_{j,0,+}(f)=\sum_k T_{jk}(\lambda)c_k.
\end{equation}
Combining formulas (\ref{resolvent diff solution form}), (\ref{resolvent diff scat matr}), (\ref{resolvent diff projection coeff}), and (\ref{special solutions via resolvent kernel}), we obtain 
$$y^S(z)-y^F(z)=\int_{X}f(z')\Big[\sum_{k,j}\mathscr{R}^S_\lambda(z,P_k)T^{-1}_{kj}(\lambda)\mathscr{R}^S_\lambda(P_j,z')\Big]h(z')dS(z').$$
Since the $f\in C_c^\infty(\dot{X};C)$ is arbitrary, we have
\begin{equation*}
\mathscr{R}_\lambda^S(z,z')-\mathscr{R}_\lambda^F(z,z')=\sum_{k,j}\mathscr{R}^S_\lambda(z,P_k)T^{-1}_{kj}(\lambda)\mathscr{R}^S_\lambda(P_j,z').
\end{equation*}
Using the bra–ket notation, we rewrite the last formula as
\begin{equation}
\label{comparing resolvents F and S}
D(\lambda)=(\Delta_S-\lambda)^{-1}-(\Delta_F-\lambda)^{-1}=\sum_{k,j}T^{-1}_{kj}(\lambda)|\mathscr{R}^S_\lambda(\cdot,P_k)\rangle \langle\mathscr{R}^S_{\overline{\lambda}}(\cdot,P_j)|,
\end{equation}
where $|f\rangle$ denotes an element of $L_2(X;C)$ while $\langle f|$ is the linear bounded functional on $L_2(X;C)$ given by $\langle f|f'\rangle:=(f',f)_{L_2(X;C)}$. In particular, $D(\lambda)$ is of finite rank. Formula (\ref{comparing resolvents F and S}) shows that the zeroes of ${\rm det}T(\lambda)$ are the eigenvalues of $\Delta_F$. In view of (\ref{derivative of T-matrix}) and the equality
$${\rm Tr}|a\rangle \langle b|=\langle b|a\rangle,$$
formula (\ref{comparing resolvents F and S}) implies
\begin{equation*}
\begin{split}
{\rm Tr}D(\lambda)=\sum_{k,j}T^{-1}_{kj}(\lambda)(\mathscr{R}^S_\lambda(\cdot,P_k),\mathscr{R}^S_{\overline{\lambda}}(\cdot,P_j))_{L_2(X;C)}=\sum_{k,j}T^{-1}_{kj}(\lambda)\partial_\lambda T_{jk}(\lambda)=\\
={\rm Tr}\big(T^{-1}(\lambda)\partial_\lambda T(\lambda)\big)=\partial_\lambda{\rm Tr}{\rm log}T(\lambda)=\partial_\lambda{\rm log}{\rm det}T(\lambda).
\end{split}
\end{equation*}
Thus, we have proved (\ref{comparing traces of resolvents F and S}). 
\end{proof}

\subsection{Comparing zeta functions}
In this subsection, we derive formula (\ref{comparing determinants F S}). To this end, we express the difference of zeta functions of $\Delta_F$ and $\Delta_S$ in terms of the trace of the difference of their resolvents, and then apply formula (\ref{comparing traces of resolvents F and S}).

\subsubsection{Some operator relations}

\

\paragraph{\it Formula for powers of operator.}
If $A$ is a positive definite operator with discrete spectrum, then
\begin{equation}
\label{powers of operator via resolvent}
A^{-s}=\frac{1}{2\pi i}\int_\gamma \xi^{-s}(A-\xi)^{-1}d\xi \qquad (\Re s>0),
\end{equation}
where $\gamma$ is the union of the paths 
\begin{equation}
\label{paths}
\gamma_\pm(t)=|t|e^{\pm i(\pi+0)} \quad (\pm t\in (\epsilon,+\infty)), \qquad \gamma_\circ(t)=\epsilon e^{i\varphi} \quad (\varphi\in (-i\pi,+i\pi)),
\end{equation}
and the number $\epsilon>0$ is sufficiently small. Indeed, since $(\zeta-\lambda)^{-1}=O(|\lambda|^{-1})$ as $\Re\lambda\to-\infty$ uniformly with respect to $\zeta\in(0,\infty)\supset {\rm Sp}(A)$, one has $\|(A-\lambda)^{-1}\|=O(|\lambda|^{-1})$ as $\Re\lambda\to-\infty$. Then the norm of the integrand is $O(|\xi|^{-1-\Re s})$ as $\gamma\ni\xi\to\infty$ and the integral in the right hand side converges (in the operator norm). In particular, the right hand side is a bounded operator, while the left hand side is a bounded operator since $|\xi^{-s}|$ is bounded on the spectrum of $A$. So, it remains to check (\ref{powers of operator via resolvent}) on elements of some orthonormal basis $\{u_k\}_k$. Let us choose this basis in such a way that $u_k$ is an eigenfunction of $A$ corresponding to the eigenvalue $\lambda_k$. Applying the both sides of (\ref{powers of operator via resolvent}) to $u_k$, and taking into account that $\Re s>0$, one obtains
$$\Big(\frac{1}{2\pi i}\int_\gamma \xi^{-s}(\lambda_k-\xi)^{-1}d\xi\Big)u_k=-{\rm Res}_{\xi=\lambda_k}(\xi^{-s}(\lambda_k-\xi)^{-1})=\lambda_k^s u_k=A^{-s}u_k.$$
Hence, formula (\ref{powers of operator via resolvent}) is proved.

\

\paragraph{\it Connection between zeta functions of self-adjoint extensions of symmetric operator.} Let $A'$ and $A$ be positive definite operators ($A,A'\ge\epsilon I>0$) with discrete spectra. Also, we suppose that $A'$ and $A$ are self-adjoint extensions of the same symmetric operator $A_0$. Then, for any $\lambda\not\in [\epsilon>0,+\infty)$, the difference of their resolvents 
$$D(\lambda):=(A'-\lambda)^{-1}-(A-\lambda)^{-1}$$
is an operator of finite rank not exceeding the deficiency index of $A_0$. According to (\ref{powers of operator via resolvent}), we have
$$A'^{-s}-A^{-s}=\frac{1}{2\pi i}\int_\gamma \xi^{-s} D(\xi)d\xi \qquad (\Re s\ge 2),$$
%see the explanation in ``net dobavochki.pdf''
where the integral in the right hand side converges in the operator norm as well as in the trace-norm $\|\cdot\|_1$. Therefore, one can apply the trace to both sides and interchange the trace and the integration in the right hand side. As a result, one obtains 
\begin{equation}
\label{difference of zeta functions}
\zeta(s|A')-\zeta(s|A)={\rm Tr}(A'^{-s}-A^{-s})=\frac{1}{2\pi i}\int_\gamma \xi^{-s} {\rm Tr}D(\xi)d\xi,
\end{equation}
where the right hand side is well-defined at least for $\Re s\ge 2$. 

\subsubsection{Comparing $\zeta(\cdot|\Delta_S)$ and $\zeta(\cdot|\Delta_F)$.}
\label{ssec comp zetas}
Let us substitute $A'=\Delta_S$, $A=\Delta_F$ into (\ref{difference of zeta functions}) and apply (\ref{comparing traces of resolvents F and S}). As a result, we obtain
\begin{align*}
\zeta(s|\Delta_S)-\zeta(s|\Delta_F)=\frac{1}{2\pi i}\int\limits_\gamma \xi^{-s} \frac{t'(\xi)}{t(\xi)}d\xi, \qquad (\Re s\ge 2),
\end{align*}
where $t:={\rm det}T$. Denote the right hand side by $J(s)$; then $J=J_1+J_2$, where
\begin{align*}
J_1(s):&=\frac{1}{2\pi i}\int_{\gamma_\circ} \xi^{-s}d{\rm log}t(\xi), \\
J_2(s):&=\frac{1}{2\pi i}\int_{\gamma_++\gamma_-} \xi^{-s}\frac{t'(\xi)}{t(\xi)}d\xi=\frac{{\rm sin}(\pi s)}{\pi}\int_{-\infty}^{-\epsilon} |\xi|^{-s}\frac{t'(\xi)}{t(\xi)}d\xi
\end{align*}
and $\gamma_\circ$, $\gamma_\pm$ are given by (\ref{paths}). Note that $s\mapsto J_1(s)$ is holomorphic for any $s\in\mathbb{C}$. Recall that the function $\lambda\mapsto t(\lambda)={\rm det}T(\lambda)$ is holomorphic and non-zero outside spectra of $\Delta_S$ and $\Delta_F$. Due to this, one can choose a branch of ${\rm log}t$ which is holomorphic outside $[2\epsilon,+\infty)$. Integrating by parts, we obtain
\begin{align*}
J_1(s)=-\frac{{\rm sin}(\pi s)}{\pi}\epsilon^{-s}{\rm log}t(-\epsilon)+sJ_3(s),
\end{align*}
where
$$J_3(s)=\frac{1}{2\pi i}\int_{\gamma_\circ}\xi^{-(s+1)}{\rm log}t(\xi)d\xi,$$
In view of the Cauchy formula, we have $J_3(0)={\rm log}t(0)$. Formula (\ref{asymptotics of detT}) can be rewritten as ${\rm log}t=p_\infty{\rm log}(-\lambda)+{\rm log}t_\infty+\tilde{q}$, where $\tilde{q}(\lambda)$ and its derivatives decay exponentially as $\Re\lambda\to -\infty$. Hence, $J_2(s)=J_4(s)+J_5(s)$, where
\begin{align*}
J_4(s):=\frac{p_\infty{\rm sin}(\pi s)}{\pi}\int_{\epsilon}^{+\infty}\xi^{-s-1}d\xi, \qquad J_5(s):=\frac{{\rm sin}(\pi s)}{\pi}\int_{-\infty}^{-\epsilon} |\xi|^{-s}\partial_\xi\tilde{q}(\xi)d\xi.
\end{align*}
Note that $J_4(s)=-p_\infty (\pi s)^{-1}{\rm sin}(\pi s)\epsilon^{-s}$ for $\Re s>0$; the same equality is valid for analitic continuation of $J_4$ on the whole $\mathbb{C}$. Since $\partial_\xi\tilde{q}(\xi)$ decays exponentially as $\Re\lambda\to -\infty$, the integral $J_5$ is analytic on $\mathbb{C}$. In addition, integration by parts yields
$$J_5(s)=\frac{{\rm sin}(\pi s)}{\pi}\Big(\epsilon^{-s}\tilde{q}(-\epsilon)-sJ_6(s)\Big)$$
where
$$J_6(s)=\int_{-\infty}^{-\epsilon}\tilde{q}(\xi)|\xi|^{-s-1}d\xi.$$
Summing up the above equalities, we obtain
\begin{align*}
J(s)=sJ_3(s)-p_\infty \frac{{\rm sin}(\pi s)}{\pi s\epsilon^{s}}-\frac{{\rm sin}(\pi s)}{\pi}\Big(\epsilon^{-s}({\rm log}t(-\epsilon)-\tilde{q}(-\epsilon))+sJ_6(s)\Big)=\\
=sJ_3(s)-(\pi s)^{-1}{\rm sin}(\pi s)\Big(\epsilon^{-s}\Big[p_\infty+s(p_\infty{\rm log}(\epsilon)+{\rm log}t_\infty)\Big]-s^2J_6(s)\Big).
\end{align*}
where the right hand side is an entire function of $s$. In particular,
$$J(0)=p_\infty, \qquad J'(0)=J_3(0)-{\rm log}t_\infty={\rm log}(t(0)/t_\infty).$$ 
Hence, 
\begin{align*}
\frac{{\rm det}\Delta_F}{{\rm det}\Delta_S}&=\frac{{\rm exp}\big(-\partial_s\zeta(s|\Delta_F)\big|_{s=0}\big)}{{\rm exp}\big(-\partial_s\zeta(s|\Delta_S)\big|_{s=0}\big)}=e^{J'(0)}=\\
=&\frac{{\rm det}T(0)}{t_\infty}=[\text{(\ref{asymptotics of det T constants})}]=\Gamma(3/4)^{4(g-1)}{\rm det}T(0).
\end{align*}
Theorem \ref{prop comparing determinants F S} is proved.

\section{Dependence of determinants on moduli}
In this section, we prove formula (\ref{determinant expression}) which describes the dependence of ${\rm det}\Delta_S$ on moduli.

\subsection{Surface families corresponding to variation of moduli}
\label{sec surface families}
The moduli space $H_g(1, \dots, 1)$ is a complex orbifold of dimension $4g-3$ (see \cite{KZ}). Let $(X,\omega)\in H_g(1, \dots, 1)$. Choose paths $a_i,b_i$ ($i=1,\dots,g$) representing a canonical basis in the homology group $H_1(X,\mathbb{Z})$. Let $l_k$ ($k=2,\dots,2g-3$) be a path which connects $P_1$ to $P_k$ and does not intersect any of $a_i,b_i$. Then $\{a_i,b_i,l_k\}$ provides a basis in a relative homology group $H_1(X,\{P_1,\dots,P_{2g-2}\},\mathbb{Z})$ (a ``marking''). One can continuously transport this basis to nearby $(X',\omega')$ in $H_g(1, \dots, 1)$ and thereby endow them with markings $a'_i,b'_i,l'_k$ (see \cite{KonZor}, some explicit way to do this is described below).

Given a marking $a_i,b_i$ ($i=1,\dots,g$), the spin structure $C$ on $X$ is determined by specifying its characteristic $(p,q)$ (where $p,q\in\{0,\frac{1}{2}\}^{2g}$) in the equation 
\begin{equation}
\label{specifying spinor C}
\mathcal{A}_{z_0}(C)=-K_{z_0}+\mathbb{B}p+q,
\end{equation}
where $\mathcal{A}_{z_0}$, $K_{z_0}$, $\mathbb{B}$ are the Abel transform, the vector of Riemann constants, and the $b$-period matrix associated with $X$, respectively. This spin structure is extended to nearby $(X',\omega')$ in such a way that 
\begin{equation}
\label{specifying spinor C prime}
\mathcal{A}'_{z'_0}(C')=-K'_{z'_0}+\mathbb{B}'p+q,
\end{equation}
where $C'$ is a spinor bundle over $X$ while $\mathcal{A}'_{z'_0}$, $K'_{z'_0}$, $\mathbb{B}'$ are the Abel transform, the vector of Riemann constants, and the period matrix associated with $X'$, respectively. It is worth noting that the marking $a_i,b_i,l_k$ and the spin structure $C$ may be changed after continuation along the closed loop in $H_g(1, \dots, 1)$ (see \cite{Atyah}).

With each triple $(X',\omega',C')$, we associate the Szeg\"o self-adjoint extension $\Delta'_S$ of the Dolbeault Laplacian on $C'$. From now on, we deal only with Szeg\"o extensions and omit the symbol $S$ in the notation of operators, resolvent kernels, Green functions e.t.c.. In addition, primed symbols denote the objects related to the pair $(X',\omega')$, while unprimed ones are associated with $(X,\omega)$.

\begin{figure}[h]
\center{\includegraphics[width=0.4\linewidth]{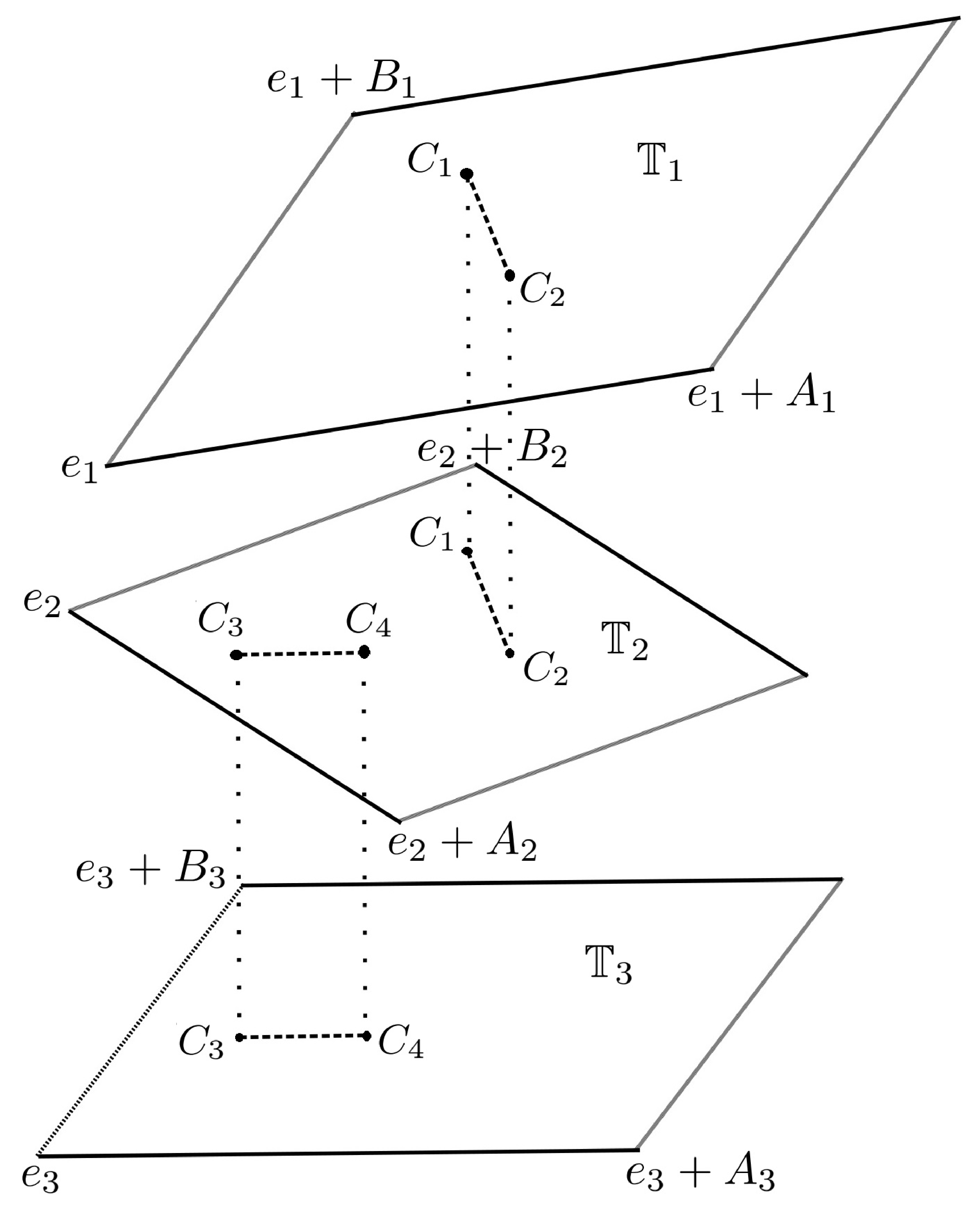}}
\caption{A genus 3 example of a surface representing coordinates (\ref{moduli space coordinates}). Here the opposite side of parallelograms are identified.}
\label{surf}
\end{figure}

The local coordinates (from now on, {\it moduli}) of $(X,\omega)$ are given by
\begin{equation}
\label{moduli space coordinates}
\begin{split}
A_i=\int\limits_{a_j}\omega, \ B_i=\int\limits_{b_j}\omega, \ C_k=\int\limits_{l_k}\omega 
\end{split}
\end{equation}
(see \cite{KonZor}, p.5). As shown on pp.5,6, \cite{KokKorot}, for generic points of $H_g(1, \dots, 1)$ (and for all the points if $g=2$, see \cite{McM}, Theorem 1.7), coordinates (\ref{moduli space coordinates}) can be visualized as follows (see Fig. \ref{surf}). Consider complex tori $\mathbb{T}_j=\mathbb{C}/\Lambda_j$, where $\Lambda_j=\mathbb{Z}A_j+\mathbb{Z}B_j$. Endow these tori with the system of cuts
\begin{eqnarray}
\label{cuts of tori}
{\small
\begin{array}{llllr}
\ [C_1,C_2]\,({\rm mod}\Lambda_1) &&& \text{ on } \mathbb{T}_1\\
\ [C_1,C_2]\,({\rm mod}\Lambda_2), &[C_3,C_4]\,({\rm mod}\Lambda_2) && \text{ on } \mathbb{T}_2\\
\ &[C_3,C_4]\,({\rm mod}\Lambda_3), &[C_5,C_6]\,({\rm mod}\Lambda_3) & \text{ on } \mathbb{T}_3\\
& \dots && \ \ \ \dots\\
\ &[C_{2g-5},C_{2g-4}]\,({\rm mod}\Lambda_{g-1}) & [C_{2g-3},C_{2g-2}]\,({\rm mod}\Lambda_{g-1})& \text{ on } \mathbb{T}_{g-1}\\
\ &&[C_{2g-3},C_{2g-2}]\,({\rm mod}\Lambda_g)& \text{ on } \mathbb{T}_g
\end{array}}
\end{eqnarray}
Making cross-gluing along the edges of the cuts, we obtain the compact Riemann surface $X$ of genus $g$ while the differential $dz$ on $\mathbb{C}$ of all $\mathbb{T}_j$ gives rise to the Abelian differential $\omega$ on $X$ (see Fig. \ref{surf}). The holomorphic local coordinates near the ends $C_k\equiv P_k$ of the cuts are given by (\ref{dist coord}). Thus, each $P_k$ is a zero of $\omega$. The basis cycles in $H_1(X,\mathbb{Z})$ are given by 
\begin{equation}
\label{cycles explicit}
a_j=[e_j,A_j+e_j]({\rm mod}\,\Lambda_j), \quad b_j=[e_j,B_j+e_j]({\rm mod}\,\Lambda_j) \quad \text{ on } \mathbb{T}_j,
\end{equation}
where $e_j$ is an arbitrary constant such that $a_j,b_j$ do not intersect the cuts on $\mathbb{T}_j$. In the construction of the perturbed surface $X'$, we will also use the ``shifted'' cycles
\begin{equation}
\label{cycles shifted}
\check{a}_j=[c_j,A_j+c_j]({\rm mod}\,\Lambda_j), \quad \check{b}_j=[c_j,B_j+c_j]({\rm mod}\,\Lambda_j) \quad \text{ on } \mathbb{T}_j,
\end{equation}
where $c_j$ is an arbitrary constant such that $c_j\ne e_j \ ({\rm mod}\,\Lambda_j)$ and $\check{a}_j,\check{b}_j$ do not intersect the cuts on $\mathbb{T}_j$.

Let $\epsilon_0>0$ be sufficiently small and denote by $\mathbb{D}_k$ the $\epsilon_0$-neighborhood of $P_k$ (in the metric $|\omega|^2$). Let $\bar{\mathbb{T}}_j$ be the domain obtained from the torus $\mathbb{T}_j$ by making the cuts (\ref{cuts of tori}) and (\ref{cycles explicit}) and removing all $\mathbb{D}_k$. Continuous sections $f$ of the bundle $C$ on $X$ can be identified with collections $\{f_i;\tilde{f}_k\}$ of continuous (up to the boundaries) functions $f_i$ on $\bar{\mathbb{T}}_i$ and $\tilde{f}_k$ on $\mathbb{D}_k$ obeying the following conditions (cf. the case of a single torus explained in \cite{Witten}, pp.275, 276):
\begin{enumerate}[a)]
\item the boundary values of $f_j$ on the different sides of the cuts $a_j$, $b_j$ are related by multiplication by the automorphy factors $\sigma(a_j|C)=\pm 1$, $\sigma(b_j|C)=\pm 1$;
\item $f_j=\pm f_{j+1}$ on $L^j_\pm$, where $L^j_\pm$ are sides of the cut $[C_{2j-1},C_{2j}]({\rm mod}\,\Lambda_j)$ in $\mathbb{T}_j$ (identified with the corresponding sides of the cut $[C_{2j-1},C_{2j}]({\rm mod}\,\Lambda_{j+1})$ in $\mathbb{T}_{j+1}$). 
\item If $\partial\mathbb{D}_k\cap\partial\bar{\mathbb{T}}_i$ is nonempty, then $\tilde{f_k}=f_j \sqrt{x_k}$ on it, where $x_k$ is given by (\ref{dist coord}). 
\end{enumerate}
\begin{lemma}
Any spinor bundle $C$ on $X$ can be constructed in the way described above. Two spinor bundles $C_1$ and $C_2$ are isomorphic {\rm(}as holomorphic bundles{\rm)} if and only if $\sigma(a_j|C_1)=\sigma(a_j|C_2)$ and $\sigma(b_j|C_1)=\sigma(b_j|C_2)$ for all $j$. 
\end{lemma}
\begin{proof}
Let $\{f_i;\tilde{f}_k\}$ obey a)-c), then $\{f_i;\tilde{f}_k\}$ defines a section $f$ of some holomorphic line bundle $C$. In view of a),b), there is a continuous function $F$ on $\dot{X}$ such that $f^2_j=F|_{\bar{\mathbb{T}}_j}$ for any $j$. If $\partial\mathbb{D}_k\cap\partial\bar{\mathbb{T}}_i$ is nonempty, then $\tilde{f_k}^2=f^2_j\, x_k=(dz_k/dx_k)F$. Then the equations
$$f^2=f^2_j dz \text{ on } \bar{\mathbb{T}}_j, \qquad f^2=\tilde{f_k}^2 dx_k \text{ on } \mathbb{D}_k$$
define a continuous section $f^2$ of the canonical bundle $K$ on $X$. Therefore, $C^2=K$, i.e., $C$ is a spinor bundle.

Next, let $C_1$ and $C_2$ be two spinor bundles and $f^{(1)}\equiv\{f^{(1)}_i,\tilde{f}_k\}$ and $f^{(2)}\equiv\{f^{(2)}_i;\tilde{f}^{(2)}_k\}$ be their sections. Then 
$$h\equiv\{h_i=f^{(2)}_i/f^{(1)}_i; \ \tilde{h}_k=\tilde{f}^{(2)}_k/\tilde{f}^{(1)}_k\}\equiv f^{(2)}/f^{(1)}$$
is a section of the bundle $B:=C_2C_1^{-1}$. The bundles $C_1$ and $C_2$ are isomorphic if and only if $B$ is trivial, i.e., if and only if there is a holomorphic section $h$ which has no zeroes on $X$.

In view of a), the boundary values of $h_j$ on the different sides of the cuts $a_j$, $b_j$ are related via multiplication by the automorphy factors 
$$\sigma(a_j|B)=\sigma(a_j|C_2)/\sigma(a_j|C_1), \qquad \sigma(b_j|B)=\sigma(b_j|C_2)/\sigma(b_j|C_1).$$ 
On the both sides of the the cut $[C_{2j-1},C_{2j}]({\rm mod}\,\Lambda_j)$ in $\mathbb{T}_j$ (identified with the corresponding sides of the cut $[C_{2j-1},C_{2j}]({\rm mod}\,\Lambda_{j+1})$), we have 
$$h_j=\mathfrak{s}_j h_{j+1} \qquad (\mathfrak{s}_j=\pm 1).$$ 
Here $\mathfrak{s}_j=+1$ if the same choice of the ``plus'' side $L_+^j$ of the cut $[C_{2j-1},C_{2j}]({\rm mod}\,\Lambda_j)$ was used in the construction of both $C_1$, $C_2$; otherwise, we have $\mathfrak{s}_j=-1$. Similarly, if $\partial\mathbb{D}_k\cap\partial\bar{\mathbb{T}}_i$ is nonempty, then 
$$\tilde{h}_k=\mathfrak{s}_{kj} h_j \qquad (\mathfrak{s}_{kj}=\pm 1)$$ 
on it. Here $\mathfrak{s}_{kj}=+1$ if the same choice of the branch of the square root $\sqrt{x_k}$ was used in the construction of both $C_1$, $C_2$; otherwise, we have $\mathfrak{s}_{kj}=-1$. Multiplying all $h_i$ with $i\le j$ by $-1$ (or multiplying arbitrary $\tilde{h}_k$ by $-1$), one obtains a section of a bundle isomorphic to $B$. Thus, one can assume that all $\mathfrak{s}_{j}$ and $\mathfrak{s}_{kj}$ are equal to $1$. 

A section $h$ of $B$ can be identified with the continuous (up to the boundary) function $H$ on $X\backslash\bigcup_{i=1}^g(a_i\cup b_i)$ such that $h_j=H|_{\bar{\mathbb{T}}_j}$ and $\tilde{h}_k=H|_{\mathbb{D}_k}$ for any $j,k$. Then the boundary values of $H$ on the different sides of the cuts $a_j$, $b_j$ are related by multiplication by the automorphy factors $\sigma(a_j|B)=\pm 1$, $\sigma(b_j|B)=\pm 1$. The bundle $B$ is trivial if and only if there is such $H$ which is holomorphic and has no zeroes on $X$. In the last case, $H^2$ can be continued to a holomorphic function on the whole $X$. Thus, $H^2$ and $H$ are a non-zero constants on $X$. The latter is possible only if all $\sigma(a_j|B)$ and $\sigma(b_j|B)$ are equal to $1$. 

Therefore, the bundle $C$ is determined (up to isomorphism) by $2g$ automorphy factors $\sigma(a_j|C)=\pm 1$, $\sigma(b_j|C)=\pm 1$. Note that different choices of the ``plus'' sides $L^j_+$ in condition c) and the branches of the square roots in condition c) lead to isomorphic bundles. To complete the proof, it remains to note that the number of all possible choices of $\{\sigma(a_j|C),\sigma(b_j|C)\}_{j=1}^g$ is equal to $2^{2g}$ and coincides with the number of all non-equivalent spinor bundles.
\end{proof}

\subsubsection{Varying $A_j$, $B_j$} 
\label{par sur fam A}
\begin{figure}[h]
\center{\includegraphics[width=1\linewidth]{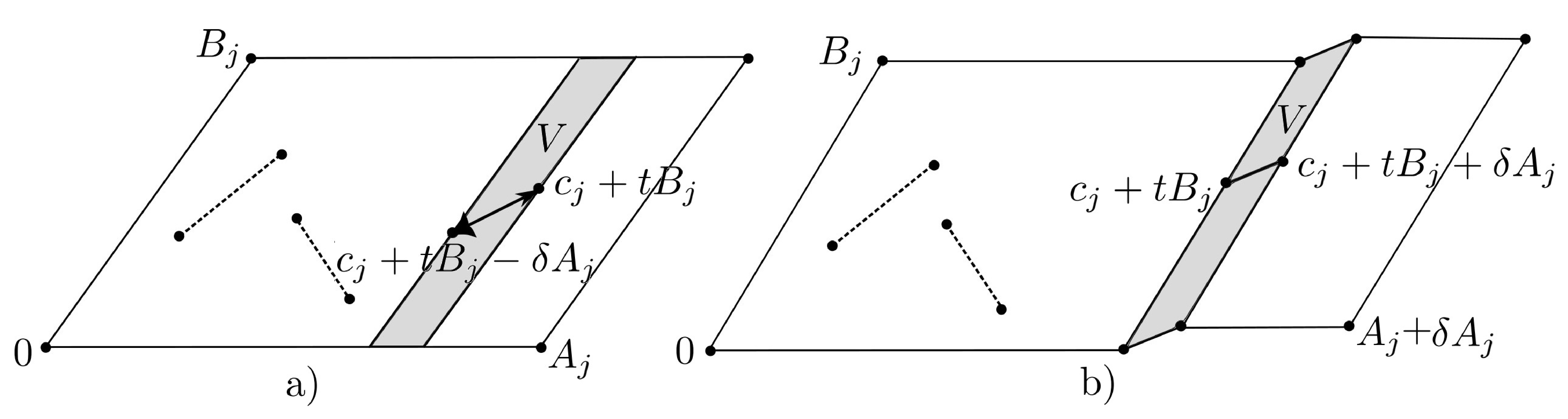}}
\caption{Deformation of the torus $\mathbb{T}_j$: a) the case $\alpha\le 0$ ($V$ is removed, the double-sided arrow denotes the identification of the points), b) the case $\alpha>0$ ($V$ is glued in). For simplicity, we assumed that $e_j=0$.}
\label{deform 1}
\end{figure}

In order to vary $A_j$ while keeping the other coordinates (\ref{moduli space coordinates}) unchanged, we implement the following deformation of $(X,\omega)$. Introduce a small variation $\delta A_j=\alpha A_j+\beta B_j$ of coordinate $A_j$, where $\alpha,\beta\in\mathbb{R}$. For $\alpha\le 0$, we remove the tubular neighborhood 
\begin{equation}
\label{removed tub nei}
V=\{[c_j+tB_j-s\delta A_j]({\rm mod}\,\Lambda_j) \ | \ t,s\in[0,1]\}
\end{equation}
of $b_j$ from the torus $\mathbb{T}_j$ of $X$. Then we identify the points
\begin{equation}
\label{identifying points alpha neg}
[c_j+tB_j+\delta A_j\,0] ({\rm mod}\,\Lambda_j)\longleftrightarrow [c_j+tB_j-\delta A_j\,(1+0)]({\rm mod}\,\Lambda_j)
\end{equation}
(see Fig. \ref{deform 1}, a)). As a result, we obtain the surface $X'$. The restriction of the $\omega$ to the domain $\acute{X}=X\backslash V$ is extended by continuity to the differential $\omega'$ on $X'$. All paths $a_i,b_i,l_k$, except $a_j$, belong to $\acute{X}$ and thus can be chosen as $a'_i,b'_i,l'_k$, respectively. At the same time, the path $a'_j$ is obtained by identification of the ends of the straight line segment $[\alpha A_j+c_j,A_j+\beta B_j+c_j]({\rm mod}\,\Lambda_j)$ in $\mathbb{T}_j\cap\acute{X}$.

For $\alpha>0$, we cut $X$ along $b_j$. Then we glue the resulting surface $\acute{X}$ and the cylinder
$$V=\{[c_j+tB_j+s\delta A_j]({\rm mod}\,B_j\mathbb{Z}) \ | \ t,s\in[0,1]\}$$
along their boundaries according to the following identification rules
\begin{align*}
\partial\acute{X}\ni [c_j+tB_j-\delta A_j\,0]({\rm mod}\,\Lambda_j)&\longleftrightarrow [c_j+tB_j]({\rm mod}\,B_j\mathbb{Z})\in\partial V, \\
\partial\acute{X}\ni [c_j+tB_j+\delta A_j\,0]({\rm mod}\,\Lambda_j)&\longleftrightarrow [c_j+tB_j+\delta A_j]({\rm mod}\,B_j\mathbb{Z})\in\partial V
\end{align*}
(see Fig. \ref{deform 1}, b)). The differential $\omega'$ on $X'$ is given by $\omega'=\omega$ on $\acute{X}$ and $\omega=dz$ on $V$. All paths $a'_i,b'_i,l'_k$, except $a_j$, coincide with $a_i,b_i,l_k$ while $a_j'$ is a union of straight line segments $[c_j,A_j+\beta(1+\alpha)^{-1} B_j]({\rm mod}\,\Lambda_j)$ and $[A_j+\beta(1+\alpha)^{-1} B_j,\delta A_j]({\rm mod}\,B_j\mathbb{Z})$ of $\acute{X}$ and $V$, respectively. 

The spinor bundle $C'$ over $X'$ is defined by conditions a)-c) with the automorphy factors 
\begin{equation}
\label{automorphy prime}
\sigma(a'_j|C')=\sigma(a_j|C), \quad \sigma(b'_j|C')=\sigma(b_j|C), \qquad (j=1,\dots,g).
\end{equation}
Note that, in both cases $\alpha\le 0$, $\alpha>0$, one can consider $\acute{X}$ as a common part of $X$ and $X'$. Moreover, the restrictions of the bundles $C$ and $C'$ on $\acute{X}$ are isomorphic,
\begin{equation}
\label{bundle restrictions}
C|_{\acute{X}}\equiv C'_{\acute{X}}.
\end{equation}
In other words, in the domain $\acute{X}$, one can add sections of $C'$ and sections of $C$. Note that a small shifts of the paths $\check{a}_j$, $\check{b}_j$, i.e., a small change of $c_j$ in (\ref{cycles shifted}), do not affect $(X',\omega',C')$.
\begin{lemma}
\label{spin deformation lemma}
The bundle $C'$ defined by {\rm(\ref{automorphy prime})} obeys {\rm (\ref{specifying spinor C prime})} for sufficiently small $\delta A_j$.
\end{lemma}
\begin{proof}
In view of the above procedure for constructing $(X',\omega',C')$, for small $\delta A_j$ there exists a near-isometric diffeomorphism $q: \ (X,|\omega|^2)\mapsto (X',|\omega'|^2)$ obeying $q\circ a_i=a'_i$, $q\circ b_i=b'_i$ ($i=1,\dots,g$) and such that $q$ is identity on the common tori $\mathbb{T}_i$ ($i\ne j$) of $X$ and $X'$ and near each $P_k$. For simplicity, one can consider $X'$ as a surface $X$ endowed with another complex structure provided by the operator $\overline{\partial'}$. Then $\beta$ is identity on $X_0$ and $\overline{\partial}=\overline{\partial'}$ on some domain $U\subset X_0$. Moreover, in each chart one has
\begin{equation}
\label{near conformal}
\overline{\partial'}=(1+a(z))\overline{\partial}+\tilde{a}(z)\partial, \text{ where } \|a\|_{C^l}+\|\tilde{a}\|_{C^l}\underset{\delta A_j\to 0}{\longrightarrow} 0.
\end{equation}
Sections $f$ of both $C$ and $C'$ can be considered as collections $\{f_i\}_{i=1}^{g},\{\tilde{f}_k\}_{k=1}^{2g-2}$ of functions $f_i$ on $\bar{\mathbb{T}}_i$ and $\tilde{f}_k$ on $\mathbb{D}_k$ obeying conditions a),b) given after (\ref{cycles explicit}). For $l=1,\dots$, denote by $\mathcal{H}^l$ the subspace in $H^{l}(\bar{\mathbb{T}}_1)\times\dots\times H^{l}(\bar{\mathbb{T}}_g)\times H^{l}(\mathbb{D}_1)\times\dots\times H^{l}(\mathbb{D}_{2g-1})$ consisting of collections obeying a)-c). Since $\overline{\partial}$ is an elliptic operator and $C$ admits no holomorphic sections, its inverse $\overline{\partial}^{-1}: \ \mathcal{H}^l\to \mathcal{H}^{l+1}$ is well defined and continuous. In view of (\ref{near conformal}), the same is true for the (elliptic) operator $\overline{\partial'}^{-1}$, and 
\begin{equation}
\label{d bars closeness}
\|\overline{\partial'}^{-1}-\overline{\partial}^{-1}\|_{\mathcal{H}^l\to \mathcal{H}^{l+1}}\underset{\delta A_j\to 0}{\longrightarrow} 0 \qquad (l=1,\dots).
\end{equation}
Let $z$ be a holomorphic (with respect to both complex structures $\overline{\partial}$ and $\overline{\partial'}$) coordinate on $U$. Let $z_0\in U$ and $\chi$ be a cut-off function equal to one near $z_0$ and vanishing outside $U$. Let $w(z)=(z-z_0)^{-1}$ on $U$. The Szego kernels of the bundles $C$ and $C'$ admit the representations (see \ref{par Szego})
$$\mathcal{S}(\cdot,z_0)=\chi \ w-\overline{\partial}^{-1}\Big([\overline{\partial},\chi]w\Big), \qquad \mathcal{S}'(\cdot,z_0)=\chi \ w-\overline{\partial'}^{-1}\Big([\overline{\partial'},\chi]w\Big).$$
In view of (\ref{near conformal}) and (\ref{d bars closeness}), we have 
$$\|\mathcal{S}'(\cdot,z_0)-\mathcal{S}(\cdot,z_0)\|_{\mathcal{H}^l}\underset{\delta A_j\to 0}{\longrightarrow} 0 \qquad (l=1,\dots).$$
Therefore, in  view of the Morrey's inequality, the zeroes $o'_1,\dots,o'_g$  of $\mathcal{S}'(\cdot,z_0)$ tend to the corresponding zeroes $o_1,\dots,o_g$ of $\mathcal{S}(\cdot,z_0)$,
\begin{equation}
\label{zeroes closeness}
o'_i\underset{\delta A_j\to 0}{\longrightarrow}o_i.
\end{equation}
Similarly, one can prove the convergences
\begin{equation}
\label{abel closeness}
\|\omega'-\omega\|_{C^{l}(X)}+\|\mathcal{A}'_{z_0}-\mathcal{A}_{z_0}\|_{C(X)}\underset{\delta A_j\to 0}{\longrightarrow} 0, \quad \mathbb{B}'\underset{\delta A_j\to 0}{\longrightarrow}\mathbb{B}, \quad K'_{z_0}\underset{\delta A_j\to 0}{\longrightarrow}K_{z_0}, 
\end{equation}
where $\{\upsilon'_s\}_{s=1}^g$ be a basis of Abelian differentials on $X\equiv (X_0,\overline{\partial})$ ($X'\equiv(X_0,\overline{\partial'})$) dual to the homology basis $\{a_i,b_i\}_{i=1}^g$, $\mathbb{B}'$ be the $b$-period matrix of $X'$, $K'_{z_0}$ be the vector of Riemann constants, and $\mathcal{A}'_{z_0}$ be the Abel map of $X'$. 

Note that $\mathcal{D}=o_1+\dots +o_g-z_0$ and $\mathcal{D}'=o'_1+\dots +o'_g-z_0$ are divisors of the bundles $C$ and $C'$, respectively. Since $C'$ is a spinor bundle, we have $$\mathcal{A}'_{z_0}(\mathcal{D}')=\mathcal{A}'_{z_0}(C')=-K'_{z_0}+\mathbb{B}'p'+q'$$
with $p',q'\in\{0,\frac{1}{2}\}^{2g}$. In view of convergences (\ref{zeroes closeness}) and (\ref{abel closeness}), the comparison of the last equation with (\ref{specifying spinor C}) yields $p'=p$, $q'=q$ for sufficiently small $\delta A_j$. Formula (\ref{specifying spinor C prime}) is proved.
\end{proof}
Variation of $B_j$ is performed similarly (one needs only to interchange $\check{a}_j, A_j$ and $\check{b}_j, B_j$ in the above construction).

\subsubsection{Varying $C_k$} 
\begin{figure}[h]
\center{\includegraphics[width=0.5\linewidth]{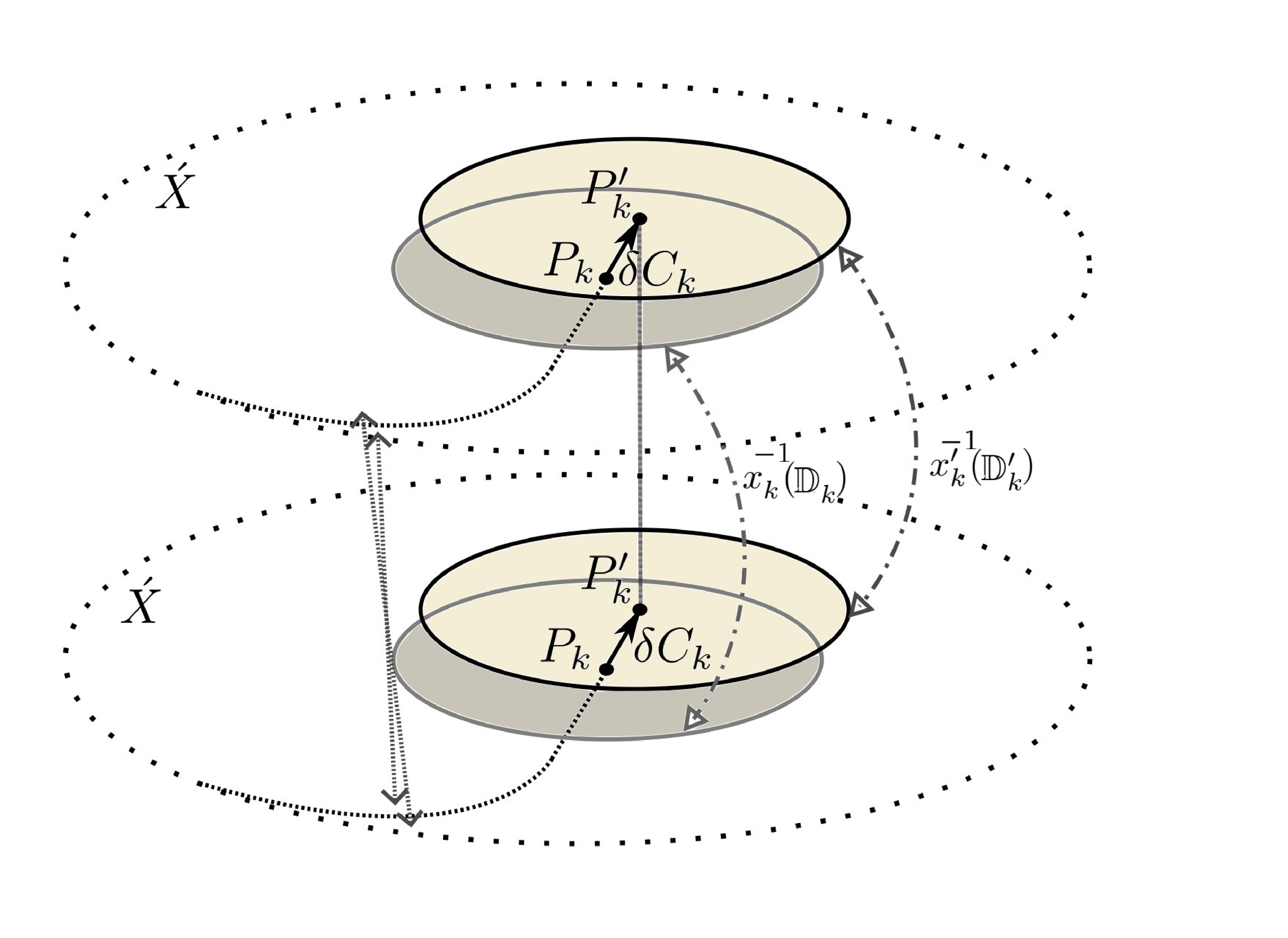}}
\caption{Deformation of the neighborhood of $P_k$ corresponding to the variation $\delta C_k$.}
\label{deform 2}
\end{figure}
In order to vary $C_k$ while keeping the other coordinates (\ref{moduli space coordinates}) unchanged, we implement the following deformation of $(X,\omega)$. Recall that the coordinates $z_k$ and $x_k=\sqrt{2z_k}$ near $P_k$ are given by (\ref{dist coord}) and (\ref{z coordinate}), respectively. One can deform $l_k$ (preserving its class in $H_1(X,\{P_1,\dots,P_{2g-2}\},\mathbb{Z})$) in such a way that it coincides with the ray $x_k/\sqrt{2\delta C_k}\ge 0$ near $P_k$. Remove from $X$ the neighborhood $\mathbb{D}_k=\{x_k\in X \ | \ |x_k|<\sqrt{2\epsilon_0}\}$ (with $\epsilon_0>|\delta C_k|$), then attach the disk $\mathbb{D}'_k:=\{x'_k\in\mathbb{C} \ | \ |x'_k|<\sqrt{2\epsilon_1}\}$ (with $\epsilon_1>\epsilon_0$) by identifying the points $x_k$ and $z_k\in X\backslash \mathbb{D}_k$ and $x'_k\in \mathbb{D}'_k$ such that
\begin{equation}
\label{identification new conical p}
x'_k=\sqrt{2(z_k-\delta C_k)}
\end{equation}
(see Fig. \ref{deform 2}). As a result, we obtain the surface $X'$. Then $\acute{X}:=X\backslash \mathbb{D}_k$ can be considered as a common part of $X$ and $X'$ while $\mathbb{D}'_k$ is considered as a domain in $X'$. The differential $\omega'$ on $X'$ is given by $\omega'=\omega$ on $\acute{X}'$ and by $\omega'=dz'_k$ on $\mathbb{D}'_k$, where $z_k=x_k^2/2$. Then the zero $P'_k$ of $\omega'$ is given by $z'_k=0$ while the other zeroes $P'_s=P_s$ ($s\ne k$) belong to $\acute{X}$. The paths $a'_i:=a_i$, $b'_i:=b_i$, $l'_s:=l_s$ ($s\ne k$) are contained in $\acute{X}$ while the path $l'_i$ coincides with $l_i$ on $\acute{X}$ and with the ray $x'_k/\sqrt{2\delta C_k}$ on $\mathbb{D}'_k$. 

The spinor bundle $C'$ on $X'$ is defined (up to isomorphism) by condition (\ref{bundle restrictions}). By repeating arguments used in the prof of Lemma \ref{spin deformation lemma}, one can verify the equality (\ref{specifying spinor C prime}). Note that $(X',\omega',C')$ does not depend on the choice of $\epsilon_0,\epsilon_1$ in the above procedure.

\subsection{Variation of resolvent kernels}
\label{sec variation of resolvent kernels}
Denote by $\mathscr{R}_{\lambda}=\mathscr{R}_{\lambda}(\cdot,\cdot|X,\omega,C)$ and $\allowbreak\mathscr{R}'_{\lambda}=\mathscr{R}_{\lambda}(\cdot,\cdot|X',\omega',C')$ the resolvent kernels of $\Delta'=\Delta'_S$ and $\Delta=\Delta_S$, respectively. In view of (\ref{bundle restrictions}), the difference
\begin{equation*}
%\label{difference}
u_\lambda(z,z')=\mathscr{R}_\lambda(z,z'|X',\omega',C')-\mathscr{R}_\lambda(z,z'|X,\omega,C)
\end{equation*}
is well defined on the joint domain $\acute{X}$ of $X$ and $X'$. In addition, the complements of $\acute{X}$ in $X$ and $X'$ shrink to conical points $P_k$ or curves $\check{a}_j$, $\check{b}_j$ as $\delta A_j,\delta B_j,\delta C_k\to 0$. Moreover, by changing the constant $c_j$ in (\ref{cycles shifted}), the curves $\check{a}_j$, $\check{b}_j$ can be slightly shifted along the torus $\mathbb{T}_j$ without affecting the resulting surface $(X',\omega',C')$. Due to these facts, one can define the derivatives 
$$\frac{\partial\mathscr{R}_{\lambda}(z,z')}{\partial\nu}, \quad \frac{\partial\mathscr{R}_{\lambda}(z,z')}{\partial\overline{\nu}} \qquad (z,z'\in\dot{X}, \ \ \nu=A_j,B_j,C_k).$$
In this subsection, we prove the following fact.
\begin{prop}
\label{variational formulas prop}
The resolvent kernel $\mathscr{R}_{\lambda}$  is differentiable with respect to $\Re\nu$, $\Im\nu$ {\rm(}$\nu=A_i,B_i,C_k${\rm)}, and
\begin{equation}
\label{variational formulas}
\begin{split}
2i\frac{\partial\mathscr{R}_{\lambda}(z,z')}{\partial\nu}=\int\limits_{\nu^\dag}\omega(\varsigma)^{-1}\Big(\lambda \mathscr{R}_\lambda(z,\varsigma)&\mathscr{R}_\lambda(\varsigma,z')\rho^{-2}(\varsigma)h(\varsigma)d\overline{\varsigma}-\\
-&\partial_{\varsigma}\mathscr{R}_\lambda(z,\varsigma)\partial_{\varsigma}(h(\varsigma)\mathscr{R}_\lambda(\varsigma,z'))d\varsigma\Big).
\end{split}
\end{equation}
Here $A_i^\dag:=-\check{b}_i$, $B_i^\dag:=\check{a}_i$, and $C_k^\dag=c_k$ is homologous to the small {\rm(}counterclockwise{\rm)} circle $|\xi_k|=\epsilon$ in $\dot{X}$.
\end{prop}
Before proving Proposition \ref{variational formulas prop}, let us note that the integrand in the right hand side of (\ref{variational formulas}) is a closed 1-form and, thus, the right hand side of (\ref{variational formulas}) depends only on the homology class of $\nu^\dag$ in $H(\dot{X},\mathbb{Z})$.

\subsubsection{Boundedness of $u_\lambda(z,z')$} 
\label{par boundedness}
For simplicity, we consider the case in which only the coordinate $\nu=A_j$ varies (the other cases are considered similarly). Suppose that $\lambda\in\mathbb{C}$ is not an eigenvalue of $\Delta=\Delta_S$. First, let us show that $\lambda$ is separated from the spectra of $\Delta'$ for sufficiently small variations. Suppose the contrary; then there exists a sequence of the triples $(X'_{(n)},\omega'_{(n)},C'_{(n)})$ constructed in paragraph \ref{par sur fam A} and representing the variations $\delta A_{j}^{(n)}\to 0$ of $A_j$, and a sequence of eiqenpairs $(\lambda'_{n},u'_{n})$ of the (Szeg\"o) Laplacians $\Delta'_{(n)}$ on $(X'_{(n)},\omega'_{(n)},C'_{(n)})$ such that $\lambda'_{n}\to \lambda$ and $\|u'_{n}\|_{L_2(X'_{(n)},C'_{(n)})}=1$.

Let $\epsilon>0$ be sufficiently small and let
$$X_{\epsilon}^1=\Big\{[(t+is)B_j+c_j]({\rm mod}\,B_j\mathbb{Z}) \ \big| \ t\in [0,1], \ s\in\Big(-\frac{\epsilon}{|B_j|},\frac{\epsilon}{|B_j|}\Big)\Big\}$$
be the $\epsilon_0$-neighborhood of the path $\check{b}_{j}$ in $(X,|\omega|^2)$. Let also $X_{\epsilon}^0=X\backslash X_{\epsilon/2}^{1}$. For small $\delta A_{j}^{(n)}$, one can consider $X_0\subset\acute{X}$ and $\check{b}_{j}\subset\overline{\acute{X}}$ as a domain and a curve in $X'_{(n)}$. The (closed) $\epsilon_0$-neighborhood of the path $\check{b}_{j}$ in $(X'_{(n)},|\omega'_{(n)}|^2)$ is isometric to the cylinder $X_{\epsilon}^1$. Then $X'_{(n)}$ can be represented as the union of $X_{\epsilon}^0$ and $X_{\epsilon}^1$ where a common point of $X_{\epsilon}^0$ and $X_{\epsilon}^1$ is identified with points
\begin{equation}
\label{idetification points rules}
[(t\pm is)B_j+c_j]({\rm mod}\,\Lambda_j) \ \longleftrightarrow \ [(t\pm is)B_j\mp\delta A_{j}^{(n)}/2+c_j]({\rm mod}\,B_j\mathbb{Z})
\end{equation}
of $\mathbb{T}_j\cap X_{\epsilon}^0$ and $X_{\epsilon}^1$, respectively.
%(where $s|B_j|\in[\epsilon,2\epsilon]$) represent the same point of $X'_{(n)}$ which belongs to both 

Note that $C'_{(n)}|_{X_{\epsilon}^m}$ is isomorphic to $C|_{X_{\epsilon}^m}$ due to (\ref{automorphy prime}). Let $\mathcal{H}^{l}_{\nu}(X_{\epsilon}^m)$ ($m=0,1$) be the space of sections of $C|_{X_{\epsilon}^m}$ with finite norms (\ref{weighted spaces}), where $X$ is replaced by $X_{\epsilon}^m$. The convergence 
$$\|(\Delta'_{(n)})^q u'_{n}\|_{L_2(X'_{(n)};C'_{(n)})}=\lambda^q_n\to \lambda^q, \qquad q=0,1,\dots$$ 
and local estimates (\ref{H l smooth incr}), (\ref{smoth inc weight}) imply that the norms $\|u'_n|_{X_{\epsilon}^m}\|_{\mathcal{H}^{l'}_{\nu'}(X_{\epsilon}^m)}$ are bounded for some $\nu'$ and (arbitrarily large) $l'$. The compactness of the embedding $\mathcal{H}^{l'}_{\nu'}(X_{\epsilon}^m)\subset\mathcal{H}^l_\nu(X_{\epsilon}^m)$ (with $\nu>\nu'$ and $l'<l$) can be proved in the same way as in \ref{par discretness spectra}. As a corollary, there is a sub-sequence $\{u_{n(s)}\}_{s=1}^{\infty}$ obeying
\begin{equation}
\label{converg from compac}
u_{n(s)}|_{X_{\epsilon}^m}\to y_m \text{ in } \mathcal{H}^l_\nu(X_{\epsilon}^m) \quad (m=0,1).
\end{equation} In view of (\ref{idetification points rules}) and (\ref{automorphy prime}), $y_0,y_1$ are restrictions on $X_{\epsilon}^0,X_{\epsilon}^1$ of some section $y$ of $C$. Then (\ref{converg from compac}) and the equations $\Delta'_{(n)} u_n=\lambda'_n u_n$ imply $\Delta y=\lambda y$ on $\dot{X}$. In addition, $u_{n(s)}$ admit asymptotics (\ref{Szego domain}) near each $O_k\in\acute{X}$. The same is true for $y$ due to (\ref{converg from compac}) and Proposition \ref{prop elliptic asymptotics} (with $X$ replaced by $X_{\epsilon}^m$). Thus, $\Delta_S y=\lambda y$, which gives a contradiction. Therefore, $\lambda$ is separated from the spectra of $\Delta'$ for sufficiently small $\delta A_j$.

Next, let $\delta A_j$ be sufficiently small and let $\dot{X}^0_\epsilon$ be a domain in $X$ obtained by removing the $\epsilon$-neighborhoods (in the metric $|\omega|^2$) of $\check{b}_j$ and conical points. Let us show that
\begin{equation}
\label{boundedness}
|u_{\lambda}(z,z')|+|\partial_{\overline{z}}u_{\lambda}(z,z')|\le C \qquad (z\in\overline{\acute{X}}, \ z'\in\dot{X}^0_\epsilon) 
\end{equation}
uniformly in $(X',\omega')$, where the coordinates $z,z'$ are given by (\ref{z coordinate}). Let $\chi$ be a cut-off function equal to one on $X\backslash X^1_{\epsilon/2}$ and equal to zero outside of $X\backslash X^1_{\epsilon/4}$. Represent $u_\lambda$ as the sum $u_{1,\lambda}+u_{2,\lambda}$, where 
$$u_{1,\lambda}(\cdot,z')=(\chi-1)\mathscr{R}_{\lambda}(\cdot,z'), \qquad u_{2,\lambda}(\cdot,z')=\mathscr{R}'_{\lambda}(\cdot,z')-\chi \mathscr{R}_{\lambda}(\cdot,z').$$
Since $\mathscr{R}_{\lambda}(z,z')$ and  $\mathscr{R}'_{\lambda}(z,z')$ are smooth outside conical points and have the same logarithmic singularity (given by (\ref{log singularity})) at the diagonal $z=z'$, the section $u_{2,\lambda}(\cdot,z')$ is smooth for $z'\in \dot{X}^0_\epsilon$. Note that $(\Delta-\lambda)u_{2,\lambda}(\cdot,z')=-[\Delta,\chi]\mathscr{R}_{\lambda}(\cdot,z')$ does not depend on $(X',\omega')$. Since the spectra of $\Delta'$ are separated from $\lambda$, the sum $\sum_{s=1}^l\|\Delta^s u_{2,\lambda}(\cdot,z')\|_{L_{2}(X';C)}$ is bounded uniformly in $z'\in \dot{X}^0_\epsilon$ and $(X',\omega')$. Thus, (\ref{boundedness}) follows from (\ref{H l smooth incr}) and (\ref{Morrey}).

\subsubsection{Varying $A_j,B_j$}
For simplicity, we consider the case in which only the coordinate $\nu=A_j$ varies, and $\delta A_j=\alpha A_j+\beta B_j$, $\alpha\le 0$ (the other cases are considered similarly). In this case $X'$ is obtained by identification (\ref{identifying points alpha neg}) of the boundary points of $\acute{X}=X\backslash V$, where $V$ is given by (\ref{removed tub nei}). 

Recall that the resolvent kernels $\mathscr{R}_{\lambda}(z,z')$, $\mathscr{R}'_{\lambda}(z,z')$ are smooth outside conical points and have the same logarithmic singularity (given by (\ref{log singularity})) at the diagonal $z=z'$. So, their difference $u_\lambda(z,z')$ is smooth and obeys $(\Delta-\lambda)u_\lambda(\cdot,z')=0$ on $\acute{X}\cap\dot{X}$. Since $\mathscr{R}_{\lambda}(z,z')$, $\mathscr{R}'_{\lambda}(z,z')$ correspond to the Szeg\"o extensions of the Laplacians on $(X,\omega,C)$ and $(X',\omega',C')$, their difference $u_\lambda(z,z')$ admits the asymptotics of the form (\ref{Szego domain}) near conical points. Thus, Green formula (\ref{Green formula}) with $f=u_\lambda(\cdot,z')$, $f'=\mathscr{R}_{\overline{\lambda}}(\cdot,z)=\overline{\mathscr{R}_{\lambda}(z,\cdot)}$, and $U=\acute{X}$ yields
\begin{equation}
\label{resolvent on boundary data}
\begin{split}
0=((\Delta-\lambda)u_\lambda(\cdot,z'),&\mathscr{R}_{\overline{\lambda}}(\cdot,z))_{L_2(\breve{X};S)}-(u_\lambda(\cdot,z'), (\Delta-\overline{\lambda})\mathscr{R}_{\overline{\lambda}}(\cdot,z))_{L_2(\breve{X};S)}=\\
=u_\lambda(z,z')+\int\limits_{\partial\acute{X}}&\frac{h(\varsigma)}{2i}\big(u_\lambda(\varsigma,z')\partial_\varsigma\mathscr{R}_{\lambda}(z,\varsigma)d\varsigma+\partial_{\overline{\varsigma}}u_\lambda(\varsigma,z')\mathscr{R}_{\lambda}(z,\varsigma)d\overline{\varsigma}\big).
\end{split}
\end{equation}
Here $\partial\acute{X}$ is a union of $\check{b}_j=[c_j,B_j+c_j]({\rm mod}\,\Lambda_j)$ and the curve $[B_j-\delta A_j+c_j,-\delta A_j+c_j]({\rm mod}\,\Lambda_j)$ obtained by shifting $\check{b}_j$ (along the torus $\mathbb{T}_j$) and reversing its orientaiton. Thus,(\ref{resolvent on boundary data}) can be rewritten as 
\begin{equation}
\label{resolvent on boundary data 1}
\begin{split}
u_\lambda(z,z')=-\int\limits_{\sigma\in\check{b}_j}\frac{h(\varsigma)}{2i}\big(u_\lambda(\varsigma,z')\partial_\varsigma\mathscr{R}_{\lambda}(z,\varsigma)d\varsigma+\partial_{\overline{\varsigma}}u_\lambda(\varsigma,z')\mathscr{R}_{\lambda}(z,\varsigma)d\overline{\varsigma}\big)\Big|_{\varsigma=\sigma-\delta A_j}^{\varsigma=\sigma},
\end{split}
\end{equation}
where $\sigma,\zeta$ are given by (\ref{z coordinate}) while $h(\varsigma)=1$. Since the bundle $C'$ obeys (\ref{bundle restrictions}) while its section $z\mapsto\mathscr{R}'_{\lambda}(z,z')$ is smooth outside $z=z'$ and conical points, we have
$$\mathscr{R}'_{\lambda}(\varsigma,z')\Big|_{\varsigma=\sigma-\delta A_j}^{\zeta=\sigma}=\partial_{\overline{\zeta}}\mathscr{R}'_{\lambda}(\varsigma,z')\Big|_{\varsigma=\sigma-\delta A_j}^{\zeta=\sigma}=0$$
and
\begin{align*}
u_\lambda(\varsigma,z')\Big|_{\varsigma=\sigma-\delta A_j}^{\varsigma=\sigma}&=-\mathscr{R}_{\lambda}(\varsigma,z')\Big|_{\varsigma=\sigma-\delta A_j}^{\varsigma=\sigma}, \\ \partial_{\overline{\varsigma}}u_\lambda(\varsigma,z')\Big|_{\varsigma=\sigma-\delta A_j}^{\varsigma=\sigma}&=-\partial_{\overline{\varsigma}}\mathscr{R}_{\lambda}(\varsigma,z')\Big|_{\varsigma=\sigma-\delta A_j}^{\varsigma=\sigma}.
\end{align*}
Since $\mathscr{R}_{\lambda}(\cdot,z')$ is smooth outside $z'=z$ on $\dot{X}$, the equalities 
\begin{align*}
\mathscr{R}_{\lambda}(\varsigma,z')\Big|_{\varsigma=\sigma-\delta A_j}^{\varsigma=\sigma}=\frac{1}{\omega(\sigma)}\big[h^{-1}(\sigma)\partial_{\sigma}h(\sigma)\big]\mathscr{R}_{\lambda}(\sigma,z')\delta A_j+\frac{1}{\overline{\omega(\sigma)}}\partial_{\overline{\sigma}}\mathscr{R}_{\lambda}(\sigma,z')\overline{\delta A_j},\\
\mathscr{R}_{\lambda}(z,\varsigma)\Big|_{\varsigma=\sigma-\delta A_j}^{\varsigma=\sigma}=\frac{1}{\omega(\sigma)}\partial_\sigma\mathscr{R}_{\lambda}(z,\sigma)\delta A_j+\frac{1}{\overline{\omega(\sigma)}}\big[h^{-1}(\sigma)\partial_{\overline{\sigma}}h(\sigma)\big]\mathscr{R}_{\lambda}(z,\sigma)\overline{\delta A_j}
\end{align*}
and
\begin{align*}
\partial_{\overline{\varsigma}}\mathscr{R}_{\lambda}(\varsigma,z')\Big|_{\varsigma=\sigma-\delta A_j}^{\varsigma=\sigma}&=\frac{1}{\omega(\sigma)}\big[h^{-1}(\sigma)\partial_\sigma h(\sigma)\partial_{\overline{\sigma}}\big]\mathscr{R}_{\lambda}(\sigma,z')\delta A_j+\\
+\big[\partial_{\overline{\sigma}}\frac{1}{\overline{\omega(\sigma)}}\partial_{\overline{\sigma}}\big]\mathscr{R}_{\lambda}&(\sigma,z')\overline{\delta A_j}=\Big[h^{-1}(\sigma)\partial_\sigma h(\sigma)\partial_{\overline{\sigma}}=
-\rho^{-2}(\sigma)\Delta_\sigma\Big]=\\
=&-\lambda\frac{\rho^{-2}(\sigma)}{\omega(\sigma)}\mathscr{R}_{\lambda}(\sigma,z')\delta A_j+\big[\partial_{\overline{\sigma}}\frac{1}{\overline{\omega(\sigma)}}\partial_{\overline{\sigma}}\big]\mathscr{R}_{\lambda}(\sigma,z')\overline{\delta A_j},\\
\partial_{\varsigma}\mathscr{R}_{\lambda}(z,\varsigma)\Big|_{\varsigma=\sigma-\delta A_j}^{\varsigma=\sigma}&=\big[\partial_{\sigma}\frac{1}{\omega(\sigma)}\partial_{\sigma}\big]\mathscr{R}_{\lambda}(z,\sigma)\delta A_j+\\
+&\frac{1}{\overline{\omega(\sigma)}}\big[h^{-1}(\sigma)\partial_{\overline{\sigma}}h(\sigma)\partial_{\sigma}\big]\mathscr{R}_{\lambda}(z,\sigma)\overline{\delta A_j}
\end{align*}
hold up to $o(|\delta A_j|)$-terms uniformly in $z'\in\acute{X}$ separated from $\check{b}_j$. Substituting the above formulas into (\ref{resolvent on boundary data 1}) and taking into account the identity
$$(ab)\big|_-^+=a(+)(b\big|_-^+)+(a\big|_-^+)b(-)=a(+)(b\big|_-^+)+(a\big|_-^+)b(+)-(a\big|_-^+)(b\big|_-^+),$$
we obtain 
\begin{equation}
\label{complete variation of resolvent kernel}
\begin{split}
u_\lambda(z,z')=\frac{\delta A_j}{2i}\int\limits_{\sigma\in -\check{b}_j}&\frac{h(\sigma)}{\omega(\sigma)}\Big[\lambda\rho^{-2}(\sigma)\mathscr{R}_{\lambda}(z,\sigma)\mathscr{R}_{\lambda}(\sigma,z')d\overline{\sigma}-\\
-&\partial_\varsigma\mathscr{R}_{\lambda}(z,\sigma)\big[h^{-1}(\sigma)\partial_{\sigma}h(\sigma)\big]\mathscr{R}_{\lambda}(\sigma,z')d\sigma\Big]+\\
+\frac{\overline{\delta A_j}}{2i}\int\limits_{\sigma\in -\check{b}_j}\frac{h(\sigma)}{\overline{\omega(\sigma)}}&\Big[\partial_{\overline{\sigma}}\mathscr{R}_{\lambda}(\sigma,z')\partial_\varsigma\mathscr{R}_{\lambda}(z,\sigma)d\sigma-\\
&-\big[\overline{\omega(\sigma)}\partial_{\overline{\sigma}}\frac{1}{\overline{\omega(\sigma)}}\partial_{\overline{\sigma}}\big]\mathscr{R}_{\lambda}(\sigma,z'))\mathscr{R}_{\lambda}(z,\sigma)d\overline{\sigma}\Big]+\\
&+O\big(|\delta A_j|(\mathcal{C}(u_\lambda)+|\delta A_j|)\big),
\end{split}
\end{equation}
where 
$$\mathcal{C}(u_\lambda):=\max_{\varsigma\in \check{b}_j}(|u(\varsigma,z')|+|\partial_{\overline{\varsigma}}u(\varsigma,z')|).$$
In view of (\ref{boundedness}), we have $\mathcal{C}(u_\lambda)\le C$ uniformly in $(X',\omega')$. Then, formula (\ref{complete variation of resolvent kernel}) implies $u_\lambda(z,z')=O(|\delta A_j|)$ uniformly in $z,z'\in\acute{X}$ separated from $\check{b}_j$ and conical points. The same estimate for the derivatives of $u_\lambda(z,z')$ with respect to $z$, $\overline{z}$ follows from the equation $\Delta u_\lambda(\cdot,z')=\lambda u_\lambda(\cdot,z')$ in $\acute{X}$ and local estimates (\ref{H l smooth incr}), (\ref{Morrey}). As mentioned after (\ref{bundle restrictions}), a slight shift of the path $\check{b}_j$ (by small change of $c_j$ in (\ref{cycles shifted})) in the construction procedure for $(X',\omega',C')$ does not affect $(X',\omega',C')$. Therefore, in the above estimates, one can omit the condition that $z,z'$ are separated from $\check{b}_j$. As a corollary, one has $\mathcal{C}(u_\lambda)=O(|\delta A_j|)$. Now, formula (\ref{complete variation of resolvent kernel}) implies that $\mathscr{R}_\lambda(z,z'|(X,\omega))$ is differentiable with respect to $A_j$ and $\overline{A_j}$, and equality (\ref{variational formulas}) is true with $\nu=A_j$.

\subsubsection{Varying $C_k$} 
In this case, the surface $X'$ is obtained by identification (\ref{identification new conical p}) of points of $\acute{X}=X\backslash \mathbb{D}_k$ and $\mathbb{D}'_k$, where $\mathbb{D}_k=\{z_k\in X \ | \ |z_k|<\epsilon_0\}$ and $\mathbb{D}'_k:=\{x'_k\in\mathbb{C} \ | \ |x'_k|<\sqrt{2\epsilon_1}\}$. The choice of (sufficiently small) $\epsilon_0$, $\epsilon_1$ does not affect $(X',\omega',C')$ until $\epsilon_1>\epsilon_0>|\delta C_k|$. In what follows, we choose $\epsilon_0=2|\delta C_k|^{p_0}$, where $p_0\in(0,1]$.

Denote $z'_k=x'^2_k/2$. Let $\chi'_k$ be a cut-off function in $\mathbb{D}'_k$ which is equal one for $|z'_k|\le\epsilon_1/3$ and equal to zero for $|z'_k|\ge 2\epsilon_1/3$. Put
\begin{equation}
\label{equation near O k prime}
\tilde{f}'_\lambda(\cdot,z'):=(\Delta-\lambda)\big(\chi'_k\mathscr{R}'_\lambda(\cdot,z'))=[\Delta,\chi'_k]\mathscr{R}'_\lambda(\cdot,z').
\end{equation}
Let the coordinates $z,z'$ on $\acute{X}$ be given by (\ref{z coordinate}). Suppose that $z$ and $z'$ are separated from each other and the conical points $P_1,\dots,P_{2g-2}$. In what follows, all estimates and asymptotics are assumed to be uniform in $z,z'$ and $(X',\omega')$. By repeating the reasoning of \ref{par boundedness}, one can show that the resolvent kernel $\mathscr{R}'_\lambda(z,z')=\mathscr{R}_\lambda(z,z')+u_\lambda(z,z')$ and its derivatives with respect to $z$ are bounded uniformly in $(X',\omega')$. In particular, $\tilde{f}'_\lambda(z'_k,z')$ and its derivatives with respect to $z'_k$ are bounded uniformly in $z'$ and $(X',\omega')$. In addition, we have $|z'_k|\in[\epsilon_0,2\epsilon_0]$ on the support of $\tilde{f}'_\lambda(\cdot,z')$. Now, applying Proposition \ref{prop elliptic asymptotics} (with $X$ replaced by $\mathbb{D}'_k$) yields the asymptotics
\begin{equation}
\label{resolvent kernel near O k}
\mathscr{R}'_\lambda(z'_k,z')=\sum_{(m,\pm)}c'^{\lambda}_{k,m,\pm}(z')f^{(k)}_{k,m,\pm}(r'_k,\varphi'_k)\Big(\sum_{n\ge 0} d(n,\pm i\mu_m)(\lambda|z'_k|^{2})^n\Big)+O(r_k'^{N}).
\end{equation}
Here $(r'_k,\varphi'_k)$ are polar coordinates near $P'_k$ (i.e., $z'_k=r'_k e^{i\varphi'_k}$) while $f^{(k)}_{k,m,\pm}$ are given by (\ref{anzats in z}). Asymptotics (\ref{resolvent kernel near O k}) admits differentiation with respect to $z'_k$ (which means that the remainder $\tilde{\mathscr{R}}'_\lambda(z'_k,z')=O(r_k'^{N})$ obeys $\partial^l_{z'_k}\tilde{\mathscr{R}}'_\lambda(z'_k,z'),\partial^l_{\overline{z'_k}}\tilde{\mathscr{R}}'_\lambda(z'_k,z')=O(r_k'^{N-l})$ for $l=1,2,\dots$). The coefficients in (\ref{resolvent kernel near O k}) are given by
\begin{equation}
\label{resolvent kernel near O k coef}
c'^{\lambda}_{k,m,\pm}(z')=(\tilde{f}'_\lambda(\cdot,z'),\chi'_k f_{k,m,\mp})_{L_2(\mathbb{D}'_k;C)}-(\mathscr{R}'_\lambda(\cdot'_k,z'),(\Delta-\overline{\lambda})[\chi'_k f_{k,m,\mp}])_{L_2(\mathbb{D}'_k;C)}.
\end{equation}
The number $N$ in (\ref{resolvent kernel near O k}) can be taken arbitrarily large while the sum in the right hand side contains only the terms decreasing slower than $O(r_k'^{N})$. Since the resolvent kernel $\mathscr{R}'_\lambda$ corresponds to the Szeg\"o extension of the Laplacian on $(X',\omega',C')$, the first four terms in (\ref{resolvent kernel near O k}) are given by
$$\mathscr{R}'_\lambda(z'_k,z')=c'^{\lambda}_{k,0,+}f^{(k)}_{k,0,+}+c'^{\lambda}_{k,1,+}f^{(k)}_{k,1,+}+c'^{\lambda}_{k,2,+}f^{(k)}_{k,2,+}+c'^{\lambda}_{k,-1,-}f^{(k)}_{k,-1,-}+O(|z'_k|^{5/4}).$$
The same formulas with omitted primes are valid for the resolvent kernel $\mathscr{R}_\lambda(\cdot,z)$ of the Szeg\"o Laplacian on $(X,\omega,C)$. 

Green formula (\ref{Green formula}) with $f=u_\lambda(\cdot,z')$, $f'=\mathscr{R}_{\overline{\lambda}}(\cdot,z)=\overline{\mathscr{R}_{\lambda}(z,\cdot)}$, and $U=\acute{X}$ yields
\begin{equation}
\label{resolvent on boundary data 2}
u_\lambda(z,z')=\int\limits_{\partial\acute{X}}\frac{h(\varsigma)}{2i}\big(u_\lambda(z_k,z')\partial_{z_k}\mathscr{R}_{\lambda}(z,z_k)dz_k+\partial_{\overline{z_k}}u_\lambda(z_k,z')\mathscr{R}_{\lambda}(z,z_k)d\overline{z_k}\big).
\end{equation}
Here $\partial\acute{X}$ is the (counterclockwise oriented) circle $|z_k|=2|\delta C_k|^{p_0}$. In view of (\ref{resolvent kernel near O k}) and (\ref{identification new conical p}), the integrand in the right hand side of (\ref{resolvent on boundary data 2}) is $O(|z_k|^{-1/2})$. Then formula (\ref{resolvent on boundary data 2}) with $p_0=1$, the equality $\Delta u_\lambda(\cdot,z')=\lambda u_\lambda(\cdot,z')$, and local estimates (\ref{H l smooth incr}), (\ref{Morrey}) imply
\begin{equation}
\label{C k variation auxillary estimate 1}
|\partial^{l_1}_z\partial_{\overline{z}}^{l_2} u_\lambda(z,z')|=O(|\delta C_k|^{\rho}),
\end{equation}
where $l_1,l_2=0,1,\dots$ and $\rho=1/2$. In view of (\ref{C k variation auxillary estimate 1}), (\ref{identification new conical p}), we have
$$|\tilde{f}'_\lambda(z'_k,z')-\tilde{f}_\lambda(z_k,z')|=O(|\delta C_k|^{\rho})$$
where $\tilde{f}'_\lambda$ is given by (\ref{equation near O k prime}) and $\tilde{f}'_\lambda(z_k,z')=[\Delta,\chi'_k(z_k)]\mathscr{R}_\lambda(z_k,z')$. Thus, (\ref{resolvent kernel near O k coef}) yields
\begin{equation}
\label{C k variation auxillary estimate 2} 
|c'^{\lambda}_{k,m,\pm}(z')-c^{\lambda}_{k,m,\pm}(z')|=O(|\delta C_k|^{\rho}),
\end{equation}
for the coefficients in asymptotics (\ref{resolvent kernel near O k}) for $\mathscr{R}'_\lambda(z'_k,z')$ and $\mathscr{R}_\lambda(z_k,z')$.  Since 
$$z'^\alpha_k-z_k^{\alpha}=z_k^\alpha\big((1-\delta C_k/z_k)^\alpha-1\big)=-\alpha\delta C_k z_k^{\alpha-1}+O(|\delta C_k|^2|z_k|^{\alpha-2})$$
as $|\delta C_k|/r_k\to 0$, we have
\begin{equation}
\label{C_k variation difference asymp 2}
\begin{split}
f^{(k)}_{k,m,\pm}(r'_k,\varphi'_k)-f^{(k)}_{k,m,\pm}(r_k,\varphi_k)=&\\
=\big(1\mp(1/2+m)\big)(\Re\delta C_k\pm i\Im &\delta C_k) f^{(k)}_{k,m\mp 2,\pm}(r_k,\varphi_k)+O(|\delta C_k|^2 r_k^{-2\pm(2m-1)/4}).
\end{split}
\end{equation}
Combining (\ref{resolvent kernel near O k}), (\ref{C k variation auxillary estimate 2}), (\ref{C_k variation difference asymp 2}), we obtain
\begin{align*}
u_\lambda(z'_k,z')=\sum_{(m,\pm)}\big(1\mp(1/&2+m)\big)c^{\lambda}_{k,m,\pm}(z')f^{(k)}_{k,m\mp 2,\pm}(z_k)(\Re\delta C_k\pm i\Im\delta C_k)+\\
&+O(|\delta C_k|^2r_k^{-\frac{9}{4}})+O(r_k^{N})+O(|\delta C_k|^\rho r_k^{-\frac{1}{4}}), \\
\partial_{\overline{z'_k}}u_\lambda(z'_k,z')=\sum_{m}(3/2+m)&c^{\lambda}_{k,m,-}(z')\partial_{\overline{z'_k}}f^{(k)}_{k,m+2,\pm}(z_k)(\Re\delta C_k-i\Im\delta C_k)+\\
&+O(|\delta C_k|^2|r_k^{-\frac{9}{4}}|)+O(r_k^{N-1})+O(|\delta C_k|^\rho r_k^{-\frac{1}{4}}).
\end{align*}
Moreover, expansion (\ref{resolvent kernel near O k}) with omitted primes implies $\partial_{z_k}\mathscr{R}_{\lambda}(z_k,z')=O(r_k^{-1/4})$. In view of the above asymptotics and the equality
\begin{equation*}
%\label{pairing of asymptotic terms}
\frac{1}{2i}\int_{r_k=c}\big(f_{k,m,s}h\partial\overline{f_{k,m',s'}}dz+\overline{f_{k,m',s'}}h\overline{\partial}f_{k,m,s}d\overline{z}\big)=\delta_{m,m'}\delta_{s,-s'},
\end{equation*}
formula (\ref{resolvent on boundary data}) implies
\begin{equation}
\label{variation in C k in terms of coefficients}
\begin{split}
u_\lambda(z,z')=\frac{1}{2}&\Big(c^{\lambda}_{k,1,+}(z')\overline{c^{\overline{\lambda}}_{k,-1,-}(z)}-c^{\lambda}_{k,0,+}(z')\overline{c^{\overline{\lambda}}_{k,-2,-}(z)}\Big)\delta C_k +\\
+\frac{1}{2}&\Big(c^{\lambda}_{k,-2,-}(z')\overline{c^{\overline{\lambda}}_{k,0,+}(z)}-c^{\lambda}_{k,-1,-}(z')\overline{c^{\overline{\lambda}}_{k,1,+}(z)}\Big)\overline{\delta C_k}+O(|\delta C_k|^{\tilde{\rho}}),
\end{split}
\end{equation}
where $\tilde{\rho}:=\min\{2-3p_0/2,1+3p_0/2,5p_0/2,\rho+p_0/2\}$. The choice $p_0=2/3$ in (\ref{variation in C k in terms of coefficients}) provides estimate (\ref{C k variation auxillary estimate 1}) with $\rho=5/6$. Due to this, one can choose $\rho=5/6$ and $p_0=1/2$ in the above formulas; then $\tilde{\rho}=13/12$ and the remainder in (\ref{variation in C k in terms of coefficients}) is $o(|\delta C_k|)$. Thus, formula (\ref{variation in C k in terms of coefficients}) provides differentiability of $\mathscr{R}_\lambda(z,z'|(X,\omega))$ with respect to $\Re C_k$, $\Im C_k$ as well as the equality 
\begin{equation}
\label{differentation in C k}
2i\frac{\partial\mathscr{R}_\lambda(z,z')}{\partial C_k}=i\Big(c^{\lambda}_{k,1,+}(z')\overline{c^{\overline{\lambda}}_{k,-1,-}(z)}-c^{\lambda}_{k,0,+}(z')\overline{c^{\overline{\lambda}}_{k,-2,-}(z)}\Big).
\end{equation}
At the same time, the substitution of asymptotics (\ref{resolvent kernel near O k}) into the right hand side of (\ref{variational formulas}) with $\nu=C_k$ and $\nu^\dag=\{\zeta_k(t)=\epsilon e^{it}\}_{t\in[0,4\pi)}$, $\epsilon\to 0$ yields the same value as in the right hand side of (\ref{differentation in C k}). Thus, formula (\ref{variational formulas}) is proved for $\nu=C_k$.

\

So, Proposition \ref{variational formulas prop} is proved. Moreover, we have proved that the asymptotics
\begin{equation}
\label{uniform differentablity}
\mathscr{R}_\lambda(z,z'|X',\omega',C')-\mathscr{R}_\lambda(z,z'|X,\omega,C)=\frac{\partial\mathscr{R}_{\lambda}(z,z')}{\partial\nu}\delta\nu+\frac{\partial\mathscr{R}_{\lambda}(z,z')}{\partial\overline{\nu}}\overline{\delta\nu}+o(|\delta\nu|)
\end{equation}
(with $\nu=A_i,B_0,C_k$ and with $\partial\mathscr{R}_{\lambda}(z,z')/\partial\nu$ given by (\ref{variational formulas})) is uniform in $(X',\omega')$ and $z,z'\in\acute{X}$ separated by a fixed distance (in the metric $|\omega|^2$) from the conical points of $X$.

\subsection{Variation of eigenvalues} 
Let $\lambda_1\le \lambda_2\le \dots$ be the sequence of eigenvalues of the Szeg\"o extension $\Delta=\Delta_S$ of the spinor Laplacian on $(X,\omega,C)$ counted with multiplicities. Let $u_n$ be the corresponding normalized eigensections. The eigenvalues and eigensections of the (Szeg\"o) Laplacian $\Delta'$ on $(X',\omega',C')$ are denoted by $\lambda'_n$, $u'_n$, respectively. For a bounded domain $U\subset\mathbb{C}$, denote 
$$\Sigma_U=\Sigma_U(A_j,B_j,C_k):=\sum_{\lambda_n\in U}\lambda_n, \qquad \Sigma'_U:=\sum_{\lambda'_n\in U}\lambda'_n.$$ 
\begin{prop}
\label{variation of eigenvalues prop}
If $\partial U$ do not intersect the spectrum of $\Delta$, then $\Sigma_U$ is differentiable with respect to $\Re\nu$ and $\Im\nu$ {\rm(}$\nu=A_i,B_i,C_k${\rm)}, and
\begin{equation}
\label{variation of eigenvalues}
\begin{split}
2i\frac{\partial \Sigma_U}{\partial\nu}=\sum_{\lambda_n\in U}\int_{\nu^\dag}\omega(\varsigma)^{-1}\Big(\partial_\varsigma\big(h(\varsigma)u_k(\varsigma)\big)\partial_\varsigma\overline{u_{k}(\varsigma)}d\varsigma-|u_k(\varsigma)|^{2}\rho^{-2}(\varsigma)h(\varsigma)\lambda_n d\overline{\varsigma}\Big).
\end{split}
\end{equation}
\end{prop}
In the remaining part of this subsection, we prove Proposition \ref{variation of eigenvalues prop}.

\

\paragraph{\it Non-concentration of eigensections near conical points.} For $\epsilon>|\delta\nu|$, denote by $\breve{X}(\epsilon,\delta\nu)$ the domain in $X$ obtained from $\acute{X}$ by removing $\epsilon$-neighborhoods (in the metric $|\omega|^2$) of all conical points. Then $\breve{X}(\epsilon,\delta\nu)$ can be also considered as a domain in $X'$. Denote $\tilde{X}(\epsilon,\delta\nu)=X\backslash \breve{X}(\epsilon,\delta\nu)$ and $\tilde{X}'(\epsilon,\delta\nu)=X'\backslash \breve{X}(\epsilon,\delta\nu)$. Let us prove the estimate
\begin{equation}
\label{non-concentration of eigenfunctions}
\|u'_{n}\|^2_{L_2(\tilde{X}'(\epsilon,\delta\nu);C')}\le c_U\epsilon \qquad (\epsilon>|\delta\nu|, \lambda_n\in U).
\end{equation}
Denote by $\tilde{X}'_k(\epsilon)$ the $\epsilon$-neighborhood (in the metric $|\omega'|^2$) of $P'_k$. Let $\epsilon_0>0$ be a sufficiently small number. Since $(\Delta')^l u'_n=(\lambda'_n)^l u'_n$ and $\|u'_n\|_{L_2(X';C')}=1$, Proposition \ref{prop elliptic asymptotics} (with $X$ replaced by $\tilde{X}'_k(\epsilon_0)$) provides the asymptotics
$$u'_n(x'_k)=c'_{n,k,0,+}+c'_{n,k,1,+}x'_k+c'_{n,k,2,+}x_k'^2+(-2/3\pi)c'_{n,k,-1,-}\overline{x'_k}|x'_k|+O(x_k'^3)$$
which is uniform in $(X',\omega')$. Here, the coordinate $x'_k\in \tilde{X}'_k(\epsilon_0)$ is given by (\ref{dist coord}) with $P_k$, $\omega$ replaced by $P'_k$, $\omega'$, respectively. Note that the above asymptotics contains only the terms that allow $u'_n$ to be an element of ${\rm Dom}\Delta'={\rm Dom}\Delta'_S$ (cf. with (\ref{Szego domain})). In particular, $u'_n(x'_k)=O(1)$ on $\tilde{X}'_k(\epsilon_0)$ uniformly in $(X',\omega')$. Since $\rho'(z'_k)=h'(z'_k)=|z'_k|^{-1}$ in coordinates $x'_k$, we have 
$$\|u'_{n}\|^2_{L_2(\tilde{X}'_k(\epsilon))}\le c_U\epsilon^3.$$
If $\nu=C_k$, then $\tilde{X}'(\epsilon,\delta\nu)$ is contained in the union of $\tilde{X}'_s(\epsilon)$ ($s=1,\dots,2g-2$) and the last inequality implies (\ref{non-concentration of eigenfunctions}). 

To prove (\ref{non-concentration of eigenfunctions}) for the case $\nu=A_j$ (or $\nu=B_j$), it remains to estimate $\|u'_{n}\|^2_{L_2(X'\backslash\acute{X};C')}$. Note that $X'\backslash\acute{X}$ is contained in $\epsilon$-neighborhood (in the metric $|\omega'|^2$) of the path $\check{b_j}\subset\overline{\acute{X}}$ (or $\check{a}_j\subset\overline{\acute{X}}$). In view of (\ref{H l smooth incr}) and (\ref{Morrey}), the equations $(\Delta')^l u'_n=(\lambda'_n)^l u'_n$ imply that $|u'_n(z')|$ is bounded on $X'\backslash\acute{X}$ uniformly in $(X',\omega')$. Here the coordinate $z'$ is given by (\ref{z coordinate}) with $\omega'$ instead of $\omega$. Since the area of $X'\backslash\acute{X}$ (in the metric $|\omega'|^2$) is $O(|\delta\nu|)=O(\epsilon)$, we have $\|u'_{n}\|^2_{L_2(X'\backslash\acute{X};C')}\le c_U\epsilon$. Thus, (\ref{non-concentration of eigenfunctions}) is proved for the case $\nu=A_j,B_j$.

Note that formula (\ref{non-concentration of eigenfunctions}) remains valid with omitted primes. 

\

\paragraph{\it Connection between eigenvalues and resolvent kernel.} Since 
\begin{equation}
\label{resolven kernel via eigensections}
\mathscr{R}'_\lambda(z,z'):=\sum_{n}\frac{u'_n(z)\overline{u'_n(z')}}{\lambda'_n-\lambda},
\end{equation}
we have
\begin{equation}
\label{resolvent-eigenvalues connection pointwise}
-\frac{1}{2\pi i}\int_{\partial U}\mathscr{R}'_\lambda(z,z')\mathscr{Q}(\lambda)d\lambda=\sum_{\lambda'_n\in U}u'_n(z)\overline{u'_n(z')}\mathscr{Q}(\lambda'_n)
\end{equation}
for any function $\mathscr{Q}$ holomorphic in the bounded domain $U\subset\mathbb{C}$ and continuous up to $\partial U$. The same formulas with omitted primes are valid for $\mathscr{R}_\lambda$.

\

\paragraph{\it Continuity of $\nu\mapsto\lambda_n(\nu)$.} Denote $\mathfrak{C}=\sharp\{\lambda_k\in U\}$ and $\mathfrak{C}'=\sharp\{\lambda'_k\in U\}$. In formula (\ref{resolvent-eigenvalues connection pointwise}), put $z'=z$ and $\mathscr{Q}=1$, then integrate over $\breve{X}(\epsilon,\delta\nu)$. In view of (\ref{non-concentration of eigenfunctions}), we obtain
\begin{align*}
-\frac{1}{2\pi i}\int\limits_{\breve{X}(\epsilon,\delta\nu)}\Big(\int\limits_{\partial U}\mathscr{R}'_\lambda(z,z')&d\lambda\Big)\Big|_{z'=z}h(z)dS(z)=\\=&\sum_{\lambda'_{n}\in U}\|u'_{n}\|^2_{L_2(\breve{X}(\epsilon,\delta\nu);C)}=\mathfrak{C}'(1+O(\epsilon)).
\end{align*}
This formula remains valid with omitted primes. Therefore,
\begin{align*}
\mathfrak{C}'-\mathfrak{C}=-\frac{1}{2\pi i}\int\limits_{\lambda\in\partial U}\int\limits_{\breve{X}(\epsilon,\delta\nu)}u_\lambda(z,z)h(z)dS(z)d&\lambda+O\big((\mathfrak{C}'+\mathfrak{C})\epsilon\big)=\\
=&O(|\delta \nu|)_{\epsilon}+O\big((\mathfrak{C}'+\mathfrak{C})\epsilon\big),
\end{align*}
where the index $\epsilon$ in $O(|\delta \nu|)_{\epsilon}$ means that the estimate is not uniform in $\epsilon$ (but still uniform in $(X',\omega')$). Since the left-hand side is integer, we have $\mathfrak{C}'=\mathfrak{C}$ for sufficiently small $|\delta \nu|$. Since $U$ is arbitrary, the eigenvalues $\lambda_n=\lambda_n(\nu)$ are continuous functions of $\nu=A_j,B_j,C_k$. 

\

\paragraph{\it Differentiability of $\Sigma_U$ with respect to $\nu$.} Without loss of generality, one can assume that $U$ contains eigenvalues $\lambda_n=\dots=\lambda_{n+m-1}$ of $\Delta$. Let $v_n,\dots,v_{n+m-1}$ be the corresponding eigenfunctions orthonormalized with respect to scalar product in $L_2(\breve{X}(\epsilon,\delta\nu);C')$ (for small $\epsilon$, the orthonormalization is possible due to (\ref{non-concentration of eigenfunctions})). Introduce the $m\times m$-matrices $\mathcal{M}'$ and $\mathcal{M}$ with entries
$$\mathcal{M}'_{pq}=(u'_{n+p},v_{n+q})_{L_2(\breve{X}(\epsilon,\delta\nu);C)}, \qquad \mathcal{M}_{pq}=(u_{n+p},v_{n+q})_{L_2(\breve{X}(\epsilon,\delta\nu);C)}.$$
Denote $\hat{\mathscr{Q}}={\rm diag}\{\mathscr{Q}(\lambda_{n}),\dots,\mathscr{Q}(\lambda_{n+m-1})\}$ and $\hat{\mathscr{Q}}'={\rm diag}\{\mathscr{Q}(\lambda'_{n}),\dots,\allowbreak \mathscr{Q}(\lambda'_{n+m-1})\}$.
Multiplying both parts of (\ref{resolvent-eigenvalues connection pointwise}) by $v_{n+p}(z')\overline{v}_{n+q}(z)$ and integrating, we obtain
\begin{align*}
\int\limits_{\breve{X}(\epsilon,\delta\nu)}\int\limits_{\breve{X}(\epsilon,\delta\nu)}\Big(\frac{-1}{2\pi i}\int\limits_{\partial U}\mathscr{R}'_\lambda(z,z')\mathscr{Q}(\lambda)d\lambda\Big)v_{n+p}(z')\overline{v}_{n+q}(z)h(z)dS(z&)h(z')dS(z')=\\
=&(\mathcal{M}^{'*}\hat{\mathscr{Q}}'\mathcal{M}')_{pq}.
\end{align*}
This equality remains valid with omitted primes. Then asymptotics (\ref{uniform differentablity}) implies
\begin{align}
\label{differentiability of eigensums}
\begin{split}
(\mathcal{M}^{'*}\hat{\mathscr{Q}}'\mathcal{M}')_{pq}-(\mathcal{M}^{*}\hat{\mathscr{Q}}\mathcal{M})_{pq}&=\\=\frac{-1}{2\pi i}\int\limits_{\lambda\in\partial U}\int\limits_{\breve{X}(\epsilon,\delta\nu)}\int\limits_{\breve{X}(\epsilon,\delta\nu)}u_\lambda(z,z')&\mathscr{Q}(\lambda)v_{n+p}(z')\overline{v}_{n+q}(z)h(z)dS(z)h(z')dS(z')d\lambda=\\
&=\mathfrak{D}(\nu | \mathscr{Q},\epsilon)_{pq}\delta\nu+\mathfrak{D}(\overline{\nu}|\mathscr{Q},\epsilon)_{pq}\overline{\delta\nu}+o(|\delta\nu|)_{\epsilon},
\end{split}
\end{align}
where
$\mathfrak{D}(\ \cdot\ |\mathscr{Q},\epsilon)$ is the $m\times m$-matrix with entries
\begin{align*}
\mathfrak{D}&(\ \cdot\ |\mathscr{Q},\epsilon)_{pq}:=\\
=&\frac{-1}{2\pi i}\int\limits_{\lambda\in\partial U}\int\limits_{\breve{X}(\epsilon,\delta\nu)}\int\limits_{\breve{X}(\epsilon,\delta\nu)}\frac{\partial\mathscr{R}_\lambda(z,z')}{\partial \ \cdot}\mathscr{Q}(\lambda)\overline{v}_{n+q}(z)v_{n+p}(z')h(z)dS(z)h(z')dS(z')d\lambda.
\end{align*}
Let $\mathcal{M}'=\mathcal{U}'\mathcal{A}'$ be a polar decomposition of $\mathcal{M}'$, where $\mathcal{U}'$ and $\mathcal{A}'$ are unitary and hermitian matrices, respectively. Let also $\mathcal{M}=\mathcal{U}\mathcal{A}$ be a polar decomposition of $\mathcal{M}$. Then (\ref{differentiability of eigensums}) with $\mathscr{Q}=1$ implies
$$\mathcal{A}'^2-\mathcal{A}^2=\mathcal{M}'^{*}\mathcal{M}'-\mathcal{M}^{*}\mathcal{M}=\mathfrak{D}(\nu | 1,\epsilon)\delta\nu+\mathfrak{D}(\overline{\nu}|1,\epsilon)\overline{\delta\nu}+o(|\delta\nu|)_{\epsilon}.$$
Thus, 
$$\mathcal{A}'^{-1}-\mathcal{A}^{-1}=-\frac{1}{2}\big(\mathfrak{D}(\nu | 1,\epsilon)\delta\nu+\mathfrak{D}(\overline{\nu}|1,\epsilon)\overline{\delta\nu}\big)+o(|\delta\nu|)_{\epsilon}=O(|\delta\nu|)_{\epsilon}.$$
Let $\mathscr{Q}(\lambda)=\lambda-\lambda_n$, then $\hat{\mathscr{Q}}=0$ and $\hat{\mathscr{Q}}'=\Sigma'_U-\Sigma_U$. In addition, (\ref{differentiability of eigensums}) yields
\begin{align*}
\mathcal{U}'^{*}\hat{\mathscr{Q}}'&\mathcal{U}'=\mathcal{A}'^{-1}\mathcal{M}'^{*}\hat{\mathscr{Q}}'\mathcal{M}'\mathcal{A}'^{-1}=\\
=&\mathcal{A}'^{-1}\Big(\mathfrak{D}(\nu | \mathscr{Q},\epsilon)\delta\nu+\mathfrak{D}(\overline{\nu}|\mathscr{Q},\epsilon)\overline{\delta\nu}+o(|\delta\nu|)_{\epsilon}\Big)\mathcal{A}'^{-1}=\\
=&\mathcal{A}^{-1}\Big(\mathfrak{D}(\nu | \mathscr{Q},\epsilon)\delta\nu+\mathfrak{D}(\overline{\nu}|\mathscr{Q},\epsilon)\overline{\delta\nu}\Big)\mathcal{A}^{-1}+o(\delta\nu)_{\epsilon}.
\end{align*}
Taking the trace from both sides, we arrive at
\begin{align*}
\Sigma'_U-\Sigma_U={\rm Tr}\big[\mathcal{A}^{-1}\mathfrak{D}(\nu | \mathscr{Q},\epsilon)\mathcal{A}^{-1}\big]\delta\nu+{\rm Tr}\big[\mathcal{A}^{-1}\mathfrak{D}(\overline{\nu}|\mathscr{Q},\epsilon)\mathcal{A}^{-1}\big]\overline{\delta\nu}+o(\delta C_k)_{\epsilon}.
\end{align*}
As a corollary, $\Sigma_U$ is differentiable with respect to $\Re\nu$, $\Im\nu$, and 
\begin{equation}
\label{derivative of eig}
\partial\Sigma_U/\partial\nu={\rm Tr}\big[\mathcal{A}^{-1}\mathfrak{D}(\nu | \mathscr{Q},\epsilon)\mathcal{A}^{-1}\big].
\end{equation}
In particular, the right hand side is independent of $\epsilon$.

\

\paragraph{\it Calculation of $\partial\Sigma_U/\partial\nu$.} In view of (\ref{non-concentration of eigenfunctions}), the eigensections $u'_n,\dots,u'_{n+m-1}$ can be chosen in such a way that $\mathcal{M}-I=O(\epsilon_0^{1/2})$. Then passing to the limit as $\epsilon\to 0$ in (\ref{derivative of eig}) yields
\begin{align}
\label{derivative of eigenv}
\begin{split}
\frac{\partial \Sigma_U}{\partial\nu}=\sum_{p}\mathfrak{D}(\nu | &\mathscr{Q},0)_{pp}=\\=-\sum_{p}\int\limits_{X}\int\limits_{X}\underset{\lambda=\lambda_n}{\rm Res}\Big[(\lambda-\lambda_n)&\frac{\partial\mathscr{R}_\lambda(z,z')}{\partial\nu}\Big]\overline{u}_{n+p}(z)u_{n+p}(z')h(z)dS(z)h(z')dS(z')=\\
=-\int\limits_{X}\underset{\lambda=\lambda_n}{\rm Res}\Big[(\lambda-\lambda_n)&\frac{\partial\mathscr{R}_\lambda(z,z')}{\partial\nu}\Big]\Big|_{z'=z}h(z)dS(z).
\end{split}
\end{align}
In view of (\ref{resolven kernel via eigensections}), the equality
\begin{equation}
\label{diff operators to reso ker}
\begin{split}
\int\limits_{X}\underset{\lambda=\lambda_n}{\rm Res}\Big((\lambda_n-\lambda)a(\lambda)\mathcal{P}_{\varsigma}\mathscr{R}_\lambda(\varsigma,z)\cdot\mathcal{Q}_{\varsigma}\mathscr{R}_\lambda(z,\varsigma)\Big)h(z)dS(z)=\\
=-\sum_{\lambda_k=\lambda_n}\mathcal{P}_{\varsigma}u_k(\varsigma)\cdot\mathcal{Q}_{\varsigma}\overline{u_{k}(\varsigma)}a(\lambda_k)
\end{split}
\end{equation}
is valid for any differential operators $\mathcal{P}$, $\mathcal{Q}$ and functions $\lambda\mapsto a(\lambda)$ holomorphic near $\lambda_n$. Combining (\ref{derivative of eigenv}), (\ref{variational formulas}), and (\ref{diff operators to reso ker}), we arrive at (\ref{variation of eigenvalues}). Proposition \ref{variation of eigenvalues prop} is proved.

\subsection{Variation of zeta function}
In this section, we prove Theorem \ref{main theorem}.

\subsubsection{Variation of $\zeta(2|\Delta-\mu)$.}
Let $z,z'$ be coordinates (\ref{z coordinate}). Formally differentiating the series for $\zeta(s|\Delta)=\zeta(s|\Delta_s(X,\omega))$ with respect to $\nu=A_j,B_j,C_k$ and taking into account formula (\ref{variation of eigenvalues}) and the equalities 
\begin{align*}
\rho^{2}(z)\partial_{\lambda}^2\big(\lambda h(z)\mathscr{R}_\lambda(z,z')\big)&=\sum_{n}\frac{2\lambda_n u_n(z)\overline{u_{n}(z')}\rho^{2}(z)h(z)}{(\lambda_n-\lambda)^{3}}\\
\partial_{\lambda}^2\partial_z\partial_{z'}\big(h(z)\mathscr{R}_\lambda(z,z')\big)&=\sum_{n}\frac{2\partial_z(h(z)u_n(z))\partial_{z'}\overline{u_n(z')}}{(\lambda_n-\lambda)^{3}},
\end{align*}
we obtain
\begin{align}
\label{zeta function s=2}
\begin{split}
2i&\partial_\nu\zeta(2|\Delta-\mu)=-2\sum_{n}(\lambda_n-\mu)^{-3}\Big(2i\frac{\partial\lambda_n}{\partial \nu}\Big)=\\
=&\int\limits_{\nu^\dag}\omega(z)^{-1}\Big[\big[\rho^{2}(z)\partial_{\mu}^2\big(\mu h(z)\mathscr{R}_\mu(z,z')\big)\big]\Big|_{z'=z}d\overline{z}-\big[\partial_{\mu}^2\partial_z\partial_{z'}\big(h(z)\mathscr{R}_\mu(z,z')\big)\big]\Big|_{z'=z}dz\Big].
\end{split}
\end{align}
To justify the above calculations, we should show that the derivatives of the resolvent kernel in the right hand side of (\ref{zeta function s=2}) are well defined on the diagonal $z'=z$. 

\ 

\paragraph{\it Near-diagonal asymptotics of resolvent kernel.}  In view of (\ref{elipticity}), the resolvent kernel $\mathscr{R}_\lambda$ of $\Delta=\Delta_S$ admits representation
\begin{equation}
\label{repres res ker near diag}
\mathscr{R}_\lambda(z,z')=\mathscr{R}_{q,\lambda}(\cdot,z')+[(\Delta_S-\lambda)^{-1}(\Delta-\lambda)\mathscr{R}_{q,\lambda}(\cdot,z')](z,z') \qquad (q=0,1,\infty)
\end{equation}
outside conical points, where
\begin{align*}
\mathscr{R}_{0,\lambda}(z,z')&=-\frac{2}{\pi}{\rm log}|z'-z|\cdot \chi_{z'}(z),\\
\mathscr{R}_{1,\lambda}(z,z')&=-\frac{2}{\pi}{\rm log}|z'-z|\cdot (1-\lambda|z'-z|^2)\chi_{z'}(z),\\
\mathscr{R}_{\infty,\lambda}(z,z')&=\frac{2}{\pi}K_0(-2i|z'-z|\sqrt{\lambda})\chi_{z'}(z).
\end{align*}
Here, $\chi_{z'}$ is a cut-off function equal to one near $z'$ and to zero, the support of $\chi_{z'}$ is sufficiently small, $K_0$ is the Macdonald function. In view of (\ref{C l smoth incr}), the second term in the right hand side of (\ref{repres res ker near diag}) is $(1+2q)$-differentiable with respect to $z,z'$ at the diagonal $z=z'$. Then the functions 
\begin{equation}
\label{Phi and Psi terms}
\begin{split}
\Phi_\lambda(z):&=\lambda \big[\mathscr{R}_\lambda(z,z')-\mathscr{R}_{0,\lambda}(z,z')\big]\Big|_{z'=z},\\
\Psi_\lambda(z):&=\partial_z\partial_{z'}\big[\mathscr{R}_\lambda(z,z')-\mathscr{R}_{1,\lambda}(z,z')\big]\Big|_{z'=z}.
\end{split}
\end{equation}
are well-defined and holomorphic with respect to $\lambda$ outside the spectrum of $\Delta$. Since $\partial_{\lambda}^2\big(\lambda\mathscr{R}_{0,\lambda}(z,z'))=0$ and $\partial_{\lambda}^2\partial_z\partial_{z'}\mathscr{R}_{1,\lambda}(z,z')=0$, and $h(z)=\rho(z)=1$, one can rewrite (\ref{zeta function s=2}) as 
\begin{equation}
\label{zeta function s=2 Phi Psi}
2i\partial_\nu\zeta(2|\Delta-\mu)=\partial_{\mu}^2\int_{z\in\nu^\dag}\big(\Phi_\mu(z)d\overline{z}-\Psi_\mu(z)dz\big).
\end{equation}

\

\paragraph{\it Asymptotics of $\Phi_\lambda$ and $\Psi_\lambda$ as $\Re\lambda\to -\infty$.} In view of (\ref{elipticity}), we have 
$$(\Delta-\lambda)\mathscr{R}_{\infty,\lambda}(\cdot,z')=[\Delta,\chi_z']\frac{2}{\pi}K_0(-2i|z'-z|\sqrt{\lambda}).$$
Since the support of $z\mapsto\chi_z'(z)$ is separated from the diagonal $z=z'$, and $K_0(\kappa)$ decays exponentially as $\Re \kappa\to +\infty$, the right hand side and all its derivatives decay exponentially as $\Re\lambda\to -\infty$. Hence, in view of (\ref{C l smoth incr}), the function 
$$z,z'\mapsto [(\Delta_S-\lambda)^{-1}(\Delta-\lambda)\mathscr{R}_{\infty,\lambda}(\cdot,z')](z,z')=\mathscr{R}_\lambda(z,z')-\mathscr{R}_{\infty,\lambda}(z,z')$$
and all its derivatives decay exponentially as $\Re\lambda\to -\infty$ uniformly with respect to $z,z'$ separated from conical points. As a corollary, we have
\begin{align}
\label{Psi near infinity}
\begin{split}
\Psi_\lambda(z)&=\partial_z\partial_{z'}\big[\mathscr{R}_\lambda(z,z')-\mathscr{R}_{\infty,\lambda}(z,z')+\mathscr{R}_{\infty,\lambda}(z,z')-\mathscr{R}_{1,\lambda}(z,z')\big]\Big|_{z'=z}=\\
=O(e^{\epsilon\Re\lambda}&)+\partial_z\partial_{z'}\Big[\frac{2}{\pi}\Big(K_0(-2i|z'-z|\sqrt{\lambda})+(1-\lambda|z'-z|^2){\rm log}|z'-z|\Big]\Big|_{z'=z}=\\
&=O(e^{\epsilon\Re\lambda})+0
\end{split}
\end{align}
uniformly with respect to $z$ separated from conical points, where $\epsilon>0$. Similarly, we have
\begin{align}
\label{Phi near infinity}
\begin{split}
\Phi_\lambda(z):=\lambda \big[\mathscr{R}_\lambda(z,z')-\mathscr{R}_{\infty,\lambda}(z,z')+\mathscr{R}_{\infty,\lambda}(z,z')-\mathscr{R}_{0,\lambda}(z,z')\big]\Big|_{z'=z}=O(e^{\epsilon\Re\lambda})+\\
+\frac{2\lambda}{\pi}\Big(K_0(-2i|z'-z|\sqrt{\lambda})+{\rm log}|z'-z|\Big)\Big|_{z'=z}=-\frac{\lambda}{\pi}{\rm log}\lambda-\lambda\Big(i+\frac{2\gamma}{\pi}\Big)+O(e^{\epsilon\Re\lambda}).
\end{split}
\end{align}

\

\paragraph{\it Calculation of $\Phi_0$ and $\Psi_0$.} Since $\lambda=0$ is not an eigenvalue of $\Delta=\Delta_S$, formulas (\ref{repres res ker near diag}) and (\ref{Phi and Psi terms}) yield
\begin{equation}
\label{Phi of zero}
\Phi_0=0.
\end{equation}  
Next, since $\Delta=\Delta_S$, formula (\ref{Green to Szego}) is valid, and
\begin{equation}
\label{Psi Szego}
\Psi_0(z)=\partial_z\partial_{z'}\Big(\mathcal{G}(z,z')+\frac{2}{\pi}{\rm log}|z'-z|\Big)\Big|_{z'=z}=-\frac{1}{\pi}\partial_z\Big(\mathcal{S}(z,z')-\frac{1}{z'-z}\Big)\Big|_{z'=z}.
\end{equation}
The asymptotics
\begin{equation}
\label{Szeco kernel asymptotics}
\begin{split}
S(z,z')-\frac{1}{z'-z}=\sum_{i=1}^{g}\partial_{\xi_i}{\rm log}\theta\left[^p _q\right](\xi=0)\upsilon_i(z)&+\\
+\left[\frac{S_0(z)}{6}+\sum_{i=1}^{g}\partial_{\xi_i}{\rm log}\theta\left[^p _q\right](\xi=0)\partial_{z}\upsilon_i(z)+\right.&\\
+\left.\sum_{i,j=1}^{g}\frac{\partial^2_{\xi_i \xi_j}\theta\left[^p _q\right](\xi=0)}{\theta\left[^p _q\right](0)}v_i(z) v_j(z)\right]&\frac{z'-z}{2}+O\big((z'-z)^2\big), \quad z'\to z.
\end{split}
\end{equation} 
is valid (see p.29, \cite{Fay}), where $\upsilon_1,\dots,\upsilon_g$ is the basis in the space of Abelian differentials on $X$ dual to the canonical homology basis. The Bergman projective connection $S_0$ transforms as $$S_0(x')=S_0(x)\Big(\frac{\partial x}{\partial x'}\Big)^2+\{x,x'\}$$
under a holomorphic change of variables $x\mapsto x'$. Due to the chain rule 
$$\{f,x'\}=\{f,x\}\Big(\frac{\partial x}{\partial x'}\Big)^2+\{x,x'\}$$ 
for Schwarzian derivative, the formula
\begin{equation}
\label{Bergman quadratic differential}
\mathcal{W}(x)=S_0(x)-\Big\{\int_{x}^{x_0}\omega,x\Big\}
\end{equation}
defines a quadratic differential such that $\mathcal{W}(z)=S_0(z)$. In view of (\ref{Psi Szego}), (\ref{Szeco kernel asymptotics}), (\ref{Bergman quadratic differential}), and (\ref{tau function}), we have 
\begin{equation}
\label{Psi of zero}
\Psi_0=\frac{\mathcal{W}}{12\pi\omega}+\frac{1}{2\pi}\Big[\sum_{i=1}^{g}\partial_{\xi_i}{\rm log}\theta\left[^p _q\right](\xi)\partial\Big(\frac{\upsilon_i}{\omega}\Big)+\sum_{i,j=1}^{g}\frac{\partial^2_{\xi_i \xi_j}\theta\left[^p _q\right](\xi)}{\theta\left[^p _q\right](\xi)}\frac{v_i v_j}{\omega}\Big]\Big|_{\xi=0}.
\end{equation}

\subsubsection{Variation of $\zeta(s|\Delta)$}

\

\paragraph{\it Connection between $\zeta(s|\Delta_S)$ and $\zeta(2|\Delta_S-\mu)$.} From (\ref{zeta function via spec part func}), the formula $K(t|A-\mu)=e^{\mu t}K(t|A)$, and the inversion formula for the Laplace transform it follows that
\begin{equation}
\label{spec part func via zeta func}
K(t|A):=\frac{\Gamma(s)}{2\pi i}\int\limits_{\mu_0-i\infty}^{\mu_0+i\infty}e^{-\mu t}t^{1-s}\zeta(s|A-\mu)d\mu.
\end{equation}
This formula is valid for $A=\Delta_F$ if $\Re s\ge 2$ and $\mu_0\in(0,\lambda_1^F)$. In view of Theorem \ref{prop comparing determinants F S}, it is also valid for $A=\Delta_S$ if $\Re s\ge 2$ and $\mu_0\in(0,\lambda_1^S)$. Put $s=2$ and substitute (\ref{spec part func via zeta func}) into (\ref{zeta function via spec part func}). As a result, we obtain 
\begin{equation}
\label{zeta s via zeta 2}
\begin{split}
(s-1)\zeta(s|A-\lambda)=\frac{s-1}{2\pi i\Gamma(s)}\int\limits_{\mu_0-i\infty}^{\mu_0+i\infty}\Big(\int\limits_{0}^{+\infty}e^{(\lambda-\mu) t}t^{s-2} dt\Big)\zeta(2|A-\mu)d\mu=\\
=\frac{1}{2\pi i}\int\limits_{\gamma}(\mu-\lambda)^{1-s}\zeta(2|A-\mu)d\mu,
\end{split}
\end{equation}
where $A=\Delta_F,\Delta_S$, $\Re s> 1$ and $\gamma$ is the union of the paths (\ref{paths}). For $\Re s\le 1$, formula (\ref{zeta s via zeta 2}) is also valid but its left hand side and right hand side are understood as analytic continuations of them from the half-plane $\Re s>1$. 

Let us make the analytic continuation of (\ref{zeta s via zeta 2}) more explicit by excluding the singular terms. Suppose that $A=\Delta_F$. Rewrite asymptotics (\ref{spec partition asymp}) as
\begin{equation}
\label{spec part function asymptotics 0+infty}
K(t|A)=e^{-\mu_0 t}\Big[{\rm Area}(X,|\omega|^2)\frac{1+\mu_0 t}{\pi t}+\frac{3g-1}{8}\Big]+\tilde{K}(t|A)
\end{equation}
where $\tilde{K}(t|A)=O(t)$ as $t\to +0$ and $\tilde{K}(t|A)=O(e^{-\mu_0 t/2})$ as $t\to -\infty$. The substitution of (\ref{spec part function asymptotics 0+infty}) into (\ref{zeta function via spec part func}) yields
\begin{equation}
\label{regularized zeta 1}
\zeta(s|A-\lambda)=\zeta_0(s,\lambda)+\tilde{\zeta}(s|A-\lambda),
\end{equation}
Here
\begin{equation}
\label{regularized zeta 1 singular}
\begin{split}
\zeta_0(s,\lambda):=\frac{{\rm Area}(X,|\omega|^2)}{\pi(\mu_0-\lambda)^{s-1}(s-1)}\Big(1+\frac{\mu_0(s-1)}{\mu_0-\lambda}&\Big)+\frac{3g-1}{8(\mu_0-\lambda)^{s}}=\\
=\frac{1}{2\pi i(s-1)}&\int\limits_{\gamma}(\mu-\lambda)^{1-s}\zeta_0(2,\mu)d\mu.
\end{split}
\end{equation}
The remainder
\begin{equation}
\label{regularized zeta 2}
\tilde{\zeta}(s|A-\lambda)=\frac{1}{\Gamma(s)}\int_{0}^{+\infty}t^{s-1}e^{\lambda t}\tilde{K}(t|A)dt=\frac{1}{2\pi i(s-1)}\int\limits_{\gamma}(\mu-\lambda)^{1-s}\tilde{\zeta}(2|A-\mu)d\mu
\end{equation} 
is well-defined and holomorphic in $s,\lambda$ for $\Re s>-1$, $\Re\lambda<\mu_0$. In view of the Main Theorem, \cite{Mooers} (still valid for $1/2$-forms), the same formulas are valid for $A=\Delta_S$.

\

\paragraph{\it Variation of ${\rm det}\Delta_S$.} Let $\Re s>1$. Differentiating (\ref{zeta s via zeta 2}) with respect to $\nu$, taking into account (\ref{zeta function s=2 Phi Psi}), and changing the order of integrations, we obtain
\begin{equation*}
2i(s-1)\partial_\nu\zeta(s|\Delta_S)=\int_{z\in\nu^\dag}\big(\mathrm{J}^{\Phi}_{s}(z)d\overline{z}-\mathrm{J}^{\Psi}_{s}(z)dz\big),
\end{equation*}
where
$$\mathrm{J}^{\Phi}_{s}(z):=\frac{1}{2\pi i}\int_{\gamma}\mu^{1-s}\partial_{\mu}^2\Phi_\mu(z)d\mu, \qquad \mathrm{J}^{\Psi}_{s}(z):=\frac{1}{2\pi i}\int_{\gamma}\mu^{1-s}\partial_{\mu}^2\Psi_\mu(z)d\mu.$$
In view of (\ref{Psi near infinity}), the function $\mathrm{J}^{\Psi}_{s}(z)$ is analytic in $s\in \mathbb{C}$ and the representation 
$$\mathrm{J}^{\Psi}_{s}=\pi^{-1}e^{-i\pi s}{\rm sin}(\pi s)\mathrm{J}^{\Psi,-,\epsilon}_{s}+\mathrm{J}^{\Psi,\circ,\epsilon}_{s},$$
holds, where
$$\mathrm{J}^{\Psi,-,\epsilon}_{s}(z):=\int_{-\infty}^{-\epsilon}\mu^{1-s}\partial_{\mu}^2\Psi_\mu(z)d\mu, \qquad \mathrm{J}^{\Psi,\circ,\epsilon}_{s}(z):=\frac{1}{2\pi i}\int_{\gamma_\circ}\mu^{1-s}\partial_{\mu}^2\Psi_\mu(z)d\mu$$
are analytic in $s\in\mathbb{C}$, and 
$$\mathrm{J}^{\Psi,-,0}_{0}=\big(\mu\partial_{\mu}\Psi_\mu-\Psi_\mu\big)\Big|_{\mu=-\infty}^{\mu=0}=-\Psi_0 \qquad \mathrm{J}^{\Psi,\circ,0}_{0}=\partial_{s}\mathrm{J}^{\Psi,\circ,0}_{s=0}=0$$
due to (\ref{Psi near infinity}) and (\ref{Psi of zero}). Hence,
\begin{equation}
\label{J Psi}
\mathrm{J}^{\Psi}_{0}=0, \qquad -\partial_s\mathrm{J}^{\Psi}_s\big|_{s=0}=\Psi_0.
\end{equation}

We have
$$\mathrm{J}^{\Phi}_{s}=\pi^{-1}e^{-i\pi s}{\rm sin}(\pi s)\mathrm{J}^{\Phi,-,\epsilon}_{s}+\mathrm{J}^{\Phi,\circ,\epsilon}_{s},$$
where
$$\mathrm{J}^{\Phi,-,\epsilon}_{s}(z):=\int_{-\infty}^{-\epsilon}\mu^{1-s}\partial_{\mu}^2\Phi_\mu(z)d\mu, \qquad \mathrm{J}^{\Phi,\circ,\epsilon}_{s}(z):=\frac{1}{2\pi i}\int_{\gamma_\circ}\mu^{1-s}\partial_{\mu}^2\Phi_\mu(z)d\mu.$$
Since the integration contour $\gamma_\circ$ is compact, the function $\mathrm{J}^{\Phi,\circ,\epsilon}_{s}$ is holomorphic in $s\in\mathbb{C}$. Also, since $\partial_{\mu}^2\Phi_\mu$ is holomorphic in $\mu$ outside the spectrum of $\Delta_S$ (hence, near zero), we have $\mathrm{J}^{\Phi,\circ,\epsilon}_{0}(z)=0$. Similarly, since $\partial_{\mu}^2\Phi_\mu,\partial_{\mu}\Phi_\mu,\Phi_\mu$ are holomorphic in $\mu$ near zero, and $\Phi_0=0$ due to (\ref{Phi of zero}), integration by parts yields
\begin{align*}
-\partial_s\mathrm{J}^{\Phi,\circ,\epsilon}_{s=0}=\frac{1}{2\pi i}\int_{\gamma_\circ}\mu{\rm log}\mu\partial_{\mu}^2\Phi_\mu d\mu=\frac{1}{2\pi i}\big[\mu{\rm log}\mu\partial_{\mu}\Phi_\mu-\Phi_\mu-{\rm log}\mu\Phi_\mu\big]\Big|_{\epsilon e^{-i(\pi-0)}}^{\epsilon e^{+i(\pi-0)}}+\\+\int_{\gamma_\circ}\mu^{-1}\Phi_\mu d\mu=-\epsilon\partial_{\mu}\Phi_{\mu=-\epsilon}-\Phi_{-\epsilon}.
\end{align*}
Let us represent $\mathrm{J}^{\Phi,-,\epsilon}_{s}$ as
\begin{equation}
\label{rpresentatin log asymp at inf}
\mathrm{J}^{\Phi,-,\epsilon}_{s}(z)=-\frac{1}{\pi}\int_{-\infty}^{-\epsilon}\mu^{-s}d\mu+\mathrm{J}^{\tilde{\Phi},-,\epsilon}_{s}(z),
\end{equation}
where
$$\mathrm{J}^{\tilde{\Phi},-,\epsilon}_{s}(z):=\int_{-\infty}^{-\epsilon}\mu^{1-s}\partial_{\mu}^2\tilde{\Phi}_\mu(z)d\mu, \qquad \tilde{\Phi}_\mu(\zeta):=\Phi_\mu(\zeta)+\pi^{-1}\mu\big({\rm log}\mu+\pi i+2\gamma\big).$$
The integral in the right hand side of (\ref{rpresentatin log asymp at inf}) is equal to $(-\epsilon)^{1-s}/(1-s)$ for $\Re s>1$, and admits the analytic continuation to the whole complex plane except the pole at $s=1$. In view of (\ref{Phi near infinity}), the integral $\mathrm{J}^{\tilde{\Phi},-,\epsilon}_{s}(z)$ is analytic in $s\in\mathbb{C}$. Integrating by parts and taking into account (\ref{Phi near infinity}), we obtain 
\begin{align*}
\mathrm{J}^{\tilde{\Phi},-,\epsilon}_{0}(z):=\Big(\mu\partial_{\mu}\tilde{\Phi}_\mu-\tilde{\Phi}_\mu\Big)\Big|_{\mu=-\infty}^{\mu=-\epsilon}=-\epsilon\partial_{\mu}\tilde{\Phi}_{\mu=-\epsilon}-\tilde{\Phi}_{-\epsilon}.
\end{align*}
Due to the above formulas, we have
\begin{align}
\label{J phy}
\begin{split}
\mathrm{J}^{\Phi}_{0}&=\mathrm{J}^{\Phi,\circ,\epsilon}_{0}=0, \\
\partial_s\mathrm{J}^{\Phi}_{s=0}&=\mathrm{J}^{\Phi,-,\epsilon}_{s=0}+\epsilon\partial_{\mu}\Phi_{\mu=-\epsilon}+\Phi_{-\epsilon}=\\=&\epsilon\partial_{\mu}\Phi_{\mu=-\epsilon}+\Phi_{-\epsilon}+\frac{\epsilon}{\pi}-\epsilon\partial_{\mu}\tilde{\Phi}_{\mu=-\epsilon}-\tilde{\Phi}_{-\epsilon}=0.
\end{split}
\end{align}
In view of (\ref{J Psi}) and (\ref{J phy}), formula (\ref{zeta s via zeta 2}) with $s=0$ yields 
$$\partial_\nu\zeta(0|\Delta_S)=0.$$
Now, differentiate (\ref{zeta s via zeta 2}) with respect to $s$, put $s=0$, and take into account (\ref{J Psi}) and (\ref{J phy}). As a result, we obtain
\begin{align}
\label{variation of the determinant}
\begin{split}
2i\partial_\nu{\rm log \ det}\Delta_S&=-\partial_\nu\partial_s\zeta(s|\Delta_S)|_{s=0}=\\
=\int\limits_{\zeta\in\nu^\dag}&\big(\partial_s\mathrm{J}^{\Phi}_{s=0}(\zeta)d\overline{\zeta}-\partial_s\mathrm{J}^{\Psi}_{s=0}(\zeta)d\zeta\big)=\int\limits_{\zeta\in\nu^\dag}\Psi_0(\zeta)d\zeta.
\end{split}
\end{align}
It remains to check that one can interchange the analytic continuation to the neighborhood of $s=0$ and the differentiation with respect to $\nu$ in the above calculations. First, note that no divergent integrals (\ref{rpresentatin log asymp at inf}) appear in the above calculations if one replaces $\zeta(s|\Delta_S)$ with $\tilde{\zeta}(s|\Delta_S)$. Indeed, the formal differentiation of both sides of (\ref{regularized zeta 2}), one obtains
\begin{equation}
\label{tilde zeta repeating}
\begin{split}
2i(s-1)[\partial_\nu\zeta(s|\Delta_S)-\partial_\nu\zeta_0(s,0)]=2i(s-1)\partial_\nu\tilde{\zeta}(s|\Delta_S)=\\
=\frac{1}{\pi}\int\limits_{\gamma}\mu^{1-s}(\partial_\nu\zeta(2|\Delta_S-\mu)-\partial_\nu\zeta_0(2,\mu))d\mu.
\end{split}
\end{equation}
In view of (\ref{regularized zeta 1 singular}) and the well-known formula 
\begin{equation}
\label{area formula}
{\rm Area\,}(X, |\omega|^2)=-\Im\sum_{k=1}^{g}A_k\overline{B_k},
\end{equation}
we have 
$$\partial_\nu\zeta_0(2,\mu)=\frac{1}{\pi\mu}\int_{\nu^\dag}\omega+O(|\mu|^{-2}) \qquad (\Re\mu\to -\infty).$$
Thus, the divergent parts of the integrals $-(1/\pi)\int_{\gamma}\mu^{1-s}\partial_\nu\zeta_0(2,\mu))d\mu$ and (\ref{rpresentatin log asymp at inf}) cancel each other in the right-hand side of (\ref{tilde zeta repeating}). So, the right-hand side of (\ref{tilde zeta repeating}) is well-defined for $\Re s>-1$. Since $\tilde{\zeta}(s|\Delta_S)$ is also well-defined for $\Re s>-1$, it is differentiable with respect to $\nu$ and obeys (\ref{tilde zeta repeating}) for $\Re s>-1$. So, the only term in (\ref{zeta s via zeta 2}) which needs the analytic continuation is $2i(s-1)\partial_\nu\zeta_0(s,0)=(1/\pi)\int_\gamma \mu^{1-s}\partial_\nu\zeta_0(2,\mu)d\mu$. Due to (\ref{regularized zeta 1 singular}) and (\ref{area formula}), the differentiation of this term commutes with its analytic continuation.

\

\paragraph{\it Completing the proof of {\rm(\ref{determinant expression})}.} Substituting (\ref{Psi of zero}) into (\ref{variation of the determinant}) and taking into account that the integral of $\partial(\upsilon_i/\omega)=d(\upsilon_i/\omega)$ over any closed cycle equals zero, one gets
\begin{equation}
\label{det log var 1}
\partial_\nu{\rm log \ det}\Delta_S=\frac{1}{2\pi i}\int_{\nu^\dag}\Big[\frac{\mathcal{W}}{12\omega}+\sum_{i,j=1}^{g}\frac{\partial^2_{\xi_i \xi_j}\theta\left[^p _q\right](\xi)}{\theta\left[^p _q\right](\xi)}\Big|_{\xi=0}\frac{v_i v_j}{2\omega}\Big].
\end{equation}
According to the definition of the Bergman tau function $\tau$ on $\mathcal{H}_g(1,\dots,1)$ (see (3.1), \cite{KokKorot}), we have 
\begin{equation}
\label{tau function}
\frac{\partial {\rm log}\tau}{\partial\nu}=\frac{-1}{12\pi i}\oint_{\nu^\dag}\frac{\mathcal{W}}{\omega}, \qquad \frac{\partial\tau}{\partial\overline{\nu}}=0.
\end{equation}
Recall that the theta-function satisfies the heat equation 
\begin{equation}
\label{theta function heat equation}
\frac{\partial \theta\left[^p _q\right](\xi|\mathbb{B})}{\partial \mathbb{B}_{ij}}=\frac{1}{4\pi i}\frac{\partial^2 \theta\left[^p _q\right](\xi|\mathbb{B})}{\partial\xi_i\partial\xi_j}.
\end{equation}
In addition, the dependence of $\mathbb{B}$ on coordinates (\ref{moduli space coordinates}) is described by the following formula (see (2,28) and the remark after (2.81), \cite{KokKorot})
\begin{equation}
\label{variation of perod matrix}
\frac{\partial \mathbb{B}_{ij}}{\partial\nu}=\oint_{\nu^\dag}\frac{\upsilon_i \upsilon_j}{\omega}, \qquad \frac{\partial \mathbb{B}_{ij}}{\partial\overline{\nu}}=0.
\end{equation}
In view of (\ref{tau function}), (\ref{theta function heat equation}), and (\ref{variation of perod matrix}), equality (\ref{det log var 1}) takes the form
\begin{align*}
\partial_\nu\Big[{\rm log \ det}&\Delta_S+\frac{1}{2}{\rm log}\tau]\Big]=\frac{1}{4\pi i\theta\left[^p _q\right](0|\mathbb{B})}\sum_{i,j=1}^{g}\partial^2_{\xi_i \xi_j}\theta\left[^p _q\right](0|B)\frac{\partial \mathbb{B}_{ij}}{\partial\nu}=\\
=&\frac{1}{\theta\left[^p _q\right](0|\mathbb{B})}\frac{\partial \theta\left[^p _q\right](0|\mathbb{B})}{\partial\nu}=\frac{\partial {\rm log}\theta\left[^p _q\right](0|\mathbb{B})}{\partial\nu}=\frac{1}{2}\partial_\nu{\rm log}|\theta\left[^p _q\right](0|\mathbb{B})|^2.
\end{align*}
The last formula can be rewritten as 
$$\partial_\nu\big(\text{left hand side of (\ref{determinant expression})}-\text{right hand side of (\ref{determinant expression})}\big)=0.$$
Since both sides of (\ref{determinant expression}) are real-valued, one can replace $\partial_\nu$ with $\partial_{\overline{\nu}}$ and, hence, with $d$ in the formula above. By this, Theorem \ref{main theorem} is proved.

\end{document}